\documentclass[11pt,a4paper]{article}
\usepackage[utf8]{inputenc}
\usepackage{amsmath, bm}
\usepackage{amsthm}
\usepackage{mathtools}
\usepackage{float}
\usepackage{multirow}
\usepackage{array}
\usepackage{amsfonts}
\usepackage{amssymb}
\usepackage{abstract}
\usepackage[numbers, sort&compress]{natbib}
\usepackage{graphicx}
\usepackage{authblk} 
\usepackage[labelfont=bf, font=footnotesize]{caption}
\usepackage{caption} 
\usepackage{hyperref}
\usepackage{multicol}
\usepackage{tikz}
\usepackage{wrapfig}
\usepackage{natbib} 
\setcitestyle{authoryear}
\usetikzlibrary{decorations.pathreplacing}
\usetikzlibrary{shapes.geometric}
\usetikzlibrary{backgrounds}
\usepackage{booktabs}
\usepackage{caption}
\usepackage[left=2cm,right=2cm,top=2cm,bottom=2cm]{geometry}
\usepackage{titlesec}
\numberwithin{equation}{section}
\newtheorem{theorem}{Theorem}[section]

 \let\realcitealp\citealp
\bgroup\catcode`\#=12\egroup \renewcommand*{\citealp}[1]{{ (\color{blue}\footnotesize\realcitealp{#1})}}
\graphicspath{%
	{myfig/}
}
\definecolor{blue(ryb)}{rgb}{0.01, 0.28, 1.0}
\definecolor{deepskyblue}{rgb}{0.0, 0.75, 1.0}
\definecolor{crimson}{rgb}{0.86, 0.08, 0.24}
\definecolor{cordovan}{rgb}{0.54, 0.25, 0.27}
\definecolor{daffodil}{rgb}{1.0, 1.0, 0.19}
\definecolor{green-yellow}{rgb}{0.68, 1.0, 0.18}
\tikzstyle{block} = [rectangle, 
    text width=10em, text centered, node distance=2cm,minimum height=1em] \label{tikz: block}
\tikzstyle{cloud} = [draw, ellipse, 
    minimum height=1em] \label{tikz: ellipse}
\tikzstyle{decision} = [diamond, text badly centered, aspect=2, inner sep=0.1pt] \label{tikz: diamond}
\theoremstyle{definition}
\newtheorem{PD}{Problem}[subsection]
\title{Reduced order methods for parametric optimal flow control in coronary bypass grafts, towards patient-specific data assimilation}

\author[1]{Zakia Zainib \thanks{Email: \href{mailto:zakia.zainib@sissa.it}{zakia.zainib@sissa.it}}}
\author[1]{Francesco Ballarin \thanks{Email: \href{mailto:francesco.ballarin@sissa.it}{francesco.ballarin@sissa.it}}}
\author[2]{Stephen Fremes \thanks{Email: \href{mailto:stephen.fremes@sunnybrook.ca}{stephen.fremes@sunnybrook.ca}}}
\author[3]{Piero Triverio \thanks{Email: \href{mailto:piero.triverio@utoronto.ca}{piero.triverio@utoronto.ca}}}
\author[2]{Laura Jim\'{e}nez-Juan \thanks{Email: \href{mailto:laura.jimenezjuan@sunnybrook.ca}{laura.jimenezjuan@sunnybrook.ca}}}
\author[1]{Gianluigi Rozza \thanks{Email: \href{mailto:gianluigi.rozza@sissa.it}{gianluigi.rozza@sissa.it}}}

\affil[1]{ mathLab, Mathematics area, SISSA - International School for Advance Studies, Trieste, Italy}
\affil[2]{ Sunnybrook Health Sciences Centre, Toronto, Canada}
\affil[3]{ Department of Electrical \& Computer Engineering,
Institute of Biomaterials \& Biomedical Engineering, University of Toronto, Toronto, Canada}

\date {}

\begin{document}
 \maketitle
\vspace*{-1cm}
\begin{abstract}

Coronary artery bypass grafts (CABG) surgery is an invasive procedure performed to circumvent partial or complete blood flow blockage in coronary artery disease (CAD). 
In this work, we apply a numerical optimal flow control model to patient-specific geometries of CABG, reconstructed from clinical images of real-life surgical cases, in parameterized settings.
The aim of these applications is to match known physiological data with numerical hemodynamics corresponding to different scenarios, arisen by tuning some parameters. Such applications are an initial step towards matching patient-specific physiological data in patient-specific vascular geometries as best as possible.

Two critical challenges that reportedly arise in such problems are, $\left( i \right)$. lack of robust quantification of meaningful boundary conditions required to match known data as best as possible and $\left( ii \right)$. high computational cost. In this work, we utilize unknown control variables in the optimal flow control problems to take care of the first challenge. Moreover, to address the second challenge, we propose a time-efficient and reliable computational environment for such parameterized problems by projecting them onto a low-dimensional solution manifold through proper orthogonal decomposition ({POD})--Galerkin.
\end{abstract}
\section{Introduction and Motivation} \label{sec: Introduction}

In the last few decades, several efforts have been made to bring patient-specific computational hemodynamics modelling closer to real-life clinical scenarios. This is generally done through many-query parameterized computational fluid dynamics (CFD) modeling \citealp{AuricchioEtAl2018, KabinejadianEtAl2010, AgoshkovEtAl2006JSC, LassilaEtAl2013} combined with non-invasive medical imaging techniques. In this work, we will propose a numerical framework to model patient-specific hemodynamics with special focus on addressing two commonly arisen challenges in this field, that are, lack of robust quantification of meaningful boundary conditions and high computational cost.

Accuracy of boundary conditions is vital to the accuracy of patient-specific numerical hemodynamics modeling. In literature, according to the complexity of the problem and its aim, the boundary conditions are implemented directly or through surrogate models. For instance, the inflow profiles for coronary arteries are calculated from the known aortic flow rate or from the patient-specific pressure or flow waveforms. Comparatively, the lack of availability of patient-specific outflow waveforms for coronary arteries makes it harder to prescribe outflow boundary conditions. Many studies have used restrictive conditions at the outlets such as zero or fixed pressure or traction-free boundary conditions \citealp{PolitisEtAl2007, BertolottiEtAl2000, KuEtAl2002, SchrauwenEtAl2016, LeuprechtEtAl2002}, however, implementation of such conditions affects the accuracy of numerical simulations. Furthermore, several works have used surrogate models \citealp{SankaranarayananEtAl2005, ClementelEtAl2006, SankaranEtAl2012, TaylorEtAl2013} at the outlets to accurately approximate the corresponding boundary conditions. In the implementation of these models, often some parameters are manually tuned until the user-desired accuracy is obtained. Cardiovascular modeling and simulation is still lacking, in the authors' opinion, of fast, robust, automated tuning procedures for boundary conditions, even though great effort on this has been devoted in recent years. 

In this work, we propose automated quantification of problem-specific boundary conditions through an optimal flow control paradigm. By definition, the optimal flow control problems involve minimization of an objective functional while solving for an unknown control, responsible for changes in fluid flow behavior \citealp{Gunzburger, HinzeEtAl, Troltzsch, Quarteroni, QuarteroniRozza2003, HouRavindran1999, Dede2007}. Recently, some attempts have been made at filling the gap between data assimilation and boundary conditions quantification in hemodynamics modeling via these problems \citealp{TiagoEtAl2017, Veneziani2018, KoltukluogluEtAl2018}. These attempts are based upon the idea that through the unknown control variables one can solve for problem-specific boundary conditions while matching known data in {\it least-square} sense through minimization of a defined objective functional. In this work, we make a similar attempt while focusing on hard-to-quantify outflow conditions. This pipeline is specifically beneficial with respect to matching the desired data as with recent advances in imaging, this data may well be patient-specific experimental measurements, acquired e.g. through catheterization or 4D flow MRI. The main goal of this work is to setup the optimal control framework in patient-specific configurations, and show a methodological validation of the proposed numerical approach on a patient-specific configuration. Thus, for what concerns target data, we will limit ourselves to simpler analytical expressions. Integration of the proposed framework with real-life clinical measurements is currently ongoing \citealp{Francesca}. 

Furthermore, in the works by Tiago et al.\citealp{TiagoEtAl2017}, Romorowski et al.\citealp{Veneziani2018} and Koltukluo{\u g}lu et al.\citealp{KoltukluogluEtAl2018}, the authors commonly claim to face the challenge of high computational cost, even for non-parameterized optimal flow control problems. The computing cost in these problems is high owing to the coupling in the optimality system resulting in a mathematical problem with increased complexity \citealp{NegriEtAlStokes, HinzeEtAl} and the dense mesh required for the accuracy of numerical methods. In this work, we consider the optimal flow control problem in many-query parameterized settings \citealp{BonertEtAl2002, ClementelEtAl2006, BoutsianisEtAl2004, PolitisEtAl2008}, to take into account different hemodynamics scenarios modeled by parameter tuning. Consequently, computations are carried out in a repetitive environment, thus adding to the required computing effort. We will use reduced order methods (ROMs) \citealp{HesthavenEtAl2015}, specifically proper orthogonal decomposition ({POD})--Galerkin, to lower the computational cost. The ROMs have been used in many applications, such as in shape optimization of cardiovascular configurations\citealp{BallarinEtAl2016, ManzoniEtAl2012, Rozza2006}, in inverse problems arising in hemodynamics modeling\citealp{LassilaEtAl2013, BallarinEtAl2017} and in fluid-structure interaction problems\citealp{BallarinRozza2016}. Furthermore, these techniques have recently been further extended to optimal flow control problems \citealp{ItoKunisch2008,ItoRavindran1998,ItoRavindran2001,Dede2012, NegriEtAlElliptic, NegriEtAlStokes, KarcherEtAl2017, BaderEtAl2016} with applications in oceanographic problems\citealp{StrazzulloEtAl}.

In this work, we will utilize high-fidelity numerical methods as building blocks to construct the reduced order spaces. We will refer to the full-order discretized problems as {\it truth problems} and to their solutions as {\it snapshots}. The resulting reduced order spaces will retain the reliability  and accuracy of full order spaces with comparatively lower dimension. Furthermore, the computational efficiency of this reduced order model will arise from the decomposition of the computational procedure into two separate stages, namely, an offline phase and an online phase. The phase decomposition is based upon an affine decomposition assumption, that is, smooth parameter-dependent and continuous parameter-independent components can be separated. Computational cost in the offline phase depends upon dimension of full-order solution manifold and is therefore high, however, this phase needs to be carried out only once. The online phase has a computational cost independent of the dimension of the full order solution manifold and is, therefore, very inexpensive. To consider different hemodynamics scenarios, only the online phase needs to be repeated for different parameter values and therefore, the computational cost is anticipated to be reduced from order of days to order of seconds.

At the end of this article, we will demonstrate the applications of the reduced order parametrized optimal flow control pipeline to realistic 3D geometries representing coronary artery bypass grafts (CABGs), constructed from medical images of real-life clinical cases. In these cases, the grafts can be internal thoracic artery, radial artery extracted from human arm and saphenuous vein extracted from human leg\citealp{GhistaEtAl2013, OwidaEtAl2012}, and are connected to stenosed coronary arteries through bypass graft surgery, an invasive procedure commonly performed to overcome the blood flow blockage resulting from atherosclerosis. The methodology presented in this work shall contribute to render patient-specific computational cardiovascular modeling one step closer to clinical practice.

This work is organized in the following way: in section 1, we have introduced the focus and aim of this research work, the motivation for opting for an optimal flow control pipeline and combining it with reduced order methods. Section 2 will summarize mathematical details of the applied pipeline and the steps required for its application. Results for numerical tests ran on patient-specific CABGs and corresponding computational details will be reported in section 3. We will conclude this article and discuss future perspectives in section 4.

\section{Methodology} \label{sec: methodology}
In this section, we will discuss methodological details concerning application of the reduced order parametrized optimal flow control problems to patient-specific cardiovascular configurations. The complete pipeline is applied essentially in three steps: (i) extraction of patient-specific cardiovascular configurations from clinically acquired images (ii) data assimilation in parametrized optimal flow control regime, (iii) incorporation of reduced order methods. The forthcoming discussion is divided in three parts according to the mentioned steps.
\subsection{Geometrical Reconstruction} \label{sec: geometrical reconstruction}

\begin{figure}
\begin{minipage}{7cm}
\centering
\includegraphics[width = 3.5cm, keepaspectratio]{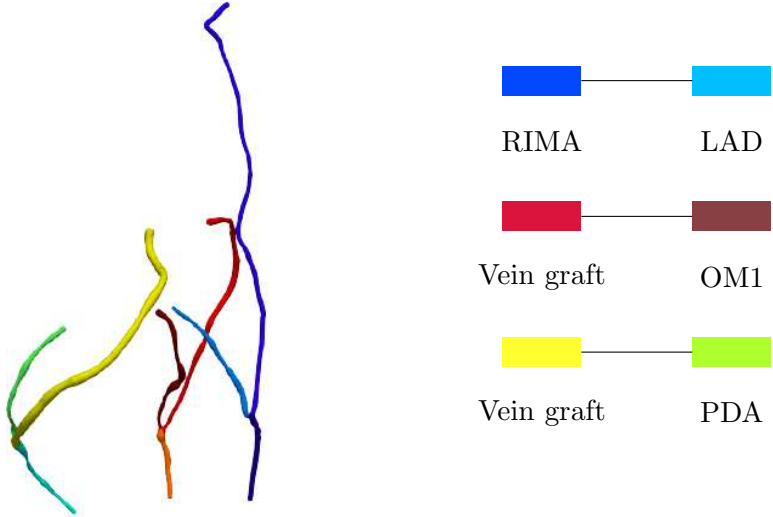}

\end{minipage}
\begin{minipage}{7cm}
\begin{tikzpicture}
\node[block, fill = blue(ryb), text width = 2em](node1){};
\node[block, fill = white, text width = 7em, below of = node1, node distance = 0.8cm](node2){RIMA};
\node[block, fill = deepskyblue, text width = 2em, right of = node1, node distance = 2.5cm](node3){};
\node[block, fill = white, text width = 5em, below of = node3, node distance = 0.8cm](node4){LAD};

\node[block, fill = crimson, text width = 2em, below of = node2, node distance = 1cm](node5){};
\node[block, fill = white, text width = 8em, below of = node5, node distance = 0.8cm](node6){Vein graft};
\node[block, fill = cordovan, text width = 2em, right of = node5, node distance = 2.5cm](node7){};
\node[block, fill = white, text width = 5em, below of = node7, node distance = 0.8cm](node8){OM1};
\node[block, fill = daffodil, text width = 2em, node distance = 1cm, below of = node6](node9){};
\node[block, fill = white, text width = 8em, below of = node9, node distance = 0.8cm](node10){Vein graft};
\node[block, fill = green-yellow, text width = 2em, right of = node9, node distance = 2.5cm](node11){};
\node[block, fill = white, text width = 5em, below of = node11, node distance = 0.8cm](node12){PDA};

\draw[-] (node1) -- (node3);
\draw[-] (node5) -- (node7);
\draw[-] (node9) -- (node11);

\end{tikzpicture}
\end{minipage}
\caption{Clinical case: triple coronary arteries bypass graft surgery}
\label{figure1}
\end{figure}
In this section, we will summarize the algorithm to construct patient-specific geometrical models of coronary artery bypass grafts from acquired medical image. For details of the reconstruction procedure, we refer the reader to Antiga et al. \citealp{AntigaEtAlChapter,VMTK1,VMTK2}  and Ballarin et al.\citealp{BallarinEtAl2016, BallarinEtAl2017}. This algorithm is employed in Python using the open-source 3-dimensional modeling libraries {\it Visualization Toolkit} (VTK)\citealp{VTK} and {\it Vascular Modeling Toolkit} (VMTK) \citealp{VMTK1}. Here, we will summarize it with respect to its application to the clinical case of a triple coronary artery bypass grafts surgery, performed at the {\it Sunnybrook Health Sciences Centre} in {\it Toronto, Canada}.

The data provided by the {\it Sunnybrook Health Sciences Centre} corresponding to this case, is a post-surgery computed tomography (CT) scan. The patient had at least three coronary arteries with occlusion ranging from mild to severe, thus, three different graft connections were made. These connections include right internal mammary artery (RIMA) grafted to bypass the blockage in left anterior descending artery (LAD) and two more connections made to first obtuse marginal artery (OM1) and posterior descending artery (PDA) separately using saphenuous veins (SV). Diseased coronary arteries and corresponding grafts, marked with different colors, are shown in figure \ref{figure1}.

\begin{figure}
\centering
\begin{tikzpicture}
\node[block, text width = 8em](node1){\includegraphics[width = \textwidth, keepaspectratio]{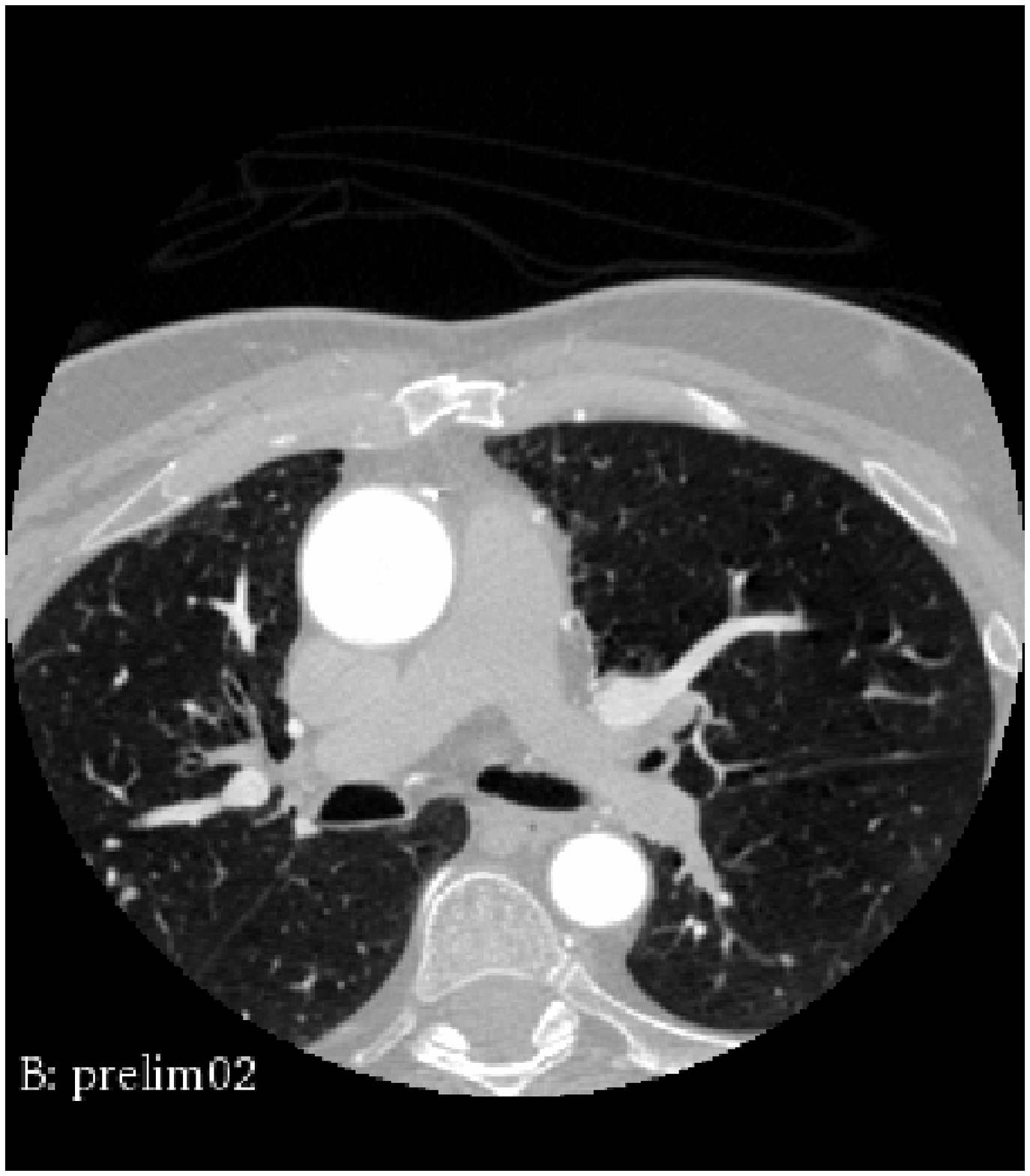}\\
{\footnotesize clinical image: CT scan}};
\node[block, right of = node1, node distance = 6cm, text width = 10em](node2){\includegraphics[width = \textwidth, keepaspectratio]{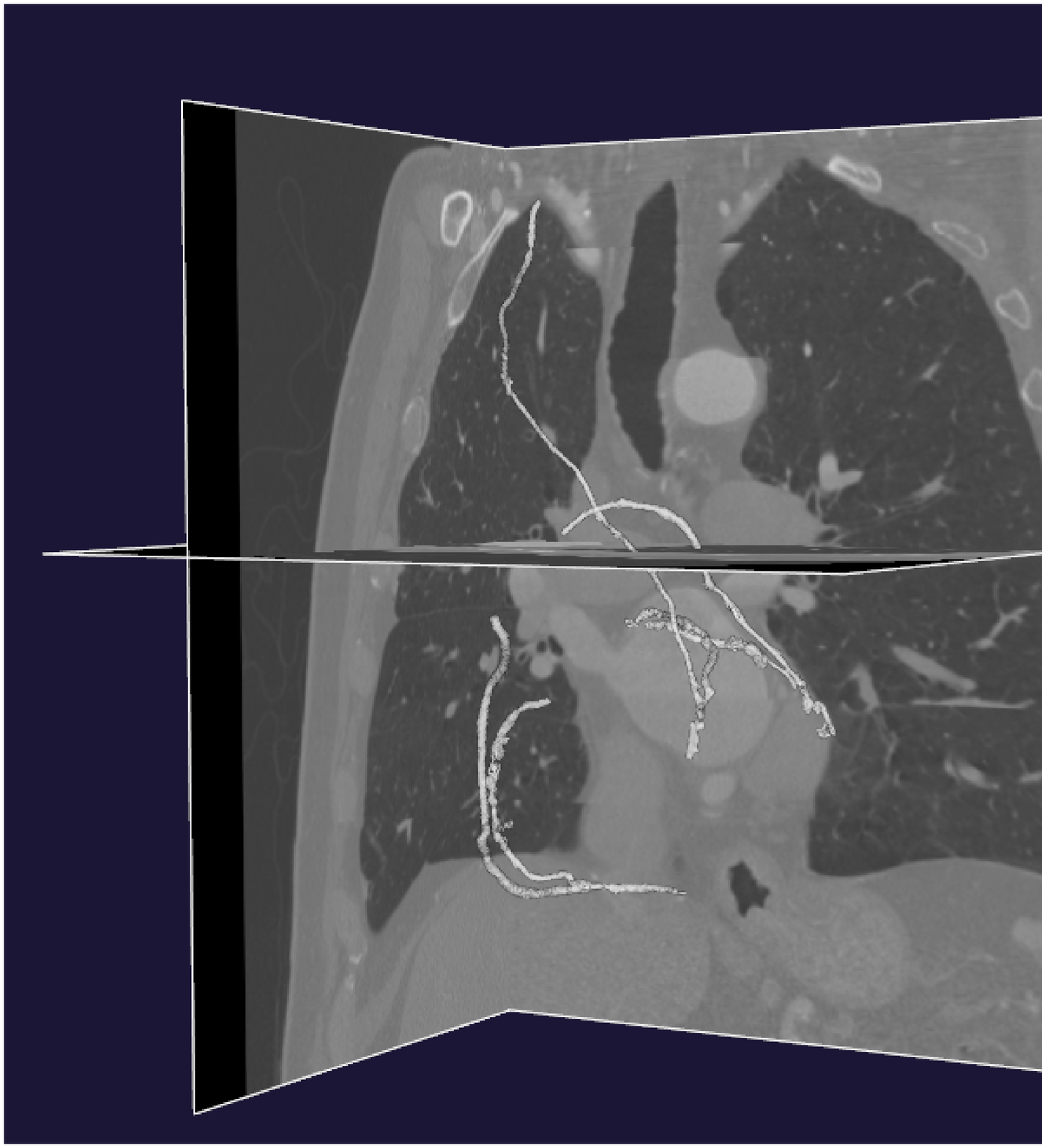}\\
{\footnotesize 3D surface reconstructed from clinical image}};
\node[block, right of = node2, node distance = 6cm, text width = 10em](node3){\includegraphics[width = 0.5\textwidth, keepaspectratio]{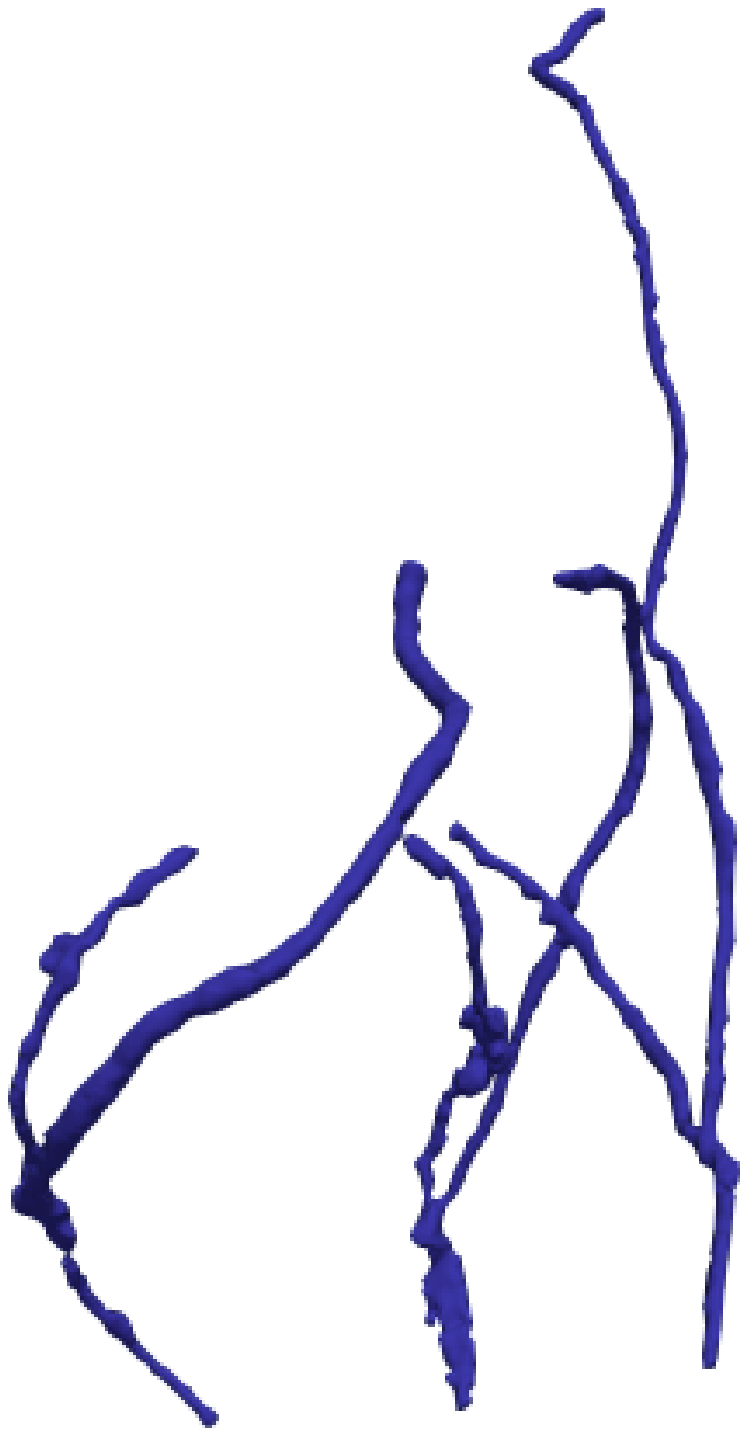}\\
{\footnotesize surface smoothed out using VMTK}};
\node[block,  below of = node3, node distance = 5.5cm, text width = 10em](node4){\includegraphics[width = 0.6\textwidth, keepaspectratio]{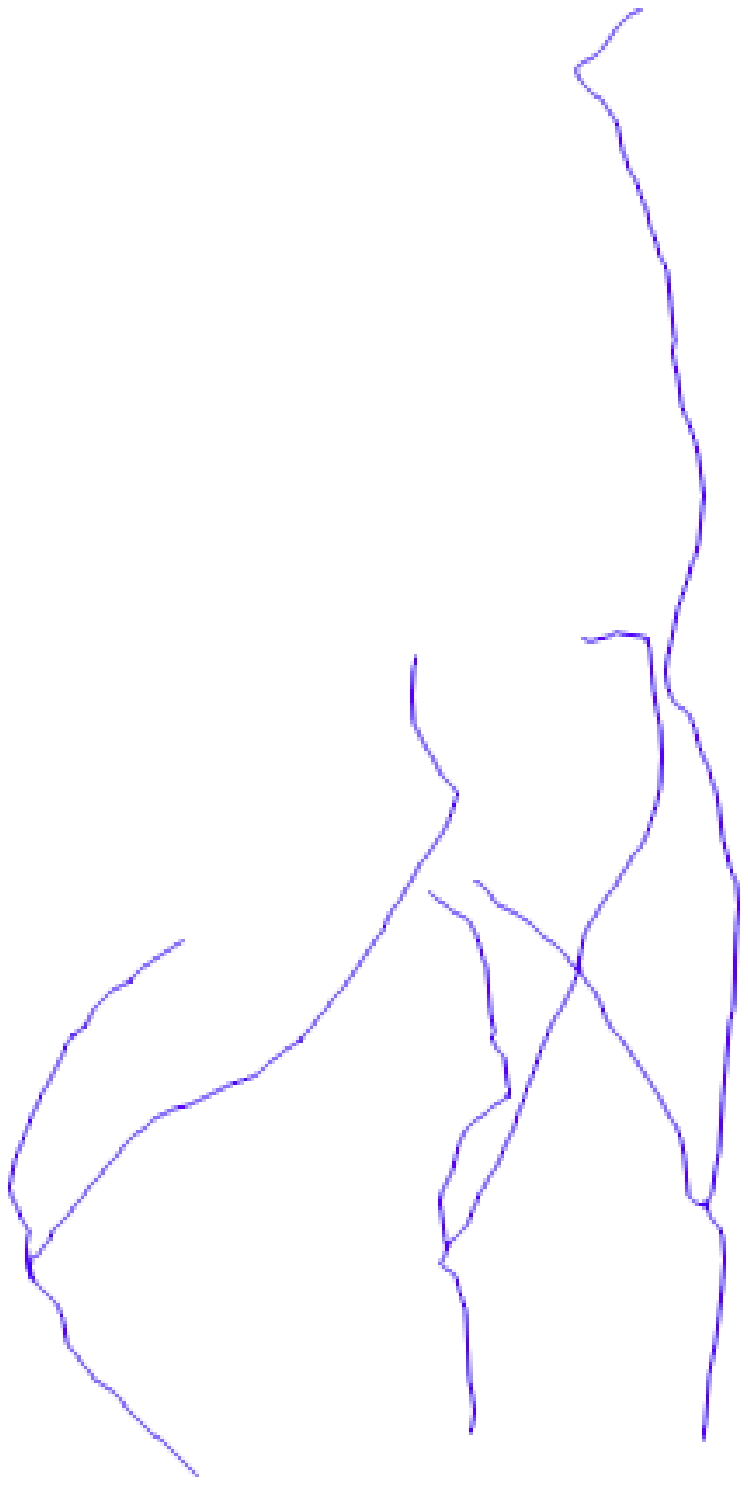}\\ {\footnotesize centerlines}};
\node[block, left of = node4, node distance = 6cm, text width = 10em](node5){\includegraphics[width = 0.7\textwidth, keepaspectratio]{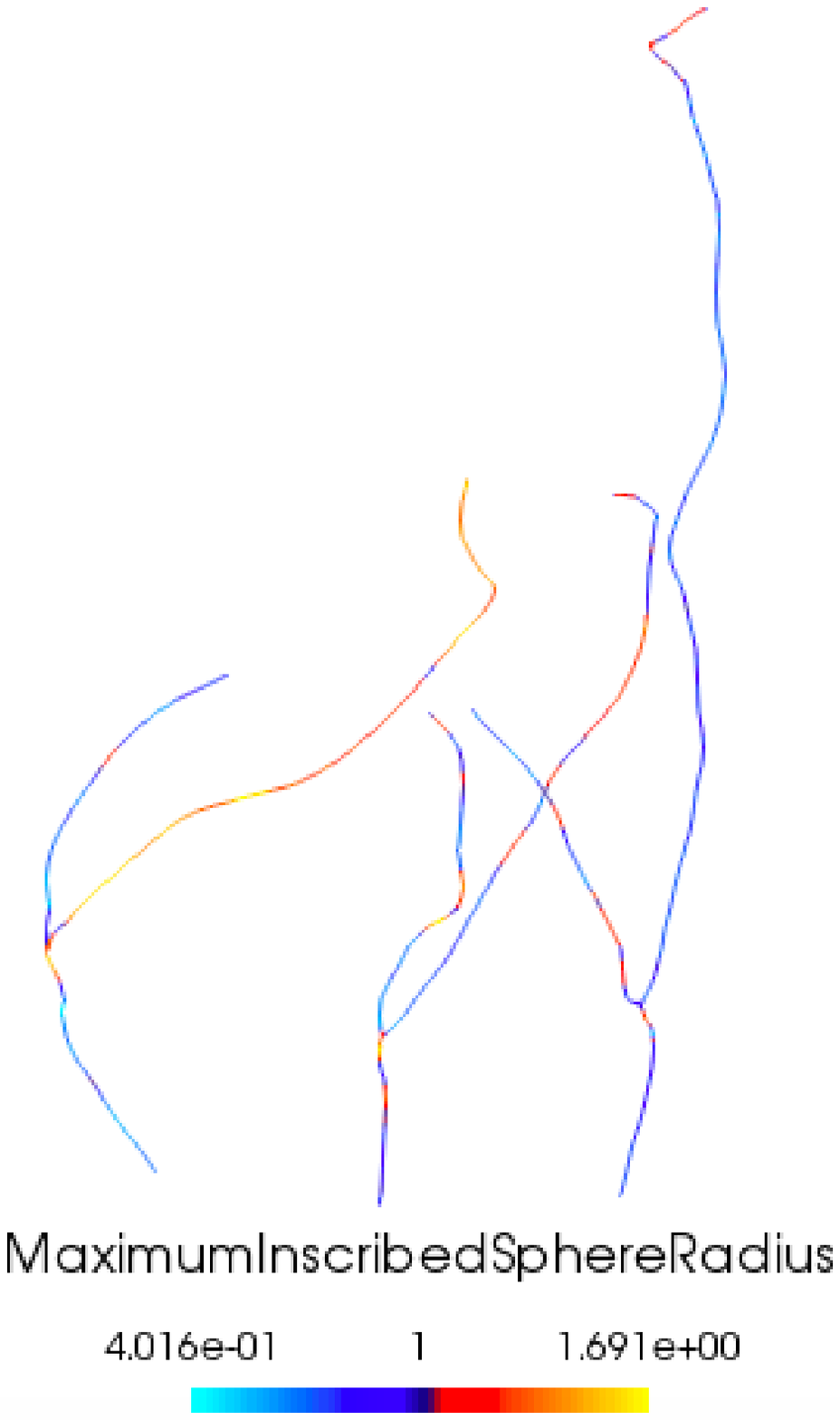}\\
{\footnotesize averagely smoothed centerlines}};
\node[block, left of = node5, node distance = 6cm, text width = 10em](node6){\includegraphics[width = 0.6\textwidth, keepaspectratio]{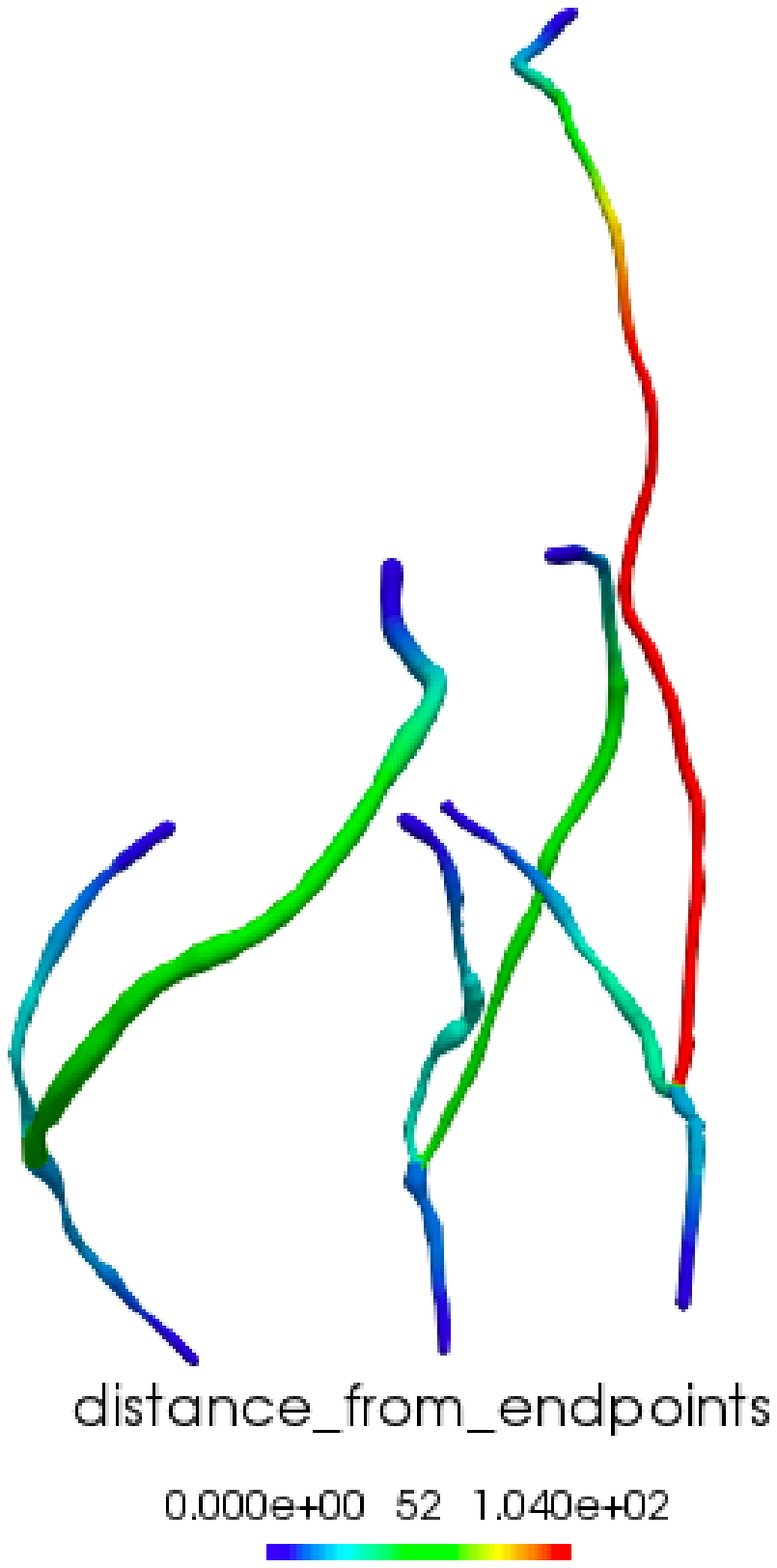}\\
{\footnotesize reconstructed geometry}};
\draw[->, red] (node1.east) [out=90, in=120] to (node2.west);
\draw[->, red] (node2.east) [out=90, in=120] to (node3.west);
\draw[->, red] (node3.east) [out=20, in=-20] to (node4.east);
\draw[->, red] (node4.west) [out=-120, in=-90] to (node5.east);
\draw[->, red] (node5.west) [out=-120, in=-90] to (node6.east);
\end{tikzpicture}
\caption{Algorithm applied for geometrical reconstruction from medical images}
\label{tikzfigure1}
\end{figure}

We pre-processed the obtained medical image in order to alleviate image segmentation and geometrical reconstruction process. The pre-processing step is done in three stages, that are resampling, smoothing and enhancement, we refer the reader to Antiga et al. \citealp{AntigaEtAlChapter} for details. The resampling stage aims at matching the resolution of the acquired image with the desired image segmentation process. Then, an anisotropic diffusion filter is used to reduce the high-frequency noise in the image and finally vessel enhancement filters are applied to enhance the visibility of vessel-shaped structures in comparison to other anatomical structures. Afterwards, using VMTK, we segment the pre-processed image into level sets, based on the colliding fronts approach. Then, a three-dimensional polygonal surface is generated through marching cube algorithm and is constructed by manually placing seed points according to visible intensity of the vessels. Due to the noise in clinical image and according to the visible intensity of coronary artery bypass grafts, the reconstruction process can add artificial artifacts corresponding to other parts of the cardiovascular system to the constructed geometries. We use VMTK smoothing filters to remove deformities to much extent, however the resulting tubular structures do not have sufficiently regular boundary, as is visible in the third figure, titled {\it surface smoothed out using VMTK}, in figure \ref{tikzfigure1}.

To generate sufficiently smooth surfaces preserving the same anatomical structures as reconstructed surfaces, we follow a centerlines-based approach \citealp{BallarinEtAl2016, BallarinEtAl2017}. We extract centerlines, that are the lines between two sections of lumen such that their minimal distance from the boundary is maximal \citealp{VMTK1, VMTK2, BallarinEtAl2017} using VMTK, which also outputs maximum inscribed spheres radii values associated to each point on the centerlines. After running average smoothing, we insert polyballs around the points on the centerlines according to the associated maximum inscribed spheres radii values. This yield a smooth 3D volume preserving patient-specific anatomical configuration \citealp{BallarinEtAl2016, BallarinEtAl2017}. The algorithm is summarized in figure \ref{tikzfigure1}.

Finally, for computational purposes we generate a tetrahedral mesh inside the reconstructed volumes and triangular mesh over the boundaries, using TetGen. To exploit Python based finite element libraries, for example, Dolfin, FEniCS and multiphenics \citealp{DOLFIN, FENICS2012, FENICS2015, multiphenics}, we read and write the mesh in {\it xml} format using VTK.

\subsection{Parametrized Optimal Flow Control Problems} \label{ss: Parametrized Optimal Flow Control Problems}
In this section, we will discuss the mathematical formulation and adopted solution strategy for parametrized optimal flow control problems. Structurally, these problems comprise of equations governing the physical state of the fluid flow, an objective functional to be minimized, parameters modeling physical and/or geometric characteristics and an unknown control that affects the fluid flow. Solving these problems, we aim at evaluating the unknown state and control variables that minimizes the objective functional. In this work, we consider PDEs-constrained optimization and utilize adjoint based techniques and a saddle-point framework to cast the problem in a monolithic structure.

\subsubsection{Problem formulation}
For the discussion henceforth, we will consider the patient-specific 3-dimensional CABGs geometrical models reconstructed from medical images following the algorithm discussed in section \ref{sec: geometrical reconstruction} as computational domain $\Omega$. We will refer to the boundaries as $\partial \Omega = \Gamma_{in} \cup \Gamma_{w} \cup \Gamma_{o}$, where $\Gamma_{in}$, $\Gamma_{o}$ and $\Gamma_{w}$ denote inlets of coronary arteries and bypass grafts, outlets of coronary arteries and vessel walls, respectively. Owing to the diameter of coronary arteries in our case, that ranges from $1\ mm$ to $2.2\ mm$, we can fairly assume that the size of blood particles such as red blood cells, white blood cells and platelets suspended in plasma, is smaller compared to the diameter of the vessel itself. Therefore, we consider the viscosity to be constant and the blood in coronary arteries and bypass grafts to behave as a Newtonian fluid. Thus, considering incompressible Navier-Stokes equations to model blood flow in $\Omega$ and physical parameters $\bm{\mu}\in \mathcal{D}\subset \mathbb{R}^d, d\in \mathbb{N}$, we define the state-constraints in strong form as below:
\begin{equation}
\label{NavierStokesEquations}
\begin{cases}
-\eta \Delta\bm{v}\left(\bm{\mu}\right) + \left(\bm{v}\left(\bm{\mu}\right)\cdot \nabla \right)\bm{v}\left(\bm{\mu}\right) + \nabla p \left(\bm{\mu}\right) = \bm{f}\left( \bm{\mu}\right),&\ \text{in}\ \Omega ,\\
\nabla\cdot\bm{v}\left( \bm{\mu}\right) = 0,&\ \text{in}\ \Omega,
\end{cases}
\end{equation}
where $\bm{v}$ and $p$ denote blood flow velocity and blood pressure respectively, $\bm{u}$ denotes control variables, $\eta$ is constant viscosity and $\bm{f}$ are body forces. Moreover, the first equation in \eqref{NavierStokesEquations} ensures conservation of momentum while second equation is the continuity equation for incompressible fluids. Furthermore, we assume $\Gamma_{w}$ to be rigid and non-permeable and consider the following conditions on $\partial\Omega$:
\begin{itemize}
\item[(i)] No-slip conditions at $\Gamma_{w}$.
\item[(ii)] Non-homogenuous Dirichlet inflow conditions at $\Gamma_{in}$.
\item[(iii)] Neumann outflow conditions at $\Gamma_{o}$.
\end{itemize}
These boundary conditions are mathematically written as:
\begin{equation}
\label{boundaryConditions}
\begin{cases}
\bm{v}\left(\bm{\mu}\right) = \bm{v}_{in}\left(\bm{\mu}\right),&\ \text{on}\ \Gamma_{in},\\
\bm{v}\left(\bm{\mu}\right) = \bm{0},&\ \text{on}\ \Gamma_{w},\\
-\eta\left(\nabla\bm{v}\left(\bm{\mu}\right)\right)\bm{n} + p\left(\bm{\mu}\right)\bm{n} = \bm{u}\left(\bm{\mu}\right), &\ \text{on}\ \Gamma_{o},
\end{cases}
\end{equation}
where $\bm{n}$ is outward unit normal at $\Gamma_o$. We set objective functional $\mathcal{J}$ to be a tracking-type functional, that is,
\begin{equation}
\label{costFunctional}
\mathcal{J}\left( \bm{v}, \bm{u};\bm{\mu}\right) = \frac{1}{2}\int_{\Omega} | \bm{v}\left(\bm{\mu}\right) - \bm{v_o}|^2 d\Omega+ \frac{\alpha}{2}\int_{\Gamma_{o}} |\bm{u}\left(\bm{\mu}\right)|^2 d\Gamma_{o},
\end{equation}
such that $\bm{v_o}$ is the patient-specific blood flow velocity provided by hospitals/clinics. $\mathcal{J}$ in equation \eqref{costFunctional}, thus, is the distance between measurements made through computational hemodynamics modeling and real-life patient-specific physiological data. In other words, we are interested in matching clinically-provided physiological data (for example, $\bm{v_o}$) with CFD variables (for example, $\bm{v}$) in $\mathcal{J}$. The second term on the right hand side of equation \eqref{costFunctional} is related to energy of control $\bm{u}$ with $\alpha > 0$ as penalization parameter. Furthermore, through control implementation in outlet boundary conditions \eqref{boundaryConditions}, we address the problem of automated quantification of meaningful boundary conditions required to match the desired data through simulation variables. The parametrized optimal flow control problem for this work reads:
\begin{PD}
\label{POFCP}
\it Given $\bm{\mu}\in \mathcal{D}$, find $\left( \bm{v}\left(\bm{\mu}\right), p\left(\bm{\mu}\right), \bm{u}\left(\bm{\mu}\right)\right)$ such that \eqref{costFunctional} is minimized while \eqref{NavierStokesEquations} and \eqref{boundaryConditions} are satisfied.
\end{PD}

We consider Hilbert spaces $V\left(\Omega\right)$ and $P\left(\Omega\right)$ for $\bm{v}$ and $p$ respectively and ease the notation by using $S = V\times P$ and $\bm{s} = \left( \bm{v}, p\right)$. For boundary control, we consider the Hilbert space $U\left(\Gamma_o\right)$, that is, $\bm{u}\in U\left(\Gamma_o\right)$\footnote{
\label{footnote1}
In case of control distributed across computational domain, $\bm{u}\in U\left(\Omega\right)$.
}. Furthermore, we consider the observation velocity $\bm{v_o}$ in another Hilbert space $Q\supseteq S$.

We adopt the adjoint-based Lagrangian approach\citealp{Quarteroni, Gunzburger, HinzeEtAl}, commonly used in PDEs-constrained optimal flow control problems, to solve problem \ref{POFCP}. First order optimality conditions ensure existence of a global minimum to the problem. These conditions will give rise to a coupled optimality system which can be solved through so-called {\it one-shot} or {\it all-at-once} approaches, discussed, for example, by Gunzburger\citealp{Gunzburger} and Stoll et al.\citealp{StollEtAl2013}. We desire to obtain the optimality system in the form of a monolithic-structured algebraic system in order to ease the computations for coupled problems with increased complexity, such as optimal control problems constrained by Stokes and Navier-Stokes equations. This is achieved when the optimality system is derived from the saddle-point form of parameterized optimal flow control problems\citealp{NegriEtAlElliptic, NegriEtAlStokes}. To cast problem \ref{POFCP} in saddle-point form, we rewrite cost functional $\mathcal{J}\left(\bm{v}, \bm{u}; \bm{\mu}\right)$ as:
\begin{equation}
\label{bilinearFormsCostFunctional}
\mathcal{J}\left(\bm{v}, \bm{u}; \bm{\mu}\right) = \frac{1}{2}\mathrm{m}\left( \bm{v}\left( \bm{\mu}\right) - \bm{v_o}, \bm{v}\left( \bm{\mu}\right) - \bm{v_o}\right) + \frac{\alpha}{2} \mathrm{n}\left(\bm{u}\left( \bm{\mu}\right), \bm{u}\left( \bm{\mu}\right)\right),
\end{equation}
such that $\mathrm{m}\left( \bm{v}\left( \bm{\mu}\right) - \bm{v_o}, \bm{v}\left( \bm{\mu}\right) - \bm{v_o}\right) =  \|\bm{v}\left( \bm{\mu}\right) - \bm{v_o}\|^2_{V\left(\Omega\right)}$ is a symmetric, continuous and non-negative form while $\mathrm{n}\left(\bm{u}\left( \bm{\mu}\right), \bm{u}\left( \bm{\mu}\right)\right)= \|\bm{u}\left( \bm{\mu}\right)\|^2_{U\left(\Gamma_{o}\right)}$ is symmetric, coercive and bounded bilinear form associated with energy of control and the norms are considered in the respective Hilbert solution spaces. Moreover, we consider Hilbert adjoint spaces $Z = Z_V\times Z_P$ such that $\bm{z} = \left(\bm{w}, q\right)$ is the ordered pair of adjoint velocity and adjoint pressure and we take $Z\equiv S$. State-constraints \eqref{NavierStokesEquations} and boundary conditions \eqref{boundaryConditions} can be written in weak formulation as:
\begin{equation}
\label{bilinearFormsNavierStokesConstraints}
\begin{cases}
\mathrm{a}\left( \bm{v}, \bm{w}; \bm{\mu} \right) + \mathrm{e}\left(\bm{v}, \bm{v}, \bm{w}; \bm{\mu}\right) + \mathrm{b}\left( p, \bm{w}; \bm{\mu} \right) + \mathrm{c}\left(\bm{u}, \bm{w}; \bm{\mu}\right) = \left\langle \bm{f}\left(\bm{\mu}\right), \bm{w}\right\rangle, \qquad \forall\ \bm{w}\in V,\\
\mathrm{b}\left( q, \bm{v}; \bm{\mu}\right) = 0, \qquad \forall\ q\in P.
\end{cases}
\end{equation}
Here $\mathrm{a} : V\times Z_V\rightarrow \mathbb{R}$ and $\mathrm{b} : V\times Z_P\rightarrow\mathbb{R}$ are the bilinear terms associated to diffusion and divergence operators respectively and are defined as:
\begin{equation*}
\mathrm{a}\left( \bm{v}, \bm{w}; \bm{\mu} \right) = \eta\int_{\Omega} \nabla\bm{v}\left(\bm{\mu}\right)\cdot \bm{w}\ d\Omega, \qquad \qquad \mathrm{b}\left( q, \bm{v}; \bm{\mu} \right) = -\int_{\Omega} q\left(\nabla\cdot \bm{v}\left(\bm{\mu}\right)\right)\ d\Omega
\end{equation*}
Moreover, $\mathrm{e}: V\times V\times Z_V\rightarrow \mathbb{R}$ is the non-linear convection term and $\mathrm{c} : U\times Z_V\rightarrow\mathbb{R}$ is the bilinear term associated to the control implemented through Neumann boundary conditions. These terms are defined below:
\begin{equation*}
\mathrm{e}\left( \bm{v}, \bm{v}, \bm{w}; \bm{\mu} \right) = \int_{\Omega}\left(\bm{v}\left(\bm{\mu}\right)\cdot\nabla\right)\bm{v}\left(\bm{\mu}\right)\cdot\bm{w}\ d\Omega, \qquad \qquad \mathrm{c}\left(\bm{u}, \bm{w}; \bm{\mu}\right) = -\int_{\Gamma_{o}}\bm{u}\left(\bm{\mu}\right) \cdot \bm{w}\ d\Gamma_{o}.
\end{equation*}

\noindent We introduce $X = S\times U$ such that for $\mathbf{x} = \left( \bm{v}, p, \bm{u}\right)$ and $\mathbf{y} = \left(\mathbf{y}_{\bm{v}}, \mathbf{y}_{p}, \mathbf{y}_{\bm{u}}\right)$, state constraints \eqref{bilinearFormsNavierStokesConstraints} are re-written as:
\begin{equation}
\label{stateConstraintsInGeneralForm}
\mathcal{B}\left( \mathbf{x}, \bm{z}; \bm{\mu}\right)+ \mathrm{e}\left( \bm{v}, \bm{v}, \bm{w}; \bm{\mu} \right) - \left\langle \bm{f}\left(\bm{\mu}\right),  \bm{z}\right\rangle = 0 , \qquad \forall\ \bm{z}\in Z.
\end{equation}
Here $\mathcal{B} : X\times Z\rightarrow \mathbb{R}$ is the operator associated to the linear part of state constraints. The saddle-point form of the parametrized optimal flow control problem, then, reads:
\begin{PD}
\label{SaddlepointProblem}
\it Given $\bm{\mu}\in \mathcal{D}$, find saddle-points $\left( \mathbf{x}\left(\bm{\mu}\right), \bm{z}\left(\bm{\mu}\right)\right)\in X\times Z$ of the following Lagrangian:
\begin{equation}
\label{SaddlepointForm}
\mathcal{L}\left( \mathbf{x}, \bm{z} ;\bm{\mu}\right) = J\left( \mathbf{x}; \bm{\mu}\right) + \mathcal{B}\left( \mathbf{x}, \bm{z}; \bm{\mu}\right)  + \mathrm{e}\left( \bm{v}, \bm{v}, \bm{w}; \bm{\mu} \right) - \left\langle \bm{f}\left(\bm{\mu}\right),  \bm{z}\right\rangle,
\end{equation}
where, 
\begin{equation*}
J\left( \mathbf{x}; \bm{\mu}\right) = \frac{1}{2}\mathcal{A}\left( \mathbf{x}, \mathbf{x} ; \bm{\mu}\right) - \left\langle\mathcal{H}\left(\bm{\mu}\right), \mathbf{x}\right\rangle,
\end{equation*}
such that $\mathcal{A}\left( \mathbf{x}, \mathbf{y}; \bm{\mu}\right) = \mathrm{m}\left(\bm{s}\left(\bm{\mu}\right), \bm{\xi}\right) + \alpha \mathrm{n}\left(\bm{u}\left(\bm{\mu}\right), \mathbf{y}_{\bm{u}}\right)$ and $\left\langle\mathcal{H}\left(\bm{\mu}\right), \mathbf{y}\right\rangle = \mathrm{m}\left(\bm{\xi}, \bm{s_o}\right),\ \bm{\xi} = \left(\mathbf{y}_{\bm{v}}, \mathbf{y}_{p}\right)$
\footnote{
\label{footnote3}
Note that $\mathcal{J}\left(\bm{v}, \bm{u}; \bm{\mu}\right) = J\left( \mathbf{x}; \bm{\mu}\right) + \frac{1}{2} \mathrm{m}\left(\bm{v_o}\left(\bm{\mu}\right), \bm{v_o}\left(\bm{\mu}\right)\right)$, where $\mathrm{m}\left(\bm{v_o}\left(\bm{\mu}\right), \bm{v_o}\left(\bm{\mu}\right)\right)$ being a constant term, has no impact on formulation and therefore, can be ignored.
}.
\end{PD}

Note that $\mathrm{e}\left( \bm{v}, \bm{v}, \bm{w}; \bm{\mu} \right) = 0$ leads to an optimization problem constrained by Stokes equations. In such case, a unique solution to the saddle-point problem \ref{SaddlepointProblem} exists if the following theorem holds true \citealp{NegriEtAlStokes}.
\begin{theorem}[Brezzi's theorem]
\label{BrezziTheorem}
A unique solution to saddle-point problem \eqref{SaddlepointProblem} will exist if the following conditions are satisfied,
\begin{itemize}
\item[$\left( i\right)$] $\mathcal{A}: X\times X\rightarrow\mathbb{R}$ is continuous and satisfies,
\begin{equation*}
\exists\ a_0 > 0,\ \text{such that}\ \inf_{\mathbf{y}\setminus\left\lbrace 0 \right\rbrace\in X_0}\frac{\mathcal{A}\left(\mathbf{y}, \mathbf{y}; \bm{\mu}\right)}{\| \mathbf{y}\|_X^2} \geq a_0,\ \forall\ \bm{\mu}\in \mathcal{D},
\end{equation*}
where, $X_0 = \left\lbrace \mathbf{y}\in X\ \text{such that}\ \mathcal{B}\left( \mathbf{y}, \bm{\kappa}; \bm{\mu}\right) = 0\ \forall\ \bm{\kappa}\in Z\right\rbrace$.
\item[$\left( ii\right)$] $\mathcal{B}: X\times Z\rightarrow\mathbb{R}$ is continuous and satisfies the following inf-sup condition,
\begin{equation*}
\exists\ b_0 > 0,\ \text{such that}\ \inf_{\bm{\kappa}\setminus\left\lbrace 0 \right\rbrace\in Z} \sup_{\mathbf{y}\setminus\left\lbrace 0 \right\rbrace\in X} \frac{\mathcal{B}\left(\mathbf{y}, \bm{\kappa}; \bm{\mu}\right)}{\|\mathbf{y}\|_X \|\bm{\kappa}\|_Z} \geq b_0,\ \forall\ \bm{\mu}\in \mathcal{D}.
\end{equation*}
\end{itemize}
\end{theorem}
 
\noindent We opt for an {\it optimize-then-discretize} approach, that is, we will first derive the coupled optimality system and, afterwards, we will introduce the numerical discretization of the problem. The optimization step can be performed by requiring first order optimality conditions, that is, $\nabla \mathcal{L}\left(\mathbf{x}, \bm{z}; \bm{\mu}\right)\left[ \mathbf{y}, \bm{\kappa}\right] = 0,\ \forall \mathbf{y}\in X, \bm{\kappa} = \left(\bm{\kappa}_{\bm{w}}, \bm{\kappa}_{q}\right)\in Z$. This gives rise to the following coupled KKT optimality system, which is then intended to be solved in one shot for $\left(\mathbf{x}\left(\bm{\mu}\right), \bm{z}\left(\bm{\mu}\right)\right)$:
\begin{equation}
\label{infiniteDimensionalOptimalitySystem}
\begin{cases}
\mathcal{A}\left( \mathbf{x}, \mathbf{y} ;\bm{\mu} \right) + \mathcal{B}\left(\mathbf{y}, \bm{z};\bm{\mu} \right) + \mathrm{e}\left(\mathbf{y}_{\bm{v}}, \bm{v}, \bm{w}; \bm{\mu}\right) + \mathrm{e}\left(\bm{v}, \mathbf{y}_{\bm{v}}, \bm{w}; \bm{\mu}\right) = \left\langle \mathcal{H}\left(\bm{\mu}\right), \mathbf{y}\right\rangle, &\ \forall\ \mathbf{y}\in X, \\
\mathcal{B}\left(\mathbf{x}, \bm{\kappa} ;\bm{\mu} \right) + \mathrm{e}\left(\bm{v}, \bm{v}, \bm{\kappa}_{\bm{w}}; \bm{\mu}\right) = \left\langle \mathcal{G}\left(\bm{\mu}\right), \bm{\kappa} \right\rangle, &\ \forall\ \bm{\kappa}\in Z.
\end{cases}
\end{equation}

\noindent The optimality system is non-linear in state, that is, $\mathrm{e}\left(\bm{v}, \bm{v}, \bm{\kappa}_{\bm{w}}\right)$ denotes convection term and $\mathrm{e}\left(\bm{v}, \bm{w}, \mathbf{y}_{\bm{v}}; \bm{\mu}\right)$ and $\mathrm{e}\left(\bm{w}, \bm{v}, \mathbf{y}_{\bm{v}}; \bm{\mu}\right)$ are trilinear terms arising in adjoint equation. Continuity of $\mathcal{A}$ and $\mathcal{B}$ is implied from the assumptions on the bilinear forms in weak formulation of state constraints. Conditions $\left( i\right)$ and $\left( ii\right)$ are satisfied when $S \equiv Z$.

\subsubsection{Numerical approximations}\label{subsec: high fidelity approximations}

\noindent In this section, we will introduce the discrete version of problem \eqref{infiniteDimensionalOptimalitySystem} and will rely on Galerkin finite element methods for its solution. We will refer to this finite dimensional discretized problem as ``{\it truth problem}'' in the next sections. Let us assume that the mesh discretization of domain $\Omega$, done in section \ref{sec: geometrical reconstruction}, has size $h\in \mathbb{N}$ such that $ 0 < h < \infty$. Furthermore, we consider finite dimensional solution spaces $S_h = \left( V_h, P_h \right)\subset S,\ U_h \subset U$ and $Z_h \subset Z$. Taking $X_h = S_h\times U_h\ \subset X$, the discretized finite dimensional optimal flow control problem is defined as:
\begin{PD}
\label{FEPOFCP}
\it Given $\bm{\mu}\in\mathcal{D}$, find $\left( \mathbf{x}_h\left(\bm{\mu}\right), \bm{z}_h\left(\bm{\mu}\right)\right)\ \in X_h\times Z_h$ such that,
\begin{equation}
\label{GalerkinFEApproximations}
\begin{cases}
\mathcal{A}\left( \mathbf{x}_h, \mathbf{y}_h ;\bm{\mu}\right) + \mathcal{B}\left(\mathbf{y}_h, \bm{z}_h ;\bm{\mu}\right) + \mathrm{e}\left(\mathbf{y}_{\bm{v}_h}, \bm{v}_h, \bm{w}_h; \bm{\mu}\right) + \mathrm{e}\left(\bm{v}_h, \mathbf{y}_{\bm{v}_h}, \bm{w}_h; \bm{\mu}\right) = \left\langle \mathcal{H}\left(\bm{\mu}\right), \mathbf{y}_h \right\rangle, &\ \forall\ \mathbf{y}_h \in X_h,\\
\mathcal{B}\left(\mathbf{x}_h, \bm{\kappa}_h ;\bm{\mu}\right) + \mathrm{e}\left(\bm{v}_h, \bm{v}_h, \bm{\kappa}_{\bm{w}_h}; \bm{\mu}\right) = \left\langle \mathcal{G}\left(\bm{\mu}\right), \bm{\kappa}_h \right\rangle, &\ \forall\ \bm{\kappa}_h\in Z_h.
\end{cases}
\end{equation}
\end{PD}

\noindent Subscript $h$ indicates that the dimensions of Galerkin finite element solution spaces rely on the mesh size $h$. We denote global dimension of {\it truth problem} by $\mathcal{N} = \mathcal{N}_{\bm{v}} + \mathcal{N}_{p} + \mathcal{N}_{\bm{u}} + \mathcal{N}_{\bm{w}} + \mathcal{N}_{q}$, where $\mathcal{N}_{\bm{\delta}}$ corresponds to the dimension of individual solution spaces for $\bm{\delta} = \bm{v}, p, \bm{u}, \bm{w}, q$. Moreover, we retain the assumption of equivalent state and adjoint spaces that is, $S_h \equiv Z_h$.

Now with $\bm{\phi}, \bm{\psi}$ and $\bm{\sigma}$ chosen as bases for velocity, pressure and control spaces, respectively, such that
\begin{equation*}
V_h = span\left\lbrace\bm{\phi}_i\right\rbrace_{i=1}^{\mathcal{N}_{\bm{v}}}, \quad P_h = span\left\lbrace\bm{\psi}_k\right\rbrace_{k=1}^{\mathcal{N}_{p}}, \quad U_h = span\left\lbrace\bm{\sigma}_l\right\rbrace_{l=1}^{\mathcal{N}_{\bm{u}}},
\end{equation*}
we define bijections $V_h \leftrightarrow \mathbb{R}^{\mathcal{N}_{\bm{v}}}$,  $P_h \leftrightarrow \mathbb{R}^{\mathcal{N}_{p}}$ and $U_h \leftrightarrow \mathbb{R}^{\mathcal{N}_{\bm{u}}}$ and give algebraic expansion of non-linear optimality system \eqref{GalerkinFEApproximations}\footnote{It is to be noted that owing to $S_h \equiv Z_h$, we can take the adjoint terms, with subscript $ad$ with the same bases functions as for the state terms, thus the components of the mass matrices of state and adjoint terms can be defined in a similar way.} below:
\begin{equation}
\label{algebraicFEMatrixForm}
\begin{bmatrix}
M\left(\bm{\mu}\right) + \tilde{E}\left(\mathbf{w}\left(\bm{\mu}\right); \bm{\mu}\right) & 0 & 0 & A_{ad}\left(\bm{\mu}\right) + E_{ad}\left(\mathbf{v}\left(\bm{\mu}\right); \bm{\mu}\right)  & B_{ad}\left(\bm{\mu}\right)\\
0 & 0 & 0 & B^T_{ad}\left(\bm{\mu}\right) & 0\\
0 & 0 & N\left(\bm{\mu}\right) & C_{ad}\left(\bm{\mu}\right) & 0\\
A\left(\bm{\mu}\right) + E\left(\mathbf{v}\left(\bm{\mu}\right); \bm{\mu}\right) & B^T\left(\bm{\mu}\right) & C\left(\bm{\mu}\right) & 0 & 0\\
B\left(\bm{\mu}\right) & 0 & 0 & 0 & 0  
\end{bmatrix}
\begin{bmatrix}
\mathbf{v}\left(\bm{\mu}\right)\\
\mathbf{p}\left(\bm{\mu}\right)\\
\mathbf{u}\left(\bm{\mu}\right)\\
\mathbf{w}\left(\bm{\mu}\right)\\
\mathbf{q}\left(\bm{\mu}\right)
\end{bmatrix} = 
\begin{bmatrix}
\mathbf{h}\left(\bm{\mu}\right)\\
\mathbf{0}\\
\mathbf{0}\\
\mathbf{f}\left(\bm{\mu}\right)\\
\mathbf{g}\left(\bm{\mu}\right)

\end{bmatrix}.
\end{equation}

\noindent Here, $A\left(\bm{\mu}\right)\in \mathbb{R}^{\mathcal{N}_{\bm{v}}}\times \mathbb{R}^{\mathcal{N}_{\bm{v}}}$, $B\left(\bm{\mu}\right) \in \mathbb{R}^{\mathcal{N}_{p}}\times \mathbb{R}^{\mathcal{N}_{\bm{v}}}$, $C\left(\bm{\mu}\right)\in \mathbb{R}^{\mathcal{N}_{\bm{u}}}\times \mathbb{R}^{\mathcal{N}_{\bm{v}}}$, $M\left(\bm{\mu}\right) \in \mathbb{R}^{\mathcal{N}_{\bm{v}}}\times \mathbb{R}^{\mathcal{N}_{\bm{v}}}$ and $N\left(\bm{\mu}\right)\in \mathbb{R}^{\mathcal{N}_{\bm{u}}}\times \mathbb{R}^{\mathcal{N}_{\bm{u}}}$ are stiffness matrices associated to the finitely discretized versions of bilinear forms in \eqref{bilinearFormsNavierStokesConstraints} and for $1 \leq i, j\leq \mathcal{N}_{\bm{v}}, 1 \leq k\leq \mathcal{N}_{p}, 1 \leq l, r\leq \mathcal{N}_{\bm{u}}$, are defined as:
\begin{equation*}
\left( A\left(\bm{\mu}\right)\right)_{ij} = \mathrm{a}\left(\bm{\phi}_i, \bm{\phi}_j;\bm{\mu}\right), \quad \left( B\left(\bm{\mu}\right)\right)_{ik} = \mathrm{b}\left(\bm{\psi}_k, \bm{\phi}_i;\bm{\mu}\right), \quad \left( C\left(\bm{\mu}\right)\right)_{il} = \mathrm{c}\left(\bm{\sigma}_l, \bm{\phi}_i;\bm{\mu}\right),
\end{equation*}
\begin{equation*}
\left( M\left(\bm{\mu}\right)\right)_{ij} = \mathrm{m}\left(\bm{\phi}_i, \bm{\phi}_j;\bm{\mu}\right), \qquad \left( N\left(\bm{\mu}\right)\right)_{lr} = \mathrm{n}\left(\bm{\sigma}_r, \bm{\sigma}_l;\bm{\mu}\right).
\end{equation*}
Moreover, $\left( E\left(\mathbf{v}\left(\bm{\mu}\right); \bm{\mu}\right)\right)_{ij} = \sum\limits_{k = 1}^{\mathcal{N}_{\bm{v}}} v_h^k \left(\bm{\mu}\right) \mathrm{e}\left(\bm{\phi}_k, \bm{\phi}_j, \bm{\phi}_i; \bm{\mu}\right)$ and $\left( \tilde{E}\left(\mathbf{w}\left(\bm{\mu}\right); \bm{\mu}\right)\right)_{ij} = \sum\limits_{k = 1}^{\mathcal{N}_{\bm{v}}} w_h^k \left(\bm{\mu}\right) \mathrm{e}\left(\bm{\phi}_k, \bm{\phi}_j, \bm{\phi}_i; \bm{\mu}\right)$ are associated to the non-linear and trilinear terms arising in state and adjoint equations respectively.

Furthermore, elements of the vector on left hand side of equation \eqref{algebraicFEMatrixForm} are vectors of coefficients, defined as:
\begin{equation*}
\mathbf{v} = 
\begin{bmatrix}
v_h^1 \\
v_h^2 \\
\vdots \\
v_h^{\mathcal{N}_{\bm{v}}}
\end{bmatrix} \in \mathbb{R}^{\mathcal{N}_{\bm{v}}},
\quad
\mathbf{p} = 
\begin{bmatrix}
p_h^1 \\
p_h^2 \\
\vdots \\
p_h^{\mathcal{N}_{p}}
\end{bmatrix} \in \mathbb{R}^{\mathcal{N}_{p}},
\quad
\mathbf{u} = 
\begin{bmatrix}
u_h^1 \\
u_h^2 \\
\vdots \\
u_h^{\mathcal{N}_{\bm{u}}} 
\end{bmatrix} \in \mathbb{R}^{\mathcal{N}_{\bm{u}}},
\quad
\mathbf{w} = 
\begin{bmatrix}
w_h^1 \\
w_h^2 \\
\vdots \\
w_h^{\mathcal{N}_{\bm{v}}} 
\end{bmatrix} \in \mathbb{R}^{\mathcal{N}_{\bm{v}}},
\quad
\mathbf{q} = 
\begin{bmatrix}
q_h^1 \\
q_h^2 \\
\vdots \\
q_h^{\mathcal{N}_{p}}
\end{bmatrix} \in \mathbb{R}^{\mathcal{N}_{p}}.
\end{equation*}

\noindent The optimality system \eqref{algebraicFEMatrixForm} can be solved directly for state, control and adjoint through iterative methods for non-linear PDEs such as Newton method. In this work, we have used the open-source libraries FEniCS\citealp{FENICS2012, FENICS2015} and {\it multiphenics}\citealp{multiphenics} for full order simulations. The latter is an internal library developed and maintained at SISSA mathlab, to attain an easy prototyping of block-structured problems in FEniCS.

Uniqueness of state and adjoint pressure in Stokes and Navier-Stokes equations is ensured by Ladyzhenskaya-Brezzi-Babu{\u s}ka (LBB) inf-sup condition \citealp{AliEtAl2019,BallarinEtAlStabilization}, defined below:
\begin{equation}
\label{LBB inf-sup FE}
\exists\ \beta > 0,\ \text{such that}\ \inf\limits_{{p}_h \setminus\left\lbrace 0\right\rbrace\in P_h} \sup\limits_{{\bm{v}}_h \setminus\left\lbrace 0\right\rbrace\in V_h} \frac{\mathrm{b}\left( p_h, {\bm{v}}_h\right)}{\|{p}_h\|_{P_h} \|{\bm{v}}_h\|_{V_h}} \geq \beta.
\end{equation}
\noindent To satisfy this condition, we use stable Taylor-hood, that is, $\mathbb{P}2 - \mathbb{P}1$ solution spaces for velocity and pressure. As aforementioned, $\mathcal{N}$ depends upon discretization size and hence, number of mesh elements. The number of mesh elements is required to be sufficiently large in numerical hemodynamics modeling, in order to attain reliable numerical solution. Thus, such problems are composed of large number of degrees of freedom. The computations for a single value of $\bm{\mu}$ in such cases can take a few hours of CPU time and additionally, for accurate patient-specific hemodynamics modeling, one has to inevitably take into consideration multiple hemodynamics scenarios modeled by tuning the parameters. Thus, computations need to be performed in repetitive environment, increasing the cost to order of days, making it necessary to employ computationally efficient numerical methods while retaining reliability of the solutions. For this reason, we add another step to {\it optimize-then-discretize}, that is ``{\it reduce}''.

\subsection{Proper Orthogonal Decomposition ({POD}) -- Galerkin Approximations}
In this section, we will discuss the {\it optimize-discretize-reduce} approach for computations in repetitive environment for the parametrized optimal flow control in computational hemodynamics modeling. This can be achieved by using solution spaces with much lower dimensions, constructed by reduced order methods. These methods are essentially built upon the solution of {\it truth problems}, which we will refer to as {\it truth-approximations} or {\it snapshots} and these snapshots will be calculated through Galerkin finite element methods.

In this work, we have utilized {\it proper orthogonal decomposition} ({POD})\citealp{HesthavenEtAl2015, KunischEtAl, muller, BallarinEtAlStabilization} to construct reduced order spaces, however a {\it greedy algorithm}\citealp{HesthavenEtAl2015, NegriEtAlBook} using residual-based error estimators can be alternatively employed. Both techniques have been well-applied to optimal flow control problems\citealp{NegriEtAlStokes, NegriEtAlElliptic, StrazzulloEtAl, BaderEtAl2016, Dede2012, KarcherEtAl2017} and patient-specific computational cardiovascular modelling\citealp{BallarinEtAl2016, BallarinEtAl2017, TezzeleEtAl2018, AuricchioEtAl2018}. We summarize the algebraic details, following Rozza et al. \citealp{HesthavenEtAl2015, QuarteroniRozza2007} and Ballarin et al. \citealp{BallarinEtAlStabilization}, below.

Let $\bm{\Lambda}\subset \mathcal{D}$ be a finitely sampled subset of parameters. We construct snapshot matrices $\bm{\mathfrak{X}_{\delta}}$ that contain solutions to {\it truth problem} \eqref{GalerkinFEApproximations} by solving discretized optimality system \eqref{algebraicFEMatrixForm} for each $\bm{\mu}^\mathrm{i} \in \bm{\Lambda},\ \mathrm{i} = 1, \cdots, |\bm{\Lambda}|$. Thus,
\begin{equation}
\label{snapshotMatrices}
\bm{\mathfrak{X}_{\bm{\delta}}} =
\begin{bmatrix}
\bm{\delta}_h\left(\bm{\mu}^1\right) & \bm{\delta}_h\left(\bm{\mu}^2\right) & \cdots & \bm{\delta}_h\left(\bm{\mu}^{|\bm{\Lambda}|}\right)
\end{bmatrix}, \qquad \bm{\delta} = \bm{v}, p, \bm{u}, \bm{w}, q,\ \text{and}\ 1\leq h\leq \mathcal{N}_{\bm{\delta}}.
\end{equation}
To construct {POD} bases, we solve the following eigenvalue problems:
\begin{equation}
\label{EigenvalueProblem}
\mathbb{A}^{\bm{\delta}} \bm{\rho}^{\bm{\delta}}_{\mathrm{i}} =  \lambda_{\mathrm{i}}^{\bm{\delta}} \bm{\rho}^{\bm{\delta}}_{\mathrm{i}}, \qquad \mathrm{i} = 1, \cdots, |\bm{\Lambda}|,
\end{equation}
and keep $n = 1,\cdots, \mathrm{N}_{max}, \mathbb{N} \ni \mathrm{N}_{max} \ll |\bm{\Lambda}|$ eigenvalue-eigenvector pairs $\left( \lambda^{\bm{\delta}}_{n}, \bm{\rho}^{\bm{\delta}}_{n}\right)$, retaining sufficient average energy of the snapshots.
In equation \eqref{EigenvalueProblem}, $\mathbb{A}^{\bm{\delta}} = \frac{1}{|\bm{\Lambda}|}  \bm{\mathfrak{X}_{\delta}}^T\bm{\mathfrak{X}_{\delta}}\ \in \mathbb{R}^{|\bm{\Lambda}| \times |\bm{\Lambda}|}$ is the correlation matrix corresponding to snapshots. Orthonormal {POD} bases are constructed from the retained $\mathrm{N}_{max}$ eigenvectors and are defined as:
\label{reducedOrderBasis}
\begin{equation*}
\bm{\tilde{\phi}}^{\bm{v}}_n = \frac{1}{\sqrt{\lambda_n^{\bm{v}}}}\mathfrak{X}_{\bm{v}}\bm{\rho}_n^{\bm{v}}, \qquad \bm{\tilde{\psi}}^{p}_n = \frac{1}{\sqrt{\lambda_n^{p}}}\mathfrak{X}_{p}\bm{\rho}_n^{p},
\end{equation*}
\begin{equation*}
\bm{\tilde{\sigma}}^{\bm{u}}_n = \frac{1}{\sqrt{\lambda_n^{\bm{u}}}}\mathfrak{X}_{\bm{u}}\bm{\rho}_n^{\bm{u}},
\end{equation*}
\begin{equation*}
\bm{\tilde{\zeta}}^{\bm{w}}_n = \frac{1}{\sqrt{\lambda_n^{\bm{w}}}}\mathfrak{X}_{\bm{w}}\bm{\rho}_n^{\bm{w}}, \qquad \bm{\tilde{\zeta}}^{q}_n = \frac{1}{\sqrt{\lambda_n^{q}}}\mathfrak{X}_{q}\bm{\rho}_n^{q}.
\end{equation*}
\begin{figure}
\hspace*{-0.5cm}\begin{tikzpicture}[framed, very thick]
\node[block, draw, text width = 4em](p-snode1){\includegraphics[width = 0.5\textwidth]{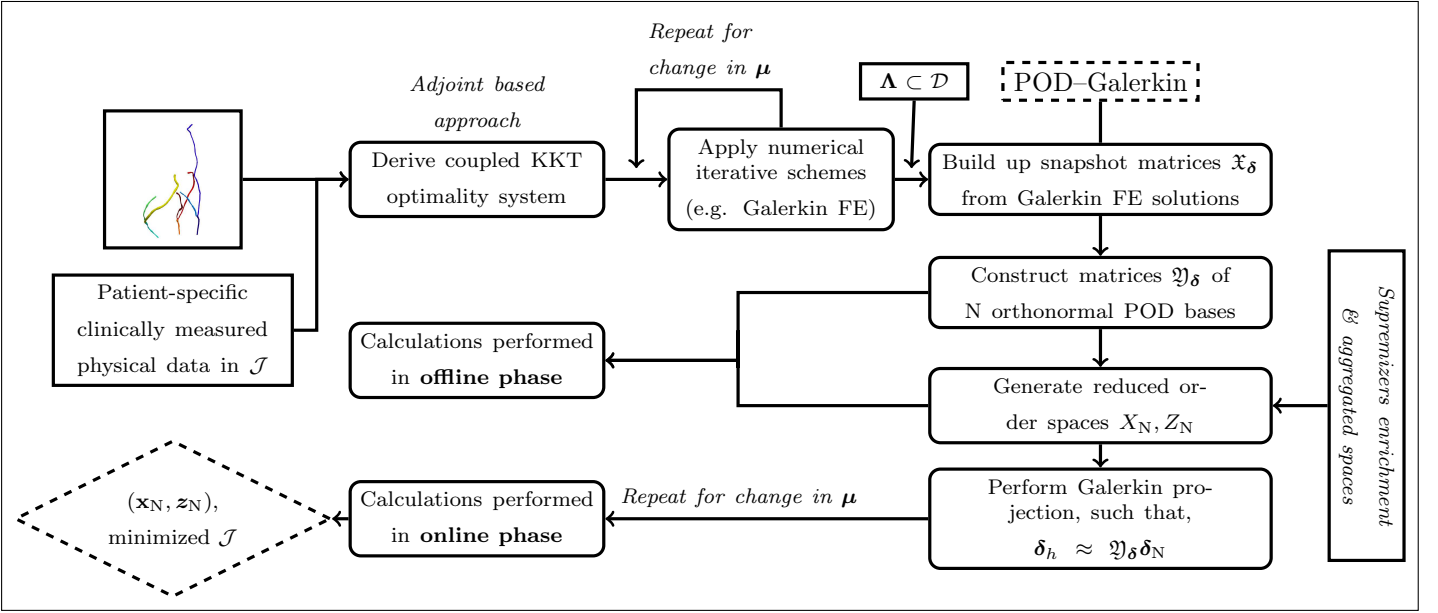}};
\node[block, draw, text width = 7.5em, below of  = p-snode1, node distance = 2cm](p-snode2){\scriptsize Patient-specific clinically measured physical data in $\mathcal{J}$};
\node[block, draw, text width = 8em, rounded corners, right of = p-snode1, node distance = 4cm](node1){\scriptsize Derive coupled KKT optimality system};
\node[block, above of = node1, node distance = 1cm, text width = 7em](node2){\scriptsize \it Adjoint based approach};
\node[block, draw, rounded corners, right of = node1, node distance = 4cm, text width = 7em](node3){\scriptsize Apply numerical iterative schemes\\ (e.g. Galerkin FE)};
\node[block, draw, rounded corners, right of = node3, node distance = 4.2cm, text width = 11em](node5){\scriptsize Build up snapshot matrices $\mathfrak{X}_{\bm{\delta}}$ from Galerkin FE solutions};
\node[block, draw, dashed, above of = node5, node distance = 1.3cm, text width = 6em](node4){\small {POD}--Galerkin};
\node[block, draw, rounded corners, below of = node5, node distance = 1.5cm, text width = 11em](node6){\scriptsize Construct matrices $\mathfrak{Y}_{\bm{\delta}}$ of $\mathrm{N}$ orthonormal {POD} bases};
\node[block, draw, rounded corners, below of = node6, node distance  = 1.5cm, text width = 11em](node7){\scriptsize Generate reduced order spaces $X_{\mathrm{N}}, Z_{\mathrm{N}}$};
\node[block, draw, rounded corners, below of = node7, node distance  = 1.5cm, text width = 11em](node8){\scriptsize Perform Galerkin projection, such that,\\ $\bm{\delta}_h \approx \mathfrak{Y}_{\bm{\delta}}\bm{\delta}_{\mathrm{N}}$};
\node[block, draw, rounded corners, below of = node1, node distance = 2.4cm, text width =8em](node9){\scriptsize Calculations performed in {\it \bf offline phase}};
\node[block, draw, rounded corners, below of = node9, node distance = 2.1cm, text width =8em](node10){\scriptsize Calculations performed in {\it \bf online phase}};
\node[block, draw, right of = node7, node distance = 3.5cm, text width = 10em, rotate = 270](node11){\scriptsize \it Supremizers enrichment\\
\& aggregated spaces};
\node[block, draw, left of = node4, node distance = 2.45cm, text width = 3em](node12){\scriptsize $\bm{\Lambda}\subset \mathcal{D}$};
\node[decision, draw, dashed, left of = node10, text width = 7em, node distance = 4cm](node13){\scriptsize $\left(\mathbf{x}_{\mathrm{N}}, \bm{z}_{\mathrm{N}}\right)$,\\
minimized $\mathcal{J}$};

\draw[->] (p-snode1.east) -- (node1.west);
\draw[->] (p-snode2.east) -- ++(0.3, 0) -- ++ (0, 1) |- (node1.west);
\draw[->] (node1) -- node[midway](arrownode1){}(node3);
\draw[->] (node3.north) -- ++(0, 0.5) -- ++(-0.5, 0) -| node[midway, above right, text width = 5em]{\it \scriptsize Repeat for change in $\bm{\mu}$}(arrownode1.north);
\draw[->] (node3) -- node[midway](arrownode2){}(node5);
\draw[-] (node4) -- (node5);
\draw[->] (node5) -- (node6);
\draw[->] (node6) -- (node7);
\draw[->] (node7) -- (node8);
\draw[->] (node6.west) -- ++(-2.5, 0) -- ++(0, -1) |- (node9.east);
\draw[->] (node7.west) -- ++(-2.5, 0) -- ++(0, 1) |- (node9.east);
\draw[->] (node8.west) -- node[midway, above, text width = 10em](arrownode3){\it \scriptsize Repeat for change in $\bm{\mu}$}(node10.east);
\draw[->] (node11.south) -- (node7.east);
\draw[->] (node12.south) -- (arrownode2.north);
\draw[->] (node10.west) -- (node13.east);
\end{tikzpicture}
\caption{Flow chart of reduced order optimal flow control pipeline}
\label{tikzfigure2}
\end{figure}
\noindent These {POD} bases generate the reduced order spaces, which will be marked by subscript `$\mathrm{N}$'. However, to ensure uniqueness and stability of the solutions, we need to satisfy Brezzi inf-sup condition \eqref{Brezzi inf-sup ROM} and Ladyzhenskaya-Brezzi-Bab{\u s}ka (LBB) inf-sup condition \eqref{LBB inf-sup ROM}. At the reduced order level, these conditions are defined as below:
\begin{equation}
\label{Brezzi inf-sup ROM}
\exists\ b_0 > 0,\ \text{such that}\ \inf\limits_{{\bm{\kappa}}_{\mathrm{N}}\setminus\left\lbrace 0 \right\rbrace\in Z_{\mathrm{N}}} \sup\limits_{{\mathbf{y}}_{\mathrm{N}}\setminus\left\lbrace 0 \right\rbrace\in X_{\mathrm{N}}} \frac{{\mathcal{B}}\left({\mathbf{y}}_{\mathrm{N}}, {\bm{\kappa}}_{\mathrm{N}}\right)}{\| {\mathbf{y}}_{\mathrm{N}}\|_{X_{\mathrm{N}}} \|{\bm{\kappa}}_{\mathrm{N}}\|_{Z_{\mathrm{N}}}} \geq b_0,
\end{equation}
\begin{equation}
\label{LBB inf-sup ROM}
\exists\ \beta_0 > 0,\ \text{such that}\ \inf\limits_{{p}_{\mathrm{N}} \setminus\left\lbrace 0\right\rbrace\in P_{\mathrm{N}}} \sup\limits_{\bm{v}_{\mathrm{N}} \setminus\left\lbrace 0\right\rbrace\in V_{\mathrm{N}}} \frac{\mathrm{b}\left(p_{\mathrm{N}}, {\bm{v}}_{\mathrm{N}}\right)}{\|{p}_{\mathrm{N}}\|_{P_{\mathrm{N}}} \|{\bm{v}}_{\mathrm{N}}\|_{V_{\mathrm{N}}}} \geq \beta_0.
\end{equation}
As discussed earlier, Brezzi's inf-sup condition \eqref{Brezzi inf-sup ROM} will be satisfied if the equivalence relation between state and adjoint spaces hold, which is not guaranteed at a reduced order level. In order to ensure the equivalence relation, we aggregate state and adjoint spaces \citealp{NegriEtAlStokes, StrazzulloEtAl}, that is, we consider equivalent state and adjoint spaces for velocity and pressure, generated by {POD} bases of both state and adjoint variables.

Furthermore, to satisfy LBB inf-sup condition \eqref{LBB inf-sup ROM} at reduced order level, two approaches are being employed in literature. One is to introduce stabilization terms in state constraints, we refer the reader to Ali et al.\citealp{AliThesis, AliEtAl2019} and Deparis et al.\citealp{DeparisEtAl2012} for more details on this. We implement the other approach in this paper, that is to use supremizer operators to enlarge the velocity spaces \citealp{RozzaAndVeroy, BallarinEtAlStabilization, NegriEtAlStokes, StrazzulloEtAl}. Thus, we define {\it supremizers} $\mathcal{T}^{\mathrm{i}}_{\bm{v}_h}: Z_{p_h}\rightarrow V_h$ and $\mathcal{T}^{\mathrm{i}}_{\bm{w}_h}: P_h\rightarrow Z_{\bm{v}_h}$ below:
\begin{equation*}
\left( \mathcal{T}^{\mathrm{i}}_{\bm{v}_h}q_h, \bm{v}_h\right) = \mathrm{b}\left(q_h, \bm{v}_h; \bm{\mu}^{\mathrm{i}}\right),\quad \mathrm{i} = 1,\cdots, |\bm{\Lambda}|
\end{equation*}
\begin{equation*}
\left( \mathcal{T}^{\mathrm{i}}_{\bm{w}_h}p_h, \bm{w}_h\right) = \mathrm{b}\left(p_h, \bm{w}_h; \bm{\mu}^\mathrm{i}\right),\quad \mathrm{i} = 1,\cdots, |\bm{\Lambda}|
\end{equation*}
and add them respectively to state and adjoint velocity spaces. State and adjoint velocity spaces are defined as below:
\begin{equation}
V_{\mathrm{N}} = span\left\lbrace \bm{\tilde{\phi}}^{\bm{v}}_n, \bm{\hat{\phi}}^{\bm{v}}_n, \bm{\tilde{\zeta}}^{\bm{w}}_n, \bm{\hat{\zeta}}^{\bm{w}}_n,\quad n = 1,\cdots, \mathrm{N}_{max}\right\rbrace \equiv Z_{{\bm{v}}_{\mathrm{N}}},
\end{equation}
where, $\bm{\hat{\phi}}^{\bm{w}}_n$ and $\bm{\hat{\zeta}}^{\bm{w}}_n$ are the {POD} modes for supremizers for state and adjoint velocity, respectively. State and adjoint pressure spaces are defined as:
\begin{equation}
P_{\mathrm{N}} = span\left\lbrace \bm{\tilde{\psi}}^{p}_n, \bm{\tilde{\zeta}}^{q}_n,\quad n = 1,\cdots, \mathrm{N}_{max}\right\rbrace \equiv Z_{p_{\mathrm{N}}}
\end{equation}
Thus, the enriched aggregated reduced order state and adjoint spaces are then defined as below:
\begin{equation}
\label{eq: supremizer enriched and aggregated reduced order spaces}
S_{\mathrm{N}} = \left( V_{\mathrm{N}} \oplus \mathcal{T}^{\bm{v}}_{\mathrm{N}}\right) \times P_{\mathrm{N}} \times \left( Z_{{\bm{v}}_{\mathrm{N}}} \oplus \mathcal{T}^{\bm{w}}_{\mathrm{N}}\right)\times Z_{p_{\mathrm{N}}} \equiv Z_{\mathrm{N}},
\end{equation}
and reduced order control space is defined as:
\begin{equation}
\label{eq: reduced order control space}
U_\mathrm{N} = span\left\lbrace \bm{\tilde{\sigma}}^{\bm{u}}_n,\quad n = 1,\cdots, \mathrm{N}_{max}\right\rbrace.
\end{equation}
We denote dimensions of reduced order spaces by $\mathrm{N}_{\bm{\delta}}$ for $\delta = \bm{v}, p, \bm{u}, \bm{w}, q$. Thus, $\mathrm{N}_{\bm{v}} = 4\mathrm{N}_{max} = \mathrm{N}_{\bm{w}}, \mathrm{N}_{p} = 2\mathrm{N}_{max} = \mathrm{N}_{q}$ and $\mathrm{N}_{\bm{u}} = \mathrm{N}_{max}$. The reduced bases matrices are defined as below:
\begin{equation*}
\mathfrak{Y}_{\bm{v}} = \left[\left(\bm{\tilde{\phi}}^{\bm{v}}_1 + \bm{\hat{\phi}}^{\bm{v}}_1\right)\ \left(\bm{\tilde{\phi}}^{\bm{v}}_2 + \bm{\hat{\phi}}^{\bm{v}}_2\right)\ \cdots\ \left(\bm{\tilde{\phi}}^{\bm{v}}_{\mathrm{N}_{max}} + \bm{\hat{\phi}}^{\bm{v}}_{\mathrm{N}_{max}}\right) \right]^T, \quad \mathfrak{Y}_{p} = \left[\bm{\tilde{\psi}}^{p}_1\ \bm{\tilde{\psi}}^{p}_2\ \cdots\ \bm{\tilde{\psi}}^{p}_{\mathrm{N}_{max}} \right]^T,
\end{equation*}
\begin{equation*}
\mathfrak{Y}_{\bm{u}} = \left[\bm{\tilde{\sigma}}^{\bm{u}}_1\ \bm{\tilde{\sigma}}^{\bm{u}}_2\ \cdots\ \bm{\tilde{\sigma}}^{\bm{u}}_{\mathrm{N}_{max}} \right]^T,
\end{equation*}
\begin{equation*}
\mathfrak{Y}_{\bm{w}} = \left[\left(\bm{\tilde{\zeta}}^{\bm{w}}_1 + \bm{\hat{\zeta}}^{\bm{w}}_1\right)\ \left(\bm{\tilde{\zeta}}^{\bm{w}}_2 + \bm{\hat{\zeta}}^{\bm{w}}_2\right)\ \cdots\ \left(\bm{\tilde{\zeta}}^{\bm{w}}_{\mathrm{N}_{max}} +  \bm{\hat{\zeta}}^{\bm{w}}_{\mathrm{N}_{max}}\right) \right]^T \qquad \text{and}\qquad  \mathfrak{Y}_{q} = \left[\bm{\tilde{\zeta}}^{q}_1\ \bm{\tilde{\zeta}}^{q}_2\ \cdots\ \bm{\tilde{\zeta}}^{q}_{\mathrm{N}_{max}} \right]^T.
\end{equation*}
Now, with  $X_{\mathrm{N}} = S_{\mathrm{N}}\times U_{\mathrm{N}}$, we write reduced order version of parametrized optimal flow control problem as:
\begin{PD}
\label{pb: rom}
\it Given $\bm{\mu}\in \mathcal{D}$, find $\left( {\mathbf{x}}_{\mathrm{N}}\left(\bm{\mu}\right), {\bm{z}}_{\mathrm{N}}\left(\bm{\mu}\right)\right)\ \in X_{\mathrm{N}}\times Z_{\mathrm{N}}$ such that for $\mathfrak{Y} = diag \left(\mathfrak{Y}_{\bm{v}}, \mathfrak{Y}_{p}, \mathfrak{Y}_{\bm{u}}, \mathfrak{Y}_{\bm{w}}, \mathfrak{Y}_{q} \right)$,
\begin{equation}
\label{algebraicPODMatrixForm, Navierstokes, mu}
\mathfrak{Y}^T\mathbb{W} = \mathbf{0},
\end{equation}
where,
\begin{equation*}
\scriptsize
\mathbb{W} = 
\begin{bmatrix}
\left( M\left(\bm{\mu}\right) + \tilde{E}\left( \mathfrak{Y}_{\bm{w}}{\bm{w}}_{\mathrm{N}}\left(\bm{\mu}\right); \bm{\mu}\right)\right) \mathfrak{Y}_{\bm{v}} \bm{v}_{\mathrm{N}}\left(\bm{\mu}\right) + \left( A_{ad}\left(\bm{\mu}\right) + E_{ad} \left( \mathfrak{Y}_{\bm{v}}\mathbf{v}_{\mathrm{N}}\left(\bm{\mu}\right); \bm{\mu}\right)\right) \mathfrak{Y}_{\bm{w}} \bm{w}_{\mathrm{N}}\left(\bm{\mu}\right) + B^T_{ad}\left(\bm{\mu}\right)\mathfrak{Y}_{q}\mathbf{q}_{\mathrm{N}}\left(\bm{\mu}\right) - \mathbf{h}\left(\bm{\mu}\right) \\
B_{ad}\left(\bm{\mu}\right)\mathfrak{Y}_{\bm{w}}\mathbf{w}_{\mathrm{N}}\left(\bm{\mu}\right)\\
N\left(\bm{\mu}\right)\mathfrak{Y}_{\bm{u}}\mathbf{u}_{\mathrm{N}}\left(\bm{\mu}\right) + C\left(\bm{\mu}\right)\mathfrak{Y}_{\bm{w}}\mathbf{w}_{\mathrm{N}}\left(\bm{\mu}\right) \\
\left( A\left(\bm{\mu}\right) + E \left( \mathfrak{Y}_{\bm{v}}\mathbf{v}_{\mathrm{N}}\left(\bm{\mu}\right); \bm{\mu}\right)\right) \mathfrak{Y}_{\bm{v}} \bm{v}_{\mathrm{N}}\left(\bm{\mu}\right) + B^T\left(\bm{\mu}\right)\mathfrak{Y}_{p}\mathbf{p}_{\mathrm{N}} + C\left(\bm{\mu}\right)\mathfrak{Y}_{\bm{u}}\mathbf{u}_{\mathrm{N}}\left(\bm{\mu}\right)\left(\bm{\mu}\right) - \mathbf{f}\left(\bm{\mu}\right) \\
B\left(\bm{\mu}\right)\mathfrak{Y}_{\bm{v}}\mathbf{v}_{\mathrm{N}}\left(\bm{\mu}\right) - \mathbf{g}\left(\bm{\mu}\right)
\end{bmatrix}
\end{equation*}
\end{PD}
Now, the reduced order mass matrices can be written as:
\begin{align*}
M_{\mathrm{N}}\left(\bm{\mu}\right) &= \mathfrak{Y}_{\bm{v}}^TM\left(\bm{\mu}\right)\mathfrak{Y}_{\bm{v}},\qquad N_{\mathrm{N}}\left(\bm{\mu}\right) = \mathfrak{Y}_{\bm{u}}^TN\left(\bm{\mu}\right)\mathfrak{Y}_{\bm{u}}, \qquad A_{\mathrm{N}}\left(\bm{\mu}\right) = \mathfrak{Y}_{\bm{w}}^TA\left(\bm{\mu}\right)\mathfrak{Y}_{\bm{v}},\\
A_{{ad}_\mathrm{N}}\left(\bm{\mu}\right) &= \mathfrak{Y}_{\bm{v}}^TA_{ad}\left(\bm{\mu}\right)\mathfrak{Y}_{\bm{w}},\qquad B_{\mathrm{N}} = \mathfrak{Y}_{q}^TB\left(\bm{\mu}\right)\mathfrak{Y}_{\bm{v}}, \qquad B_{{ad}_{\mathrm{N}}} = \mathfrak{Y}_{p}^TB_{ad}\left(\bm{\mu}\right)\mathfrak{Y}_{\bm{w}},\\
E_{\mathrm{N}} \left( \cdot ,\ \cdot ; \bm{\mu}\right) &= \mathfrak{Y}_{\bm{w}}^TE\left( \cdot ,\ \cdot; \bm{\mu}\right) \mathfrak{Y}_{\bm{v}},\qquad E_{{ad}_\mathrm{N}} \left( \cdot ,\ \cdot ; \bm{\mu}\right) =  \mathfrak{Y}_{\bm{v}}^TE_{ad}\left( \cdot ,\ \cdot; \bm{\mu}\right) \mathfrak{Y}_{\bm{w}},
\end{align*}
and 
\begin{equation*}
\tilde{E}_{\mathrm{N}} \left( \cdot ,\ \cdot ; \bm{\mu}\right) = \mathfrak{Y}_{\bm{v}}^T\tilde{E} \left( \cdot ,\ \cdot; \bm{\mu}\right) \mathfrak{Y}_{\bm{v}}.
\end{equation*}
Furthermore, for $\bm{\delta} = \bm{v}, p, \bm{u}, \bm{w}, q
$, we denote the vectors of reduced order coefficients calculated through Galerkin projection with ${\mathbf{\bm{\delta}}}_{\mathrm{N}}\in \mathbb{R}^{\mathrm{N}_{\bm{\delta}}}$  such that $
\mathfrak{Y}_{\bm{\delta}} \bm{\delta}_{\mathrm{N}} \approx \bm{\delta}_h$.

Total dimensions of reduced order problem are $\mathrm{N} = \mathrm{N}_{\bm{v}} + \mathrm{N}_{p} + \mathrm{N}_{\bm{w}} + \mathrm{N}_{q} + \mathrm{N}_{\bm{u}} = 13\mathrm{N}_{max}$, where $\mathcal{N}_{\bm{\delta}} \gg \mathrm{N}_{\bm{\delta}},\ \bm{\delta} = \bm{v}, p, \bm{u}, \bm{w}, q$, thus, $\mathrm{N} \ll \mathcal{N}$. Figure \ref{tikzfigure2} shows a flow chart summarizing implementation of the reduced order optimal flow control pipeline.

Computational efficiency of the ``{\it reduce}" step arises from decomposition of the procedure into two separate phases, namely an offline phase and an online phase. This phase decomposition relies on crucial assumption of affine decomposition for reduced order problems, generally defined as:
\begin{equation*}
\mathcal{A}\left(\mathbf{x}, \mathbf{y}; \bm{\mu}\right) = \sum\limits_{q = 1}^{Q_a} \theta^{q}\left(\bm{\mu}\right) \mathcal{A}^{q}\left(\mathbf{x}, \mathbf{y}\right), \quad \mathcal{B}\left(\mathbf{x}, \bm{\kappa}; \bm{\mu}\right) = \sum\limits_{q = 1}^{Q_b} \theta^{q}\left(\bm{\mu}\right) \mathcal{B}^{q}\left(\mathbf{x}, \bm{\kappa}\right),
\end{equation*}
\begin{equation*}
\mathrm{e}\left(\bm{v}, \bm{v}, \bm{\kappa}_{\bm{w}}; \bm{\mu}\right) = \sum\limits_{q = 1}^{Q_e} \theta^{q}\left(\bm{\mu}\right)\mathrm{e}^{q}\left(\bm{v}, \bm{v}, \bm{\kappa}_{\bm{w}}\right).
\end{equation*}

In non-affinely parametrized problems\citealp{Rozza2009}, empirical interpolation method can be used to approximate the affine decomposition. The offline phase comprises of collecting the snapshots for all $\bm{\mu}\in \bm{\Lambda}$ and the generation of reduced order spaces defined by equations \ref{eq: supremizer enriched and aggregated reduced order spaces} and \eqref{eq: reduced order control space} from {POD} bases. Computational cost of this phase relies on dimensions of full order solution spaces, that is, the time taken to solve $|\bm{\Lambda}| \times \mathcal{N}$ truth problem. Thus, this phase is computationally very expensive, however it has to be performed only once. Afterwards, for different values of $\bm{\mu}$, perform the online phase which assembles the mathematical terms by combining the corresponding parameter-dependent and already stored parameter-independent components, perform Galerkin projection to calculate reduced order coeffients and solve the reduced order problem \ref{algebraicPODMatrixForm, Navierstokes, mu}. This phase has a computational cost independent of $\mathcal{N}$ and is therefore quite inexpensive.

The reduced order optimal flow control pipeline, proposed in this work, is implemented using RBniCS\citealp{HesthavenEtAl2015, RBNICS} and {\it multiphenics}, libraries developed and maintained at SISSA mathlab. The former library is used when dealing with distributed control problems while the latter is designed specifically to deal with block-structured boundary control problems.

\section{Numerical Results} \label{sec: numerical results}

\begin{figure}
\centering
\begin{tabular}{lllllllll}
\multicolumn{1}{c}{\multirow{-12}{*}{\rotatebox{90}{(a). $\Omega_a$}}} & \multicolumn{1}{c}{\includegraphics[width = 0.089\textwidth, keepaspectratio]{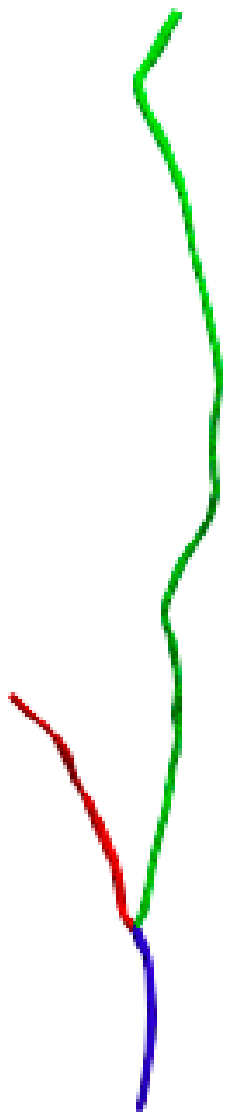}} & \multicolumn{1}{c}{\hspace*{1cm}} & \multicolumn{1}{c}{\multirow{-12}{*}{\rotatebox{90}{(b). $\Omega_b$}}} & \multicolumn{1}{c}{\includegraphics[width = 0.1\textwidth, keepaspectratio]{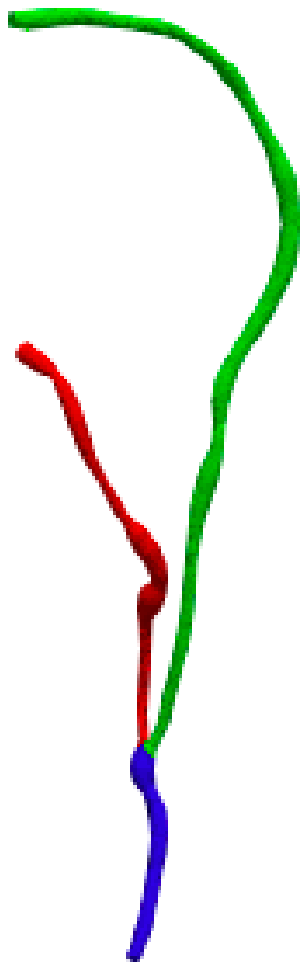}} & \multicolumn{1}{c}{\hspace*{1cm}} & \multicolumn{1}{c}{\multirow{-12}{*}{\rotatebox{90}{(c). $\Omega_c$}}} & \multicolumn{1}{c}{\includegraphics[width = 0.103\textwidth, keepaspectratio]{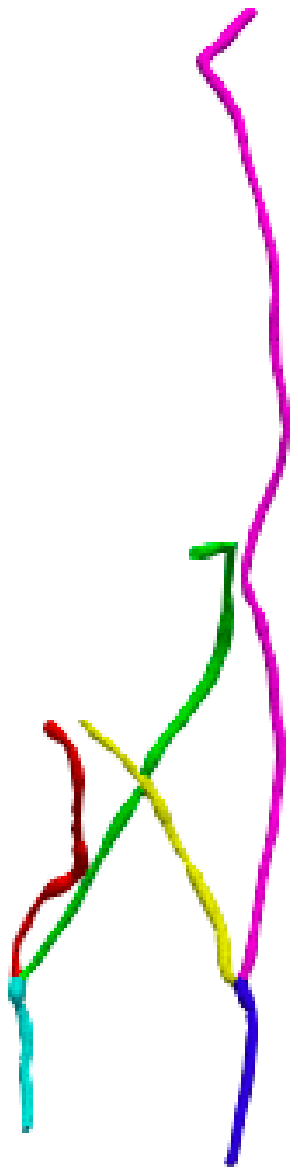}}\\
\end{tabular}
\caption{(a). Single graft connection: between right internal mammary artery (RIMA) (green) and stenosed left anterior descending artery (LAD) (red). (b). Single graft connection: between saphenuous vein (SV) (green) and stenosed first obtused marginal artery (OM1) (red). (c). Two graft connections: between right internal mammary artery (RIMA) (magenta) and stenosed left anterior descending artery (LAD) (yellow), and between saphenuous vein (SV) (green) and stenosed first obtused marginal artery (OM1) (red).}
\label{figure2}
\end{figure}

In this section, we will demonstrate application of the reduced order optimal control framework to different patient-specific geometrical models of coronary artery bypass grafts (see figure \ref{figure2}) reconstructed from CT-scan following the algorithm discussed in section \ref{sec: geometrical reconstruction} (see figure \ref{tikzfigure1}).

These numerical tests will aim at matching certain desired blood flow velocity $\bm{v_o}$ through $\mathcal{J}$ defined in equation \eqref{costFunctional}, for different parameter-induced inflow velocity scenarios\footnote{We remark that since we are using artificial target data in these numerical tests, which are carried out for methodological verification purposes, we will impose same artificial target velocity across an entire cardiovascular configuration in each numerical test. The target velocity will be imposed in parabolic shape point-wise along the centerlines. We remark also that different target velocities can be imposed in different regions of the computational domains, as is expected when matching patient-specific blood flow velocities, through the same expression \eqref{obsVelocity}.}. As observed in studies \citealp{FujiwaraEtAl1988, SankaranarayananEtAl2005} the maximum coronary flow velocity ranges between $1\ m/s$ to $4\ m/s$, we remark that one can choose any value with maximum between this range as target velocity. In this work, we choose the magnitude of desired velocity to be $0.35\ m/s \ll 4\ m/s$, that is, $v_{const} = 350\ mm/s$. This arbitrary desired blood flow velocity $\bm{v_o}$ will be distributed across the corresponding coronary artery and bypass graft through the following parabolic expression:
\begin{equation}
\label{obsVelocity}
\bm{v_o} = v_{const}\left(1 - \frac{r^2}{R^2}\right)\bm{t}_{c}.
\end{equation}
Here, for a vessel, $\bm{t}_c$ is the point-wise tangent along corresponding centerline in axial direction, $R$ is maximum radius of the vessel corresponding to points on the centerline and $r$ is the distance between mesh nodes and nearest point on the centerline. Furthermore, in these tests, the constant kinematic viscosity is $\eta = 3.6\ mm^2/s$ and the velocity at inlets will be generated by Reynolds number, $Re$, through the following expression:
\begin{equation}
\label{InflowVelocity}
\bm{v_{in}} = -\frac{\eta Re}{R_{in}}\left(1 - \frac{r^2}{R_{in}^2}\right)\bm{n}_{in}.
\end{equation}
Here, $R_{in}$ is the maximum radius of an inlet, $\bm{n}_{in}$ denotes the outward normal to the inlet and $r$ is the distance between mesh nodes and center of an inlet. Morever, these numerical tests are divided into three sub-cases, one constrained by Stokes equations and the other two constrained by Navier-Stokes equations and we will show reliability and efficiency achieved by {POD}--Galerkin in comparison to the Galerkin finite element method for each of them. Furthermore, average error for state, adjoint and control approximations, achieved by the two methods, will be calculated through the following norms:
\begin{equation}
\label{error}
\mathcal{E}_{\bm{s}} = \|\bm{s}_h\left(\bm{\mu}\right) - \bm{s}_\mathrm{N}\left(\bm{\mu}\right)\|_{S\left(\Omega\right)}, \quad \quad \mathcal{E}_{\bm{z}} = \|\bm{z}_h\left(\bm{\mu}\right) - \bm{z}_\mathrm{N}\left(\bm{\mu}\right)\|_{Z\left(\Omega\right)}, \quad \quad \mathcal{E}_{\bm{u}} = \|\bm{u}_h\left(\bm{\mu}\right) - \bm{u}_\mathrm{N}\left(\bm{\mu}\right)\|_{U\left(\Gamma_{o}\right)}.
\end{equation}
Relative error $\left(\mathcal{E}_{rel}\right)$, absolute average error $\left(\mathcal{E}_T\right)$ and absolute relative error $\left(\mathcal{E}_{T_{rel}}\right)$ will be calculated through following expressions:
\begin{equation}
\label{RelErr}
\mathcal{E}_{\bm{s}_{rel}} = \frac{\mathcal{E}_{\bm{s}}}{\|\bm{s}_h\|_{S\left(\Omega\right)}},\ 
\mathcal{E}_{\bm{z}_{rel}} = \frac{\mathcal{E}_{\bm{z}}}{\|\bm{z}_h\|_{Z\left(\Omega\right)}},\
\mathcal{E}_{\bm{u}_{rel}} = \frac{\mathcal{E}_{\bm{u}}}{\|\bm{u}_h\|_{U\left(\Gamma_o\right)}},
\end{equation}
\begin{equation}
\label{totAvgErrAndTotRelErr}
\mathcal{E}_T = \left(\mathcal{E}_{\bm{s}}^2 + \mathcal{E}_{\bm{z}}^2 + \mathcal{E}_{\bm{u}}^2\right)^{1/2},\ \mathcal{E}_{T_{rel}} = \frac{\mathcal{E}_T}{\left( \|\bm{s}_h\|^2_{ S \left(\Omega\right)} + \|\bm{u}_h\|^2_{U\left(\Gamma_o\right)} + \|\bm{z}_h\|^2_{{Z}\left(\Omega\right)} \right)^{1/2}},
\end{equation}
and the difference between Galerkin finite element and {POD}--Galerkin approximations of objective functional will be calculated as below:
\begin{equation}
\label{errorJ}
\mathcal{E}_{\mathcal{J}} = |\mathcal{J}\left(\mathbf{x}_h; \bm{\mu}\right) - \mathcal{J}\left(\mathbf{x}_\mathrm{N}; \bm{\mu}\right)|.
\end{equation}

\subsection{Stokes constrained optimal flow control: single graft connection} \label{ss: Stokes constrained optimal flow control}
In the first numerical test, we will solve a Stokes constrained optimal flow control problem over domain $\Omega_a$ (see figure \ref{figure2}(a)), to illustrate the speedup achieved by {POD}--Galerkin. We reiterate that in case of Stokes equations as state governing equations, the mathematical formulation remains the same as in previous sections of this article, with $\left(\bm{v}\left(\bm{\mu}\right)\cdot\nabla\right)\bm{v}\left(\bm{\mu}\right) = 0$ in equation \eqref{NavierStokesEquations}. In this test case, the inlets, denoted by $\Gamma_{in}$, correspond to the top openings of left anterior descending artery (LAD) and right internal mammary artery (RIMA), and the outlet, denoted by $\Gamma_o$ is the bottom opening of LAD. Moreover, in this numerical problem we aim to minimize the misfit between arbitrary desired velocity and velocity achieved through Stokes flow, for different parametrized inflow velocity scenarios.

For this purpose, as already discussed before, we consider $v_{const} = 350\ mm/s$ and consider $\bm{v}_o$ to be distributed across $\Omega_a$ through expression \eqref{obsVelocity}. The parametrized version of the velocity profile \eqref{InflowVelocity} for $\bm{\mu}= Re\in \mathcal{D} = \left[70, 80\right]$ is defined through the following expression:
\begin{equation}
\label{parametrizedInflowVelocity}
\bm{v_{in}}\left(\bm{\mu}\right) = -\frac{\eta\bm{\mu}}{R_{in}}\left(1 - \frac{r^2}{R_{in}^2}\right)\bm{n}_{in},\quad \bm{\mu}\in \mathcal{D},
\end{equation}
\noindent and is prescribed at the inlets $\Gamma_{in}$ of $\Omega_a$ through Dirichlet boundary conditions. Furthermore, as discussed in section \ref{ss: Parametrized Optimal Flow Control Problems}, no-slip boundary conditions are considered at $\Gamma_w$ and the unknown outflow boundary conditions will be quantified through the control function implemented in Neumann sense (see equations \eqref{boundaryConditions}) at $\Gamma_o$. Thus, through the parametrized optimal flow control framework applied in this numerical test, we shall automatically obtain the energy per unit mass required at the outlets, to match Stokes velocity with the desired velocity corresponding to different tunings of the inflow velocity.

\begin{figure}
\centering
\begin{minipage}{7cm}
\begin{center}
\begin{tabular}{lll}
\hline \hline
\multicolumn{1}{|c|}{Mesh size} & \multicolumn{1}{c|}{$42354$}\\
\hline
\multicolumn{1}{|c|}{$\mathcal{D}$} & \multicolumn{1}{c|}{$\left[70, 80\right]$}\\
\hline
\multicolumn{1}{|c|}{$|\Lambda |$} & \multicolumn{1}{c|}{$100$}\\
\hline
\multicolumn{1}{|c|}{Offline phase} & \multicolumn{1}{c|}{$4191.32$ seconds}\\
\hline
\multicolumn{1}{|c|}{Online phase} & \multicolumn{1}{c|}{$6.4$ seconds}\\
\hline \hline
\end{tabular}
\caption*{(a).}
\end{center}
\end{minipage}
\begin{minipage}{7cm}
\begin{center}
\includegraphics[width = 0.9\textwidth, keepaspectratio]{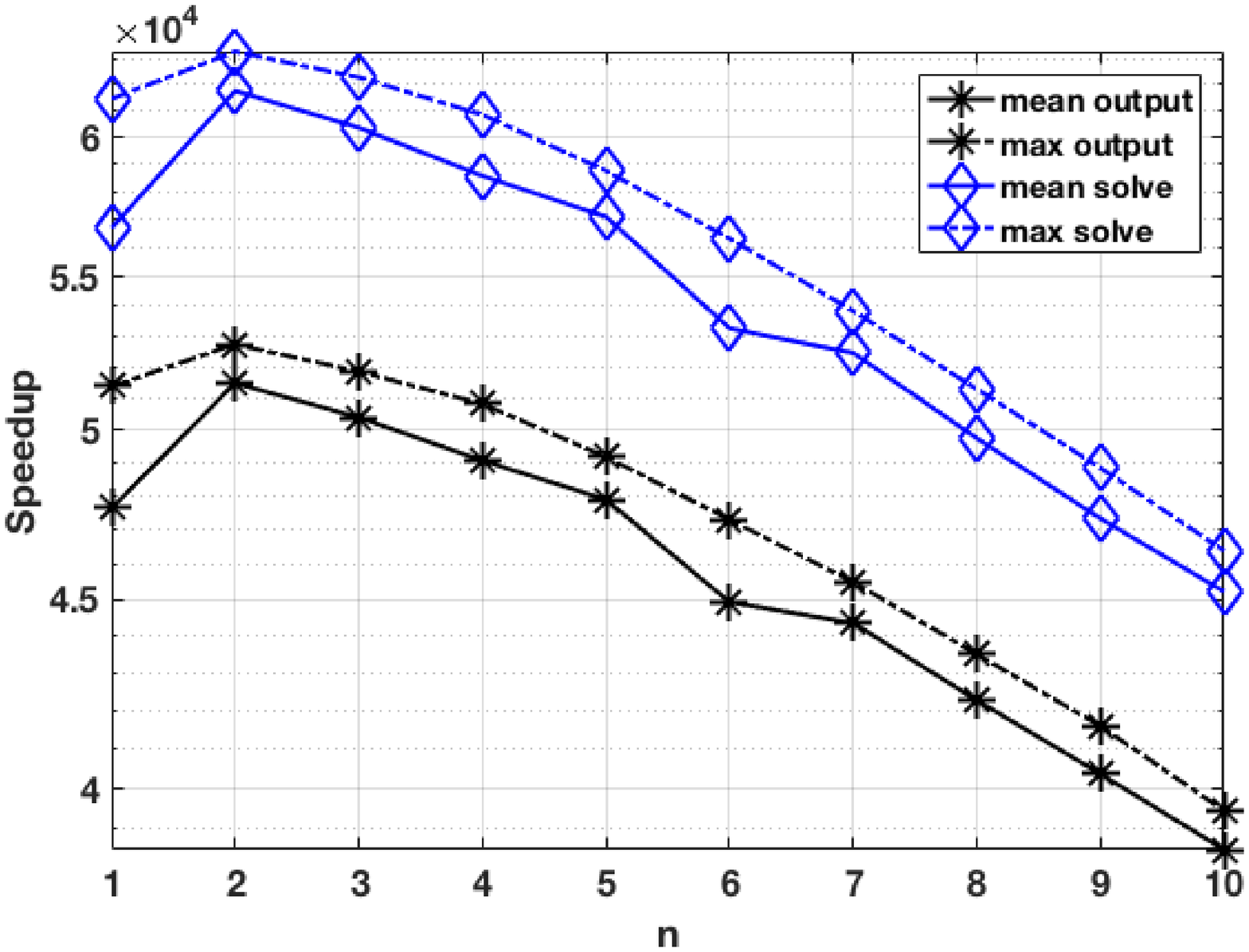}
\caption*{(b).}
\end{center}
\end{minipage}
\caption{(a). Table demonstrating computational performance of {POD}--Galerkin for Stokes constrained optimal control problem. (b). Mean and maximum speedups for solution and objective functiona| achieved by {POD}--Galerkin for Stokes constrained optimal control problem.}
\label{ComputationalPerformanceAndSpeedupForStokes}
\end{figure}

At the continuous level, we consider the following spaces for velocity, pressure and control, respectively:
\begin{equation*}
V\left(\Omega_a\right) =  H^1_{\Gamma_{in}\cup \Gamma_w}\left(\Omega_a\right) = \left[ \bm{v}\in \left[H^1\left(\Omega_a\right)\right]^3 : \bm{v}|_{\Gamma_{in}} = \bm{v_{in}}\ \text{and}\ \bm{v}|_{\Gamma_w} = \bm{0} \right],\quad P\left(\Omega_a\right)= L^2\left(\Omega_a\right),
\end{equation*}
and
\begin{equation*}
U\left(\Gamma_o\right)= \left[ L^2\left(\Gamma_o\right)\right]^3.
\end{equation*}
To construct the reduced order spaces, we sample $\bm{\Lambda}\subset\mathcal{D}$ for $100$ parameter values distributed uniformly. The offline phase, then, starts with snapshots collection for these parameter values by solving {\it truth problem}, that is the algebraic system \eqref{algebraicFEMatrixForm}, through Galerkin finite element method. At the finite element level, we utilize stable $\mathbb{P}2-\mathbb{P}1$ pair for both state and adjoint velocity and pressure and $\mathbb{P}2$ for control and therefore, fulfilling both Brezzi's and LBB inf-sup condition at this level. This phase ends with the construction of {POD} basis using $\mathrm{N}_{max} = 10$ eigenvalue-eigenvector pairs that keep relative energy greater than $1 - \epsilon_{tol},\ 0 < \epsilon_{tol}\ll 1$. We report the CPU time required for this phase to be $4191.32$ seconds (see table reported in figure \ref{ComputationalPerformanceAndSpeedupForStokes}(a)). The online phase comprises of determining reduced order coefficients through Galerkin projection onto the full order solution manifold and solving the reduced order problem \ref{pb: rom}. This phase takes only $6.4$ seconds and is repeated for different values of $\bm{\mu}$ chosen from $\mathcal{D}$.

We illustrate average and maximum speedups attained in this case through figure \ref{ComputationalPerformanceAndSpeedupForStokes}(b). An average speedup of $\mathcal{O}\left( 10^4\right)$ is achieved for both, the output objective functional and solution, as $n$ goes from $1$ to $\mathrm{N}_{max}$, with the maximum speedups falling in the same range.

\subsection{Navier-Stokes constrained optimal flow control problem}
In the numerical examples henceforth, we will solve optimal control problems with state governed by incompressible steady Navier-Stokes equations for two types of geometries, that is with single and double grafts connections\footnote{With double grafts connections/two grafts connections, we mean two separate grafts connected to two separate diseased arteries.}.

Furthermore, for single graft connection, we will report two sub-cases and will compare the computational performance of Galerkin finite element and {POD}--Galerkin in these sub-cases. For double grafts connections, we will consider combination of the geometries used in the sub-cases of single graft connection problem (see figure \ref{figure2}). Thus, the double grafts connections problem has increased computational complexity as the problem is to be solved over a larger computational domain and with increased number of parameters.

\subsubsection{Single graft connection: case I} \label{sss: Single graft connection: case I}

In the first case, we consider domain $\Omega_a$ (figure \ref{figure2}(a)) in which graft connection is between right internal mammary artery (RIMA) and stenosed left anterior descending artery (LAD). We take $\bm{\mu} = Re\in \mathcal{D} = \left[ 70, 80\right]$ and consider the same solution spaces as considered in Stokes constrained problem \ref{ss: Stokes constrained optimal flow control}. Moreover, $v_{const} = 350\ mm/s$ such that $\bm{v_o}\in \left[L^2\left(\Omega_a\right)\right]^3$. We solve a full order problem for $\mu = 80$ using $\mathbb{P}2-\mathbb{P}1$ for velocity and pressure and $\mathbb{P}2$ for control and attain $433288$ degrees of freedom $\left(\mathcal{N}\right)$. One such simulation requires $1213.3$ seconds.

For the reduced order problem, we use $\mathrm{N}_{max} = 6$ and attain $\mathrm{N} = 79$ reduced bases in total. These reduced bases are classified as $48$ bases altogether for velocity and adjoint velocity, $24$ bases altogether for pressure and adjoint pressure, $6$ bases for control and $1$ basis for lifting for non-homogenuous Dirichlet inflow condition.

\begin{figure}[h]
\centering
\begin{tabular}{lllll}
\multicolumn{1}{c}{} & \multicolumn{1}{c}{{\bf Case I}} & \multicolumn{1}{c}{\hspace*{1.3cm}} & \multicolumn{1}{c}{{\bf Case II}} \\
\multicolumn{1}{c}{\multirow{-12}{*}{\rotatebox{90}{(a).}}} & \multicolumn{1}{c}{\includegraphics[width = 0.25\textwidth , keepaspectratio]{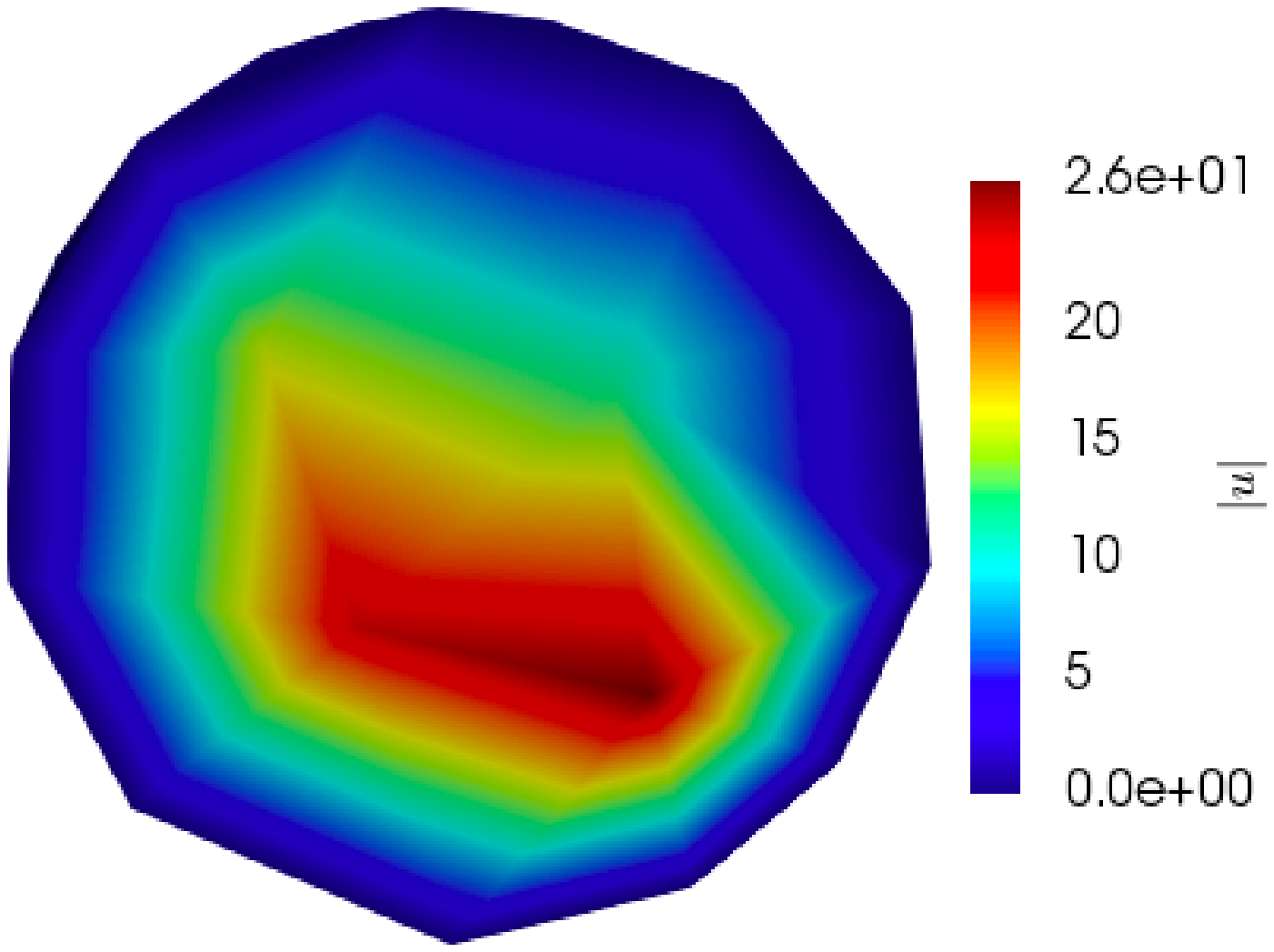}} & \multicolumn{1}{c}{\hspace*{1.3cm}} & \multicolumn{1}{c}{\includegraphics[width = 0.3\textwidth , keepaspectratio]{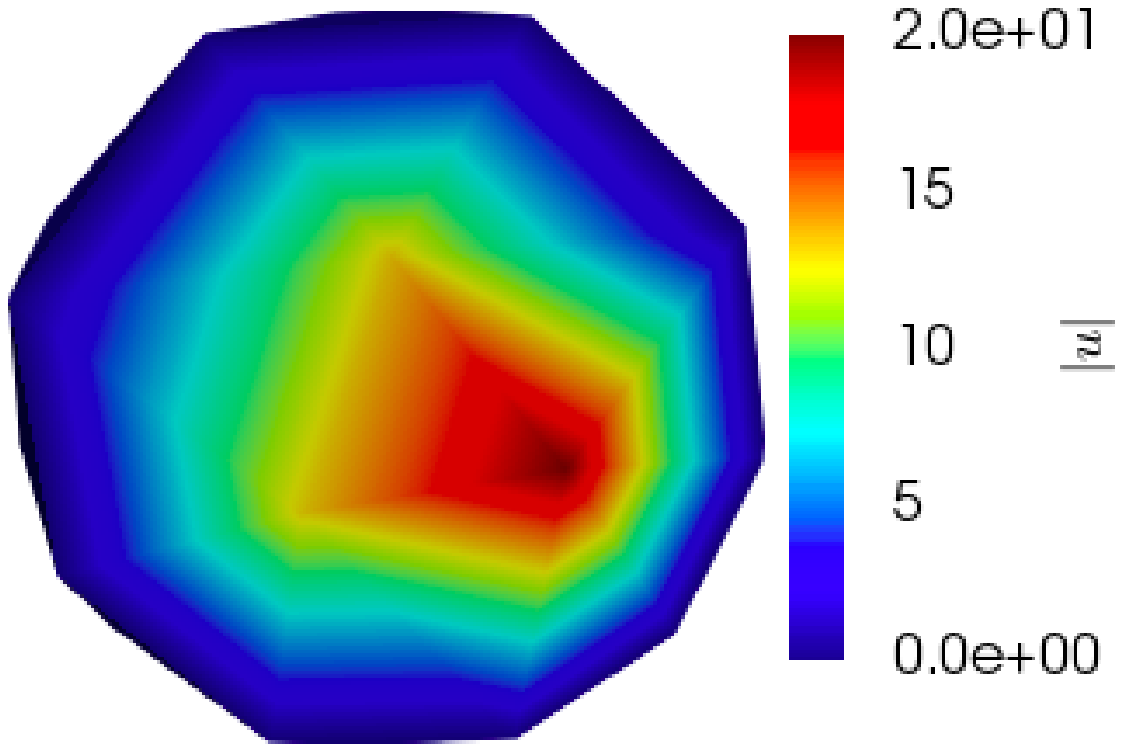}}\\
\multicolumn{1}{c}{\multirow{-12}{*}{\rotatebox{90}{(b).}}} & \multicolumn{1}{c}{\includegraphics[width = 0.4\textwidth , keepaspectratio]{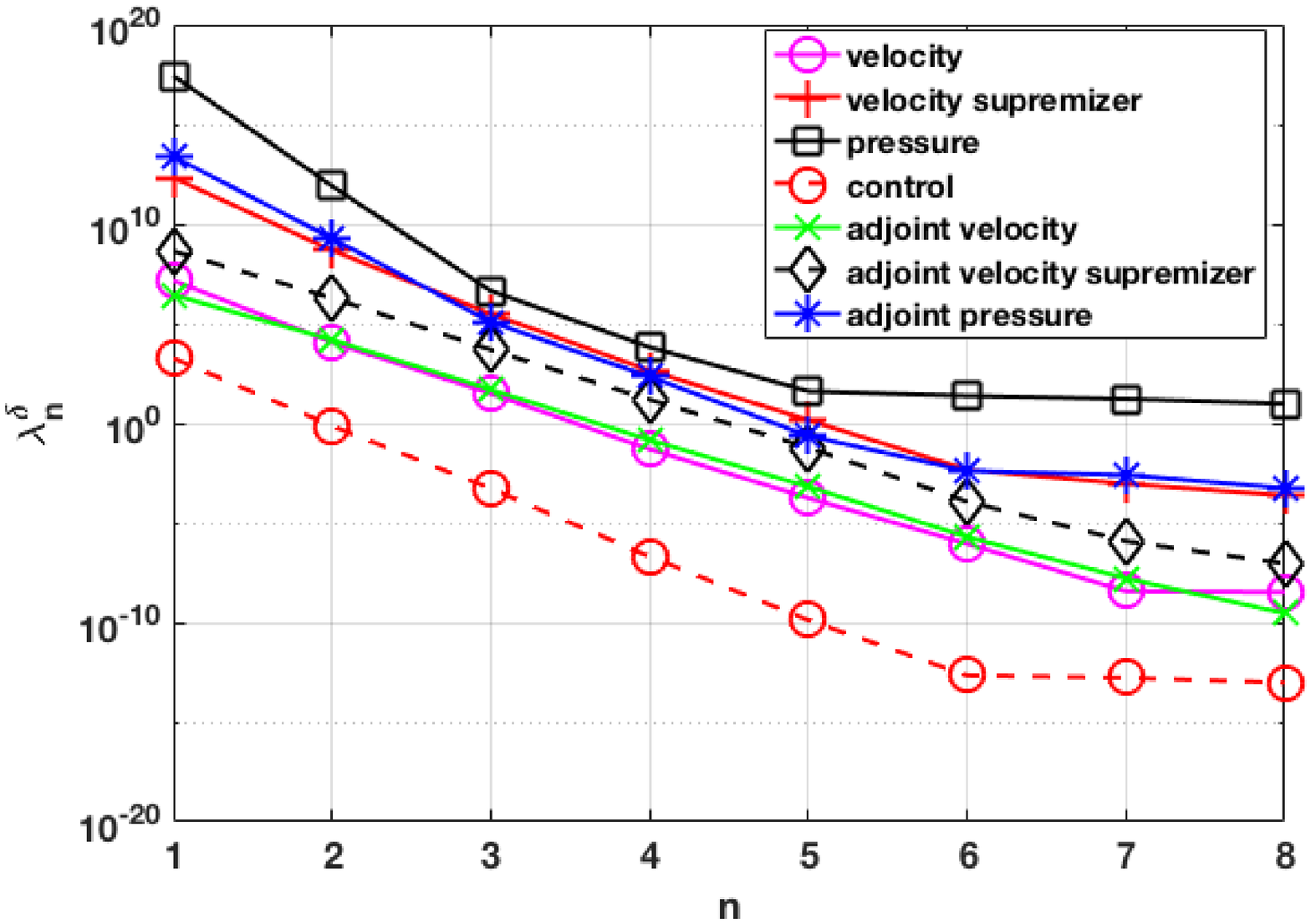}} & \multicolumn{1}{c}{\hspace*{1.3cm}} & \multicolumn{1}{c}{\includegraphics[width = 0.4\textwidth , keepaspectratio]{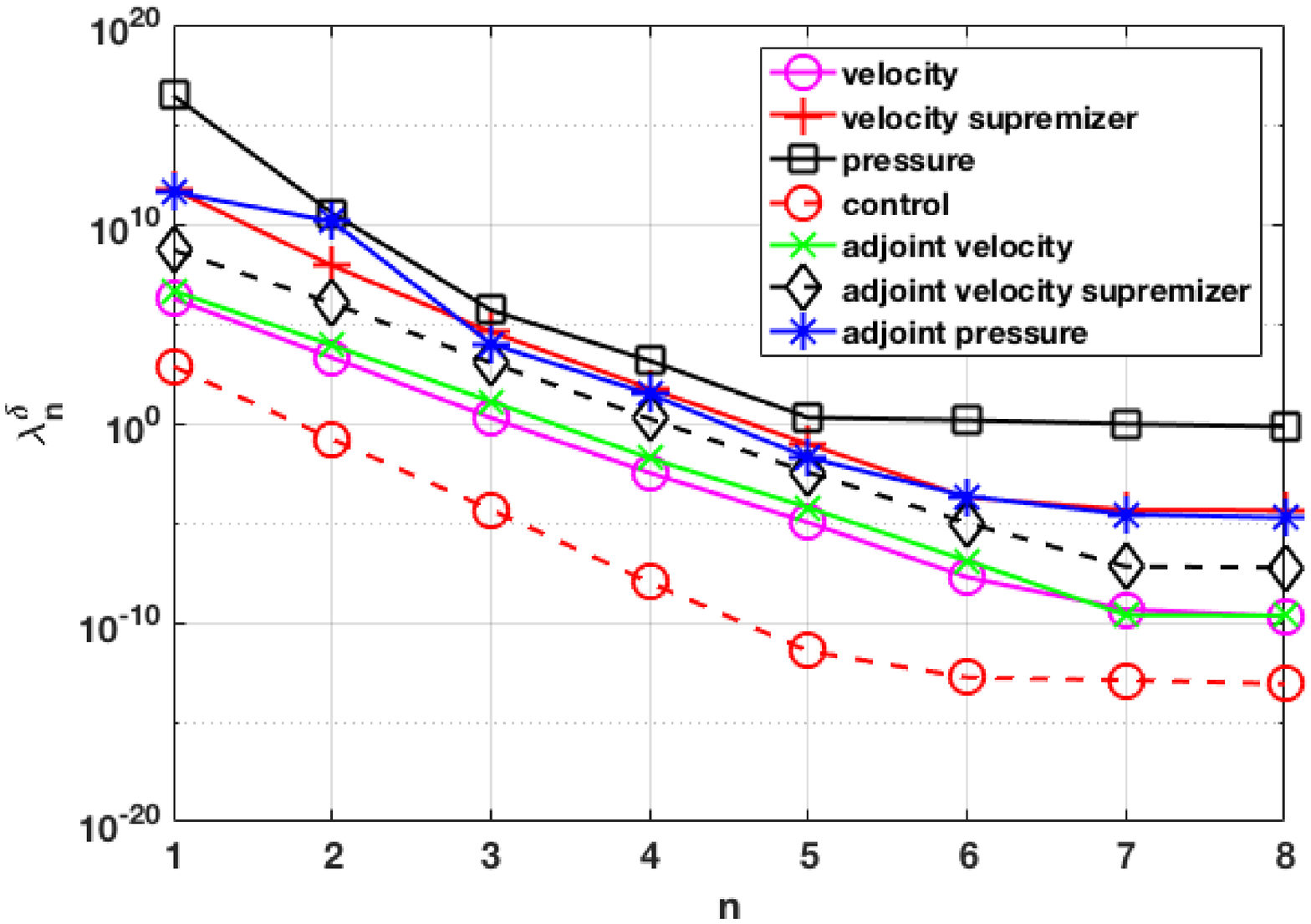}}\\
\end{tabular}
\caption{(a). Boundary control magnitude $\left( mm^2/s^2\right)$. (b). Eigenvalues' reduction.}
\label{Boundary control magnitudes and eigenvalues' reductions for case I and case II.}
\end{figure}
\noindent Figure \ref{Boundary control magnitudes and eigenvalues' reductions for case I and case II.}(b) shows that $99.99\%$ energy of full order solution manifold is captured in $\mathrm{N}_{max}$ eigenvectors. Thus in this case, we retain same accuracy and reliability as that of a full order model, which is visible through same velocity approximations achieved by Galerkin finite element that utilizes $\mathcal{O}\left( 10^5\right)$ dofs and {POD}--Galerkin approximations that use only $\mathcal{O}\left( 10^1\right)$ basis (see figures \ref{LADOM1velApprox.}(a) and \ref{LADOM1velApprox.}(b)). Moreover, we show magnitude of boundary control in figure \ref{Boundary control magnitudes and eigenvalues' reductions for case I and case II.}(a).

\begin{figure}[h]
\centering
\begin{tabular}{lllll}
\multicolumn{1}{c}{} & \multicolumn{1}{c}{{\bf Case I}} & \multicolumn{1}{c}{\hspace*{1cm}} & \multicolumn{1}{c}{{\bf Case II}} \\
\multicolumn{1}{c}{\multirow{-12}{*}{\rotatebox{90}{(a).}}} & \multicolumn{1}{c}{\includegraphics[width = 0.4\textwidth , keepaspectratio]{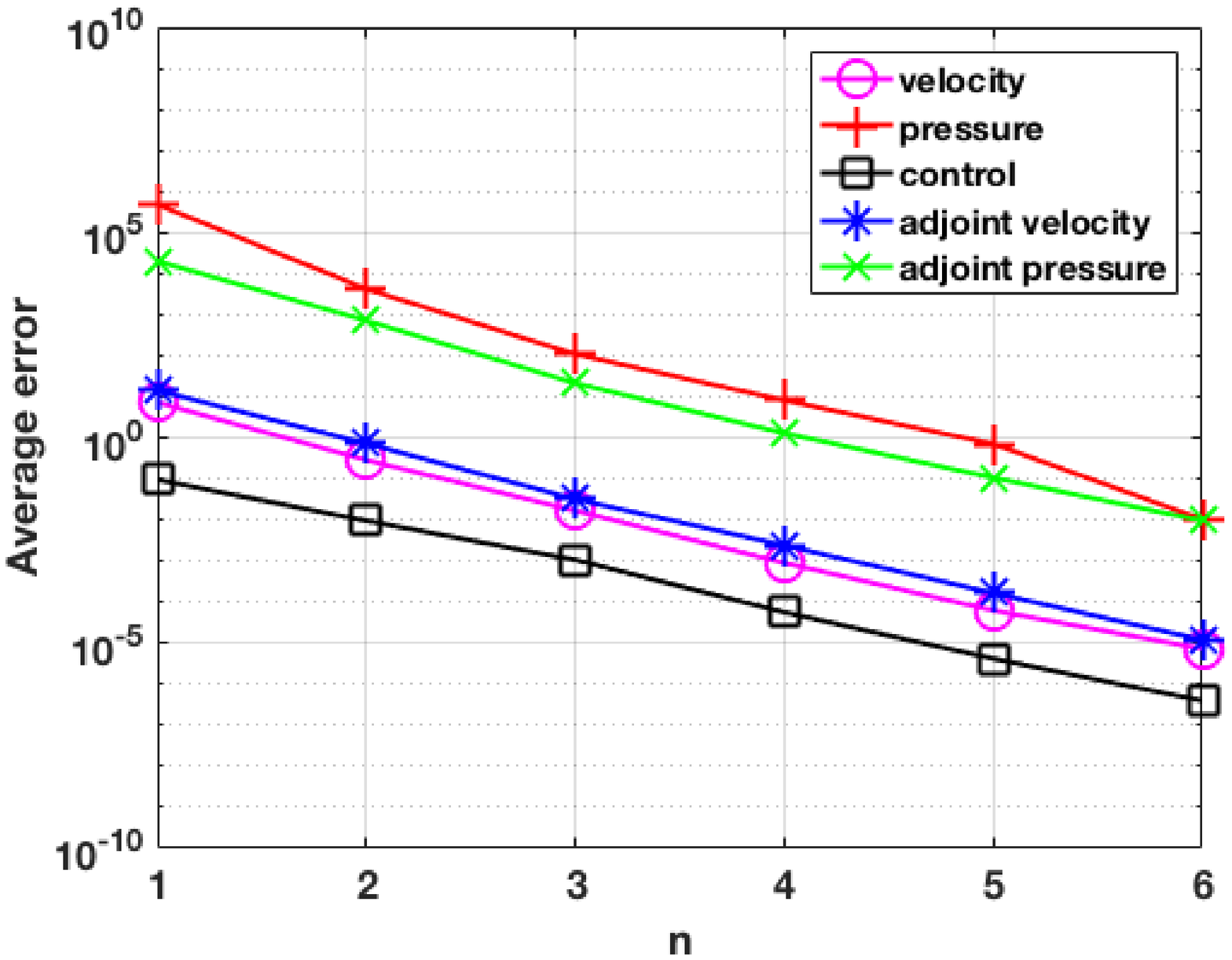}} & \multicolumn{1}{c}{\hspace*{1cm}} & \multicolumn{1}{c}{\includegraphics[width = 0.4\textwidth , keepaspectratio]{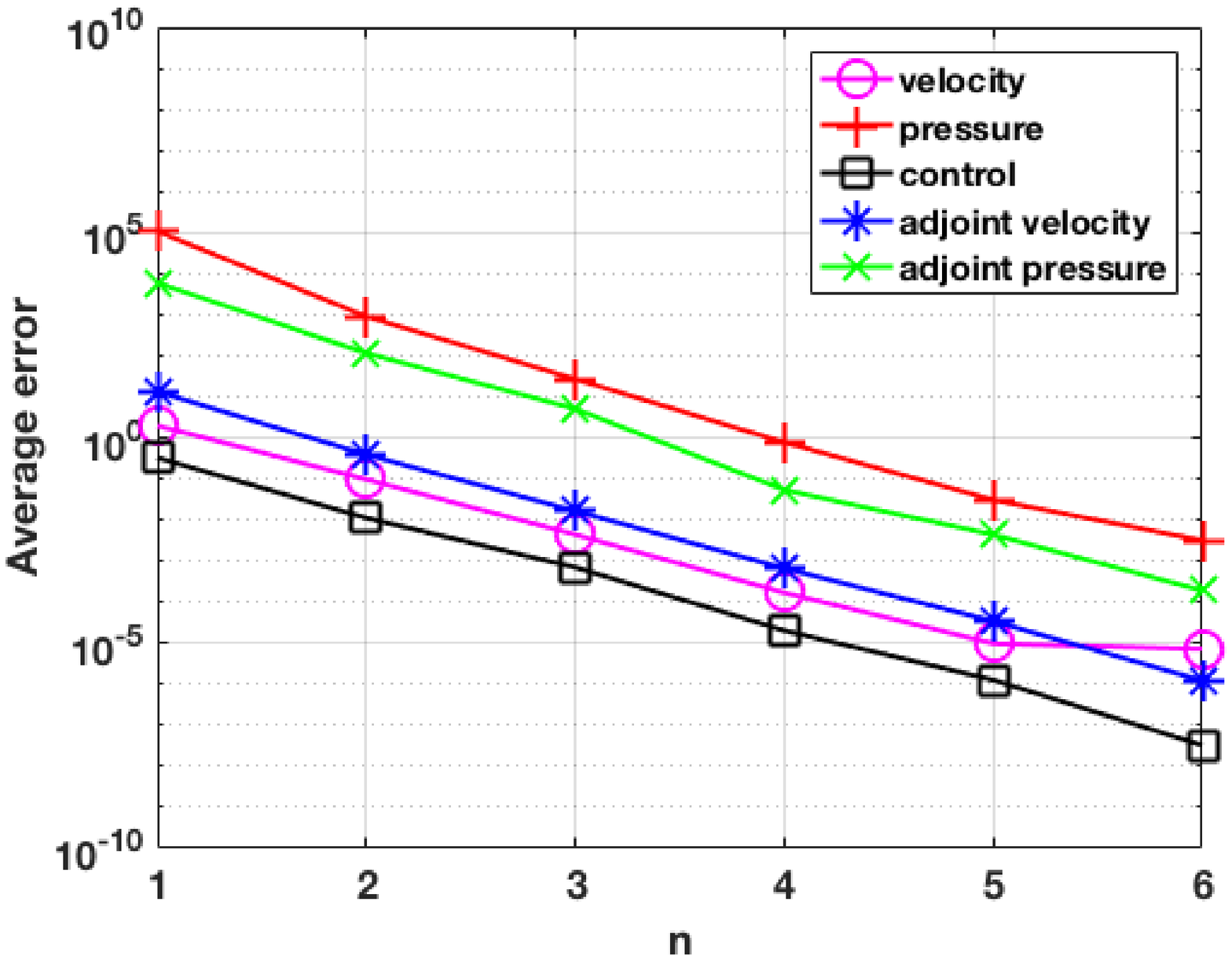}}\\
\multicolumn{1}{c}{\multirow{-12}{*}{\rotatebox{90}{(b).}}} & \multicolumn{1}{c}{\includegraphics[width = 0.4\textwidth , keepaspectratio]{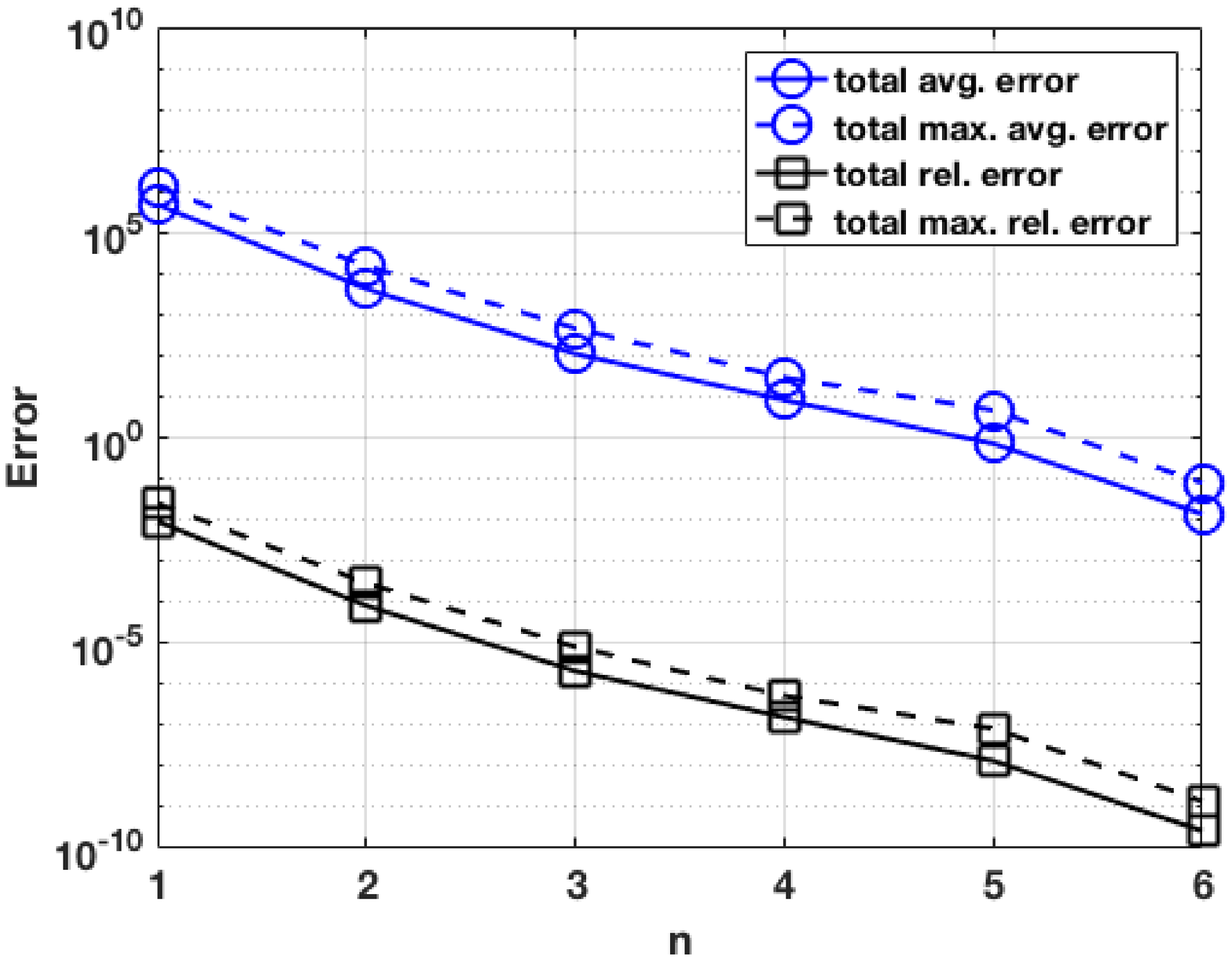}} & \multicolumn{1}{c}{\hspace*{1cm}} & \multicolumn{1}{c}{\includegraphics[width = 0.4\textwidth , keepaspectratio]{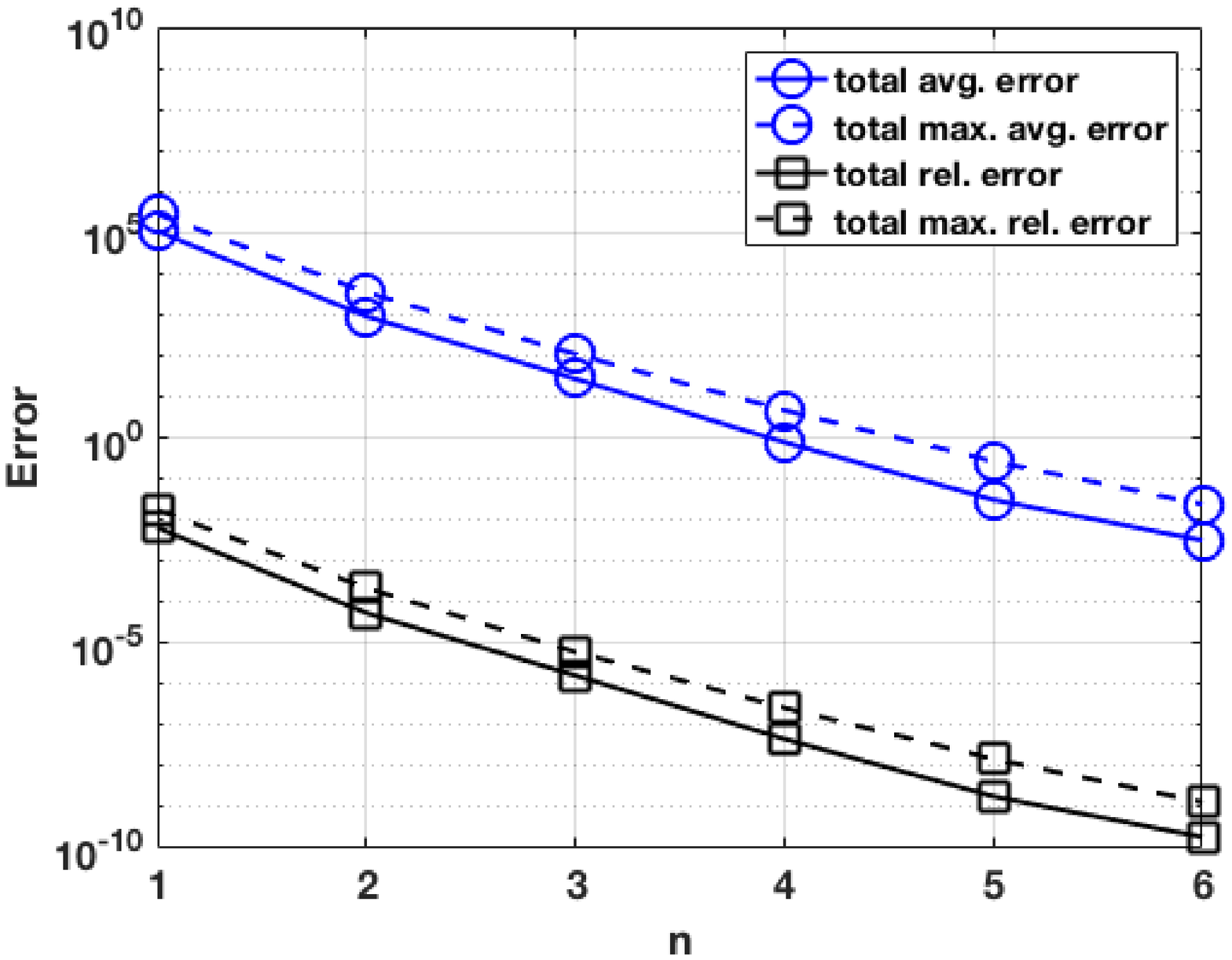}}\\
\multicolumn{1}{c}{\multirow{-12}{*}{\rotatebox{90}{(c).}}} & \multicolumn{1}{c}{\includegraphics[width = 0.4\textwidth , keepaspectratio]{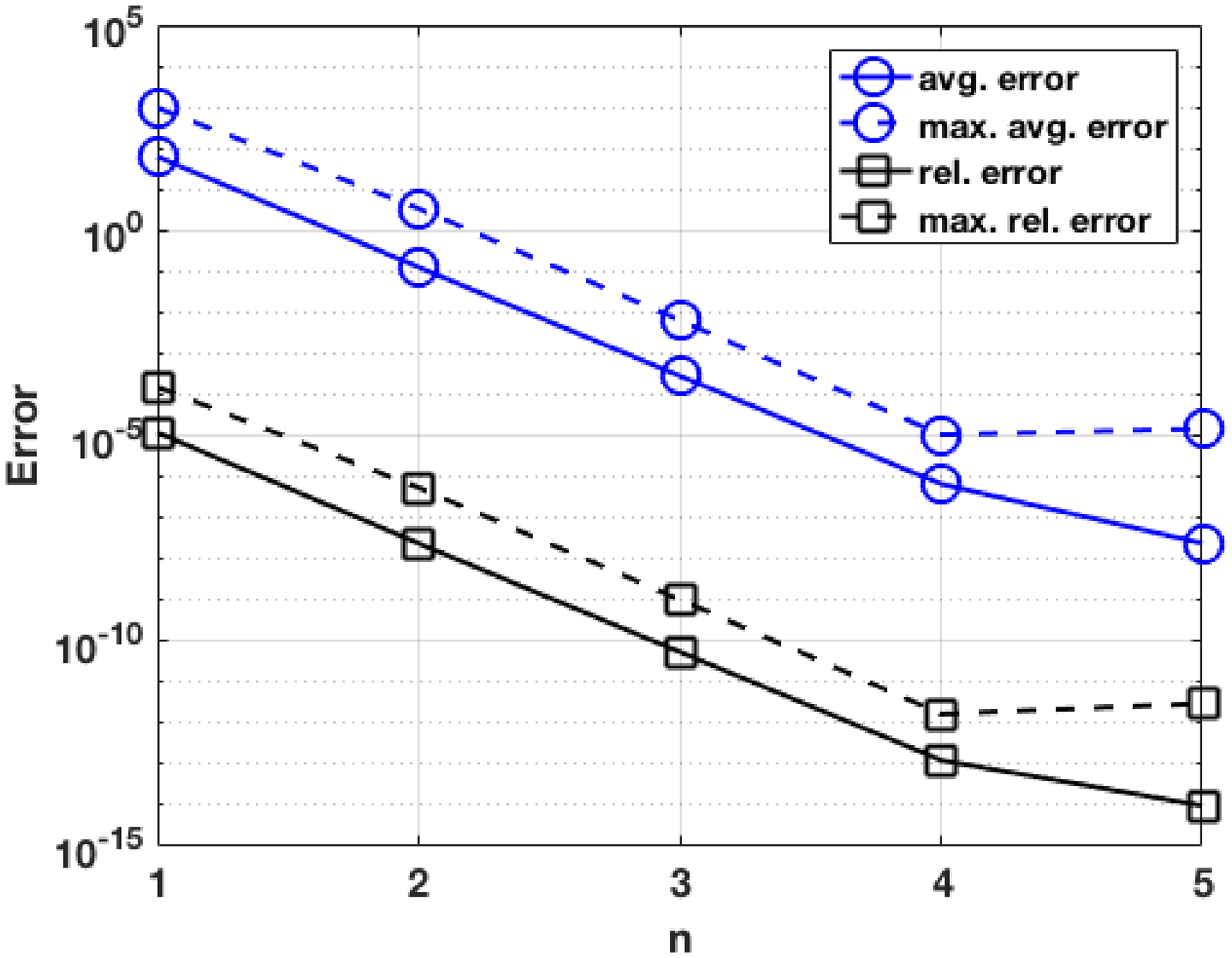}} & \multicolumn{1}{c}{\hspace*{1cm}} & \multicolumn{1}{c}{\includegraphics[width = 0.4\textwidth , keepaspectratio]{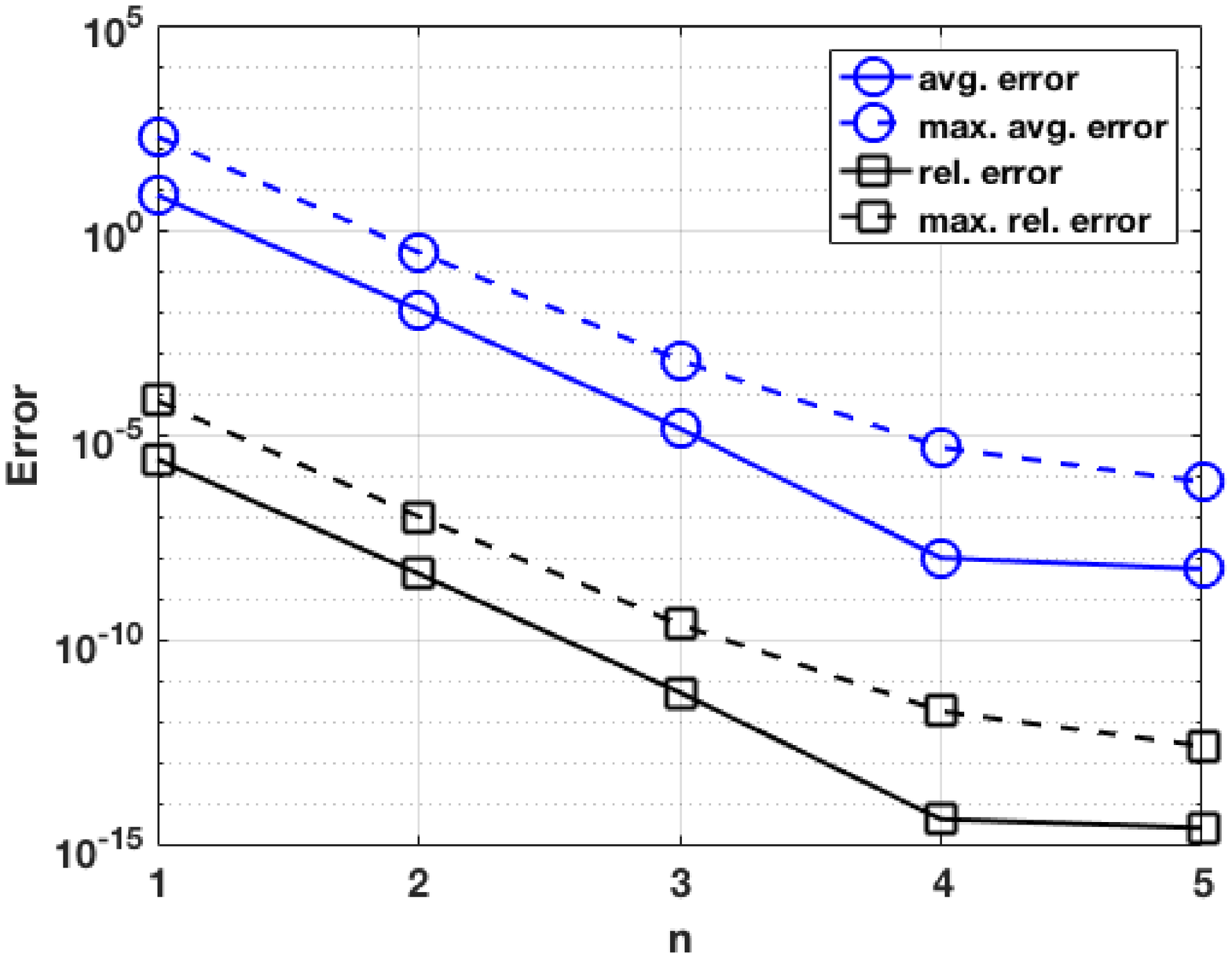}}\\
\end{tabular}
\caption{(a). Average error between FE and {POD}--Galerkin approx. of $\bm{\delta} = \bm{v}, p, \bm{u}, \bm{w}, q$. (b). Total error between Galerkin FE and {POD}--Galerkin approximations. (c). Error between FE and {POD}--Galerkin reduction of $\mathcal{J}$.}
\label{Errors for case I and case II}
\end{figure}
\begin{figure}[h]
\centering
\begin{tabular}{lllllllll}
\hline
\multicolumn{2}{|c|}{} & \multicolumn{2}{c|}{} & \multicolumn{2}{c|}{} & \multicolumn{2}{c|}{}\\
\multicolumn{1}{|c}{\multirow{-12}{*}{\rotatebox{90}{(a).}}} & \multicolumn{1}{c|}{\includegraphics[width = 0.129\textwidth, keepaspectratio]{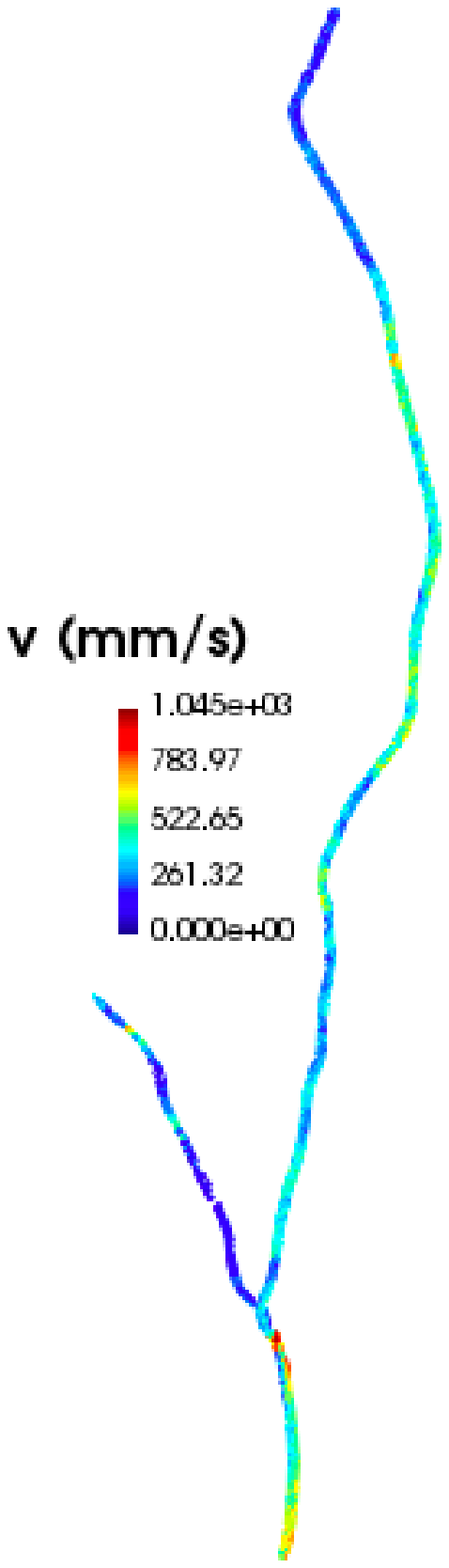}} & \multicolumn{1}{c}{\multirow{-12}{*}{\rotatebox{90}{(b).}}} & \multicolumn{1}{c|}{\includegraphics[width = 0.119\textwidth, keepaspectratio]{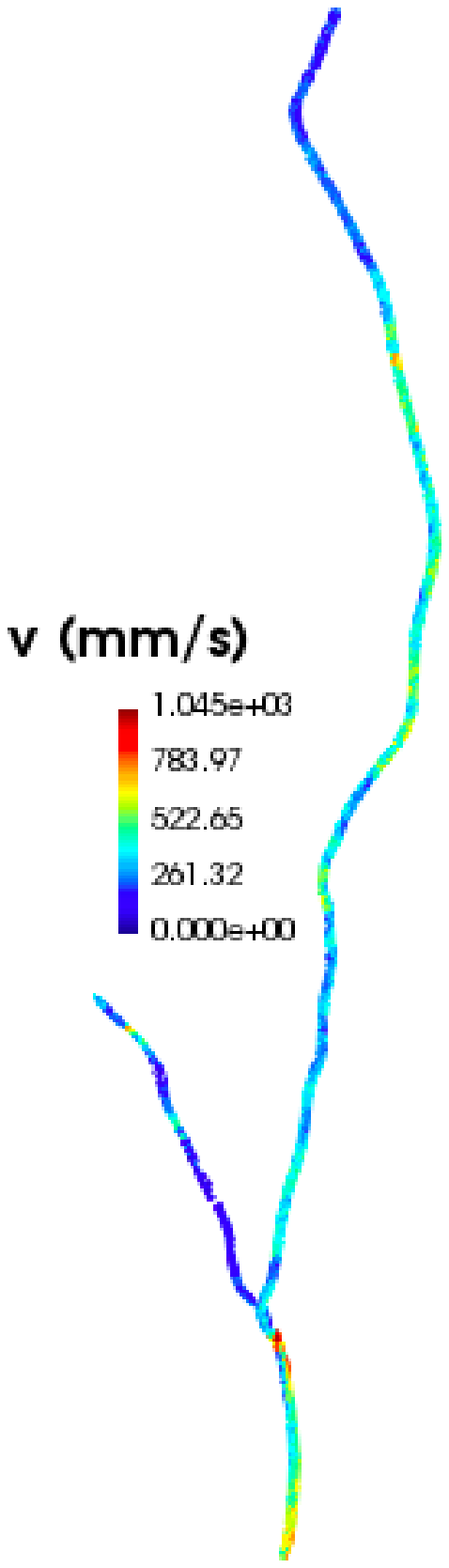}} & \multicolumn{1}{c}{\multirow{-12}{*}{\rotatebox{90}{(c).}}} & \multicolumn{1}{c|}{\includegraphics[width = 0.129\textwidth, keepaspectratio]{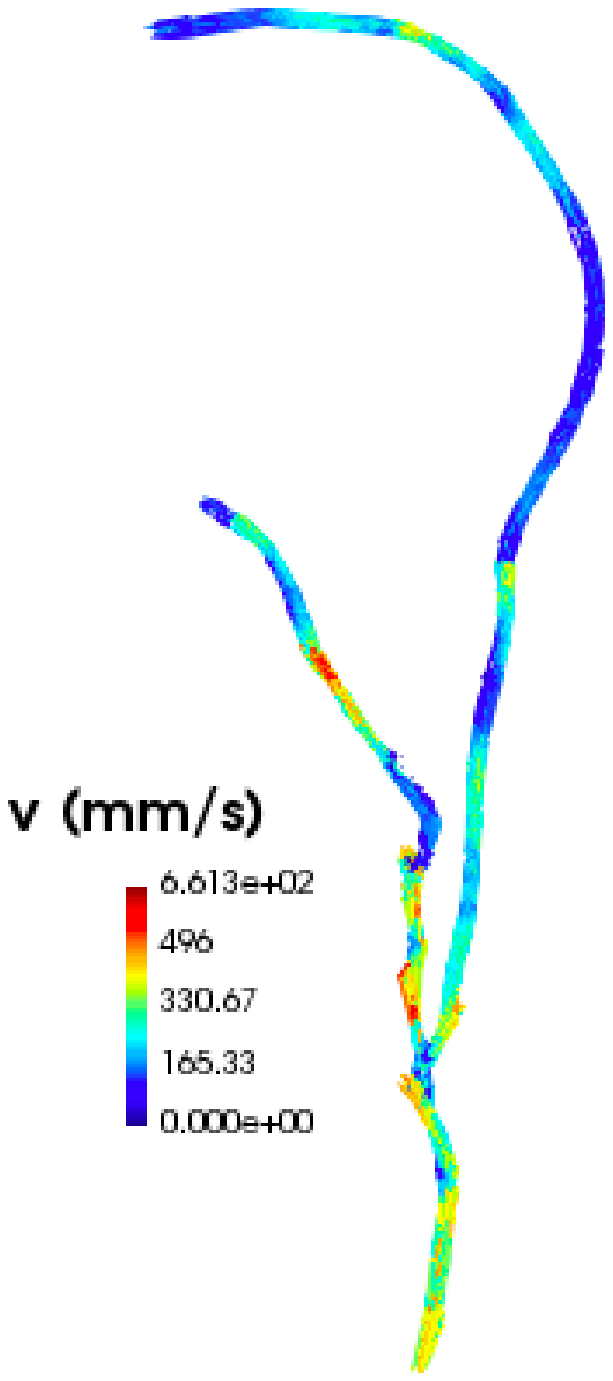}} & \multicolumn{1}{c}{\multirow{-12}{*}{\rotatebox{90}{(d).}}} & \multicolumn{1}{c|}{\includegraphics[width = 0.119\textwidth, keepaspectratio]{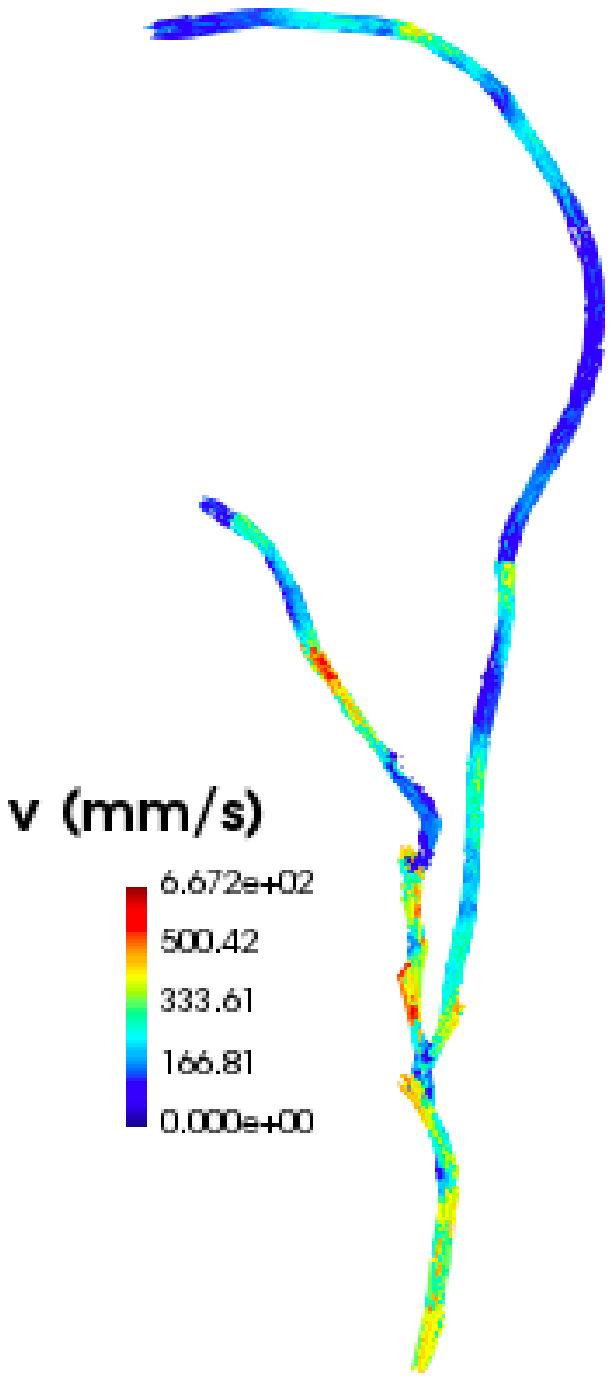}}\\
\multicolumn{2}{|c|}{} & \multicolumn{2}{c|}{} & \multicolumn{2}{c|}{} & \multicolumn{2}{c|}{}\\
\hline
\multicolumn{2}{|c|}{(a). Case I: FE approx.} & \multicolumn{2}{c|}{(b). Case I: ROM approx.} & \multicolumn{2}{c|}{(c). Case II: FE approx.} & \multicolumn{2}{c|}{(d). Case II: ROM approx.}\\
\hline

\end{tabular}
\caption{Simulation velocity $\bm{v}$}
\label{LADOM1velApprox.}
\end{figure}
 In figure \ref{Errors for case I and case II}(a), we report the error plots $\mathcal{E}_{\bm{\delta}}$ for $\bm{\delta} = \bm{v}, p, \bm{u}, \bm{w}, q$, calculated using equation \eqref{error}. $\mathcal{E}_{\bm{v}}$ decreases from $10^1$ to $10^{-5}$ approximately, $\mathcal{E}_{p}$ decreases from $10^6$ to $10^{-2}$ and $\mathcal{E}_{\bm{u}}$ decreases to $10^{-6}$ as $n$ goes from $1$ to $6$. Error reduction for adjoint variables' approximations follows similar pattern as the error reduction for state variables. Absolute error and absolute relative error are calculated through \eqref{RelErr} and \eqref{totAvgErrAndTotRelErr} and are reported in figure \ref{Errors for case I and case II}(b). $\mathcal{E}_T$ reduces from approximately $10^6$ to $10^{-2}$ for $n = 1, \cdots , \mathrm{N}_{max}$ with same behavior for maximum $\mathcal{E}_T$. A reduction from $10^{-2}$ to $10^{-10}$ is achieved for both $\mathcal{E}_{T_{rel}}$ and maximum $\mathcal{E}_{T_{rel}}$ in this case.

Figure \ref{Errors for case I and case II}(c) shows the difference between $\mathcal{J}$ computed through full order and reduced order methods. It is calculated using equation \eqref{errorJ} and a reduction from approximately $ 10^6$ to $10^{-2}$ is observed with similar behavior for relative error as $n$ increases.

\subsubsection{Single graft connection: case II}

In the second case, we consider a similar problem as case I, section \ref{sss: Single graft connection: case I} with a different coronary artery bypass graft and with different parameter range. We label the computational domain as $\Omega_b$ and take it to be saphenuous vein grafted to stenosed first obtuse marginal artery (OM1) (figure \ref{figure2}(b)). Moreover, $\bm{\mu} = Re\in \mathcal{D} = \left[ 45, 50\right]$, imposed at inlets, that are the top openings of SV and OM1, through equation \eqref{parametrizedInflowVelocity}. In this case, the observation velocity is defined across entire $\Omega_b$ by equation \eqref{obsVelocity} with $v_{const} = 350\ mm/s$. The solution spaces at the continuous level remain the same as case I, while being considered over $\Omega_b$ and we use $\mathbb{P}2-\mathbb{P}1-\mathbb{P}2$ for velocity, pressure and control respectively at the finite element level. In this case, for $\mu = 50$, we report Galerkin finite element degrees of freedom $\left(\mathcal{N}\right)$ to be $280274$, whereas, for reduced order model we report to retain $99.9\%$ energy of full order solution spaces for $\mathrm{N}_{max} = 6$ (see figure \ref{Boundary control magnitudes and eigenvalues' reductions for case I and case II.}(b)). Thus, same velocity approximation is achieved by Galerkin finite element and {POD}--Galerkin methods as shown in figures \ref{LADOM1velApprox.}(c) and \ref{LADOM1velApprox.}(d) respectively. Control magnitude in this case is illustrated in figure \ref{Boundary control magnitudes and eigenvalues' reductions for case I and case II.}(a) and is restricted to the outlet only.\\
With $\mathrm{N}_{max} = 6$, we attain $24$ reduced bases for velocity, $24$ for adjoint velocity, $1$ for lifting function for Dirichlet boundary conditions, $6$ for control, $12$ for pressure and $12$ for adjoint pressure. Therefore, the reduced state and adjoint spaces are spanned by $73$ {POD} bases and achieve a decrease to approximately $10^{-5}$ from $10^{1}$ in $\mathcal{E}_{\bm{v}}$. $\mathcal{E}_p$ decreases from $10^5$ to approximately $10^{-3}$ (see figure \ref{Errors for case I and case II}(a)) and $\left(\mathcal{E}_T\right)$ decreases to approximately $10^{-3}$ from $10^{5}$ with same order of reduction observed for total relative error $\left(\mathcal{E}_{T_{rel}}\right)$ (see figure \ref{Errors for case I and case II}(b)). The difference between Galerkin finite element and {POD}--Galerkin approximations of $\mathcal{J}$ is observed to be approximately $10^1$ at $n = 1$ and is decreased to $10^{-8}$ at $n = 5$ (see figure \ref{Errors for case I and case II}(c)).

We report the computational performance of Galerkin finite element and {POD}--Galerkin methods for case I and II in table \ref{tableForComputationalPerfomanceInCaseIAndCaseII}. Both cases have $\mathcal{O}\left(10^4\right)$ mesh size and $\mathcal{O}\left(10^5\right)$ degrees of freedom for full order solution spaces. Thanks to {POD}--Galerkin, the degrees of freedom are reduced to about $\mathcal{O}\left(10^2\right)$ with sufficient reduction in error, as discussed above. Moreover, full order simulations for single valued $\bm{\mu}$, for example, for $Re = 80$ and $Re = 50$ in case I and case II, respectively take $1214.3$ seconds and $634$ seconds. Whereas, online phases in these cases for the mentioned values of parameters, take $109.3$ seconds and $118$ seconds respectively. For different parameter-dependent cases, one has to repeat the full order simulation from scratch if dealing only with high order methods, whereas, in {POD}--Galerkin, we need to repeat only the online phase which requires much lower computational time. However, reduction in computational cost comes at a price of performing an offline phase which is expensive, as reported in table \ref{tableForComputationalPerfomanceInCaseIAndCaseII}. For cardiovascular problems, it is inevitable to consider parameter-dependent scenarios and since in reduced order methods, one has to perform offline phase only once, we consider implementation of reduced order methods to be reasonably inexpensive as compared to full order methods.
\begin{table}
\centering
\begin{tabular}{|c|c|c|c|c|c|ccc|}
\hline \hline
\multirow{2}{*}{} & \multirow{2}{*}{Mesh size}  & \multirow{2}{*}{FE dofs $\left(\mathcal{N}\right)$}  & \multirow{2}{*}{No. of RB $\left(\mathrm{N}\right)$} & \multirow{2}{*}{$\mathcal{D}$} & \multirow{2}{*}{$| \Lambda |$} & \multicolumn{3}{c|}{Computational time (seconds)}\\\cline{7-9}
& & & & & & \multicolumn{1}{c|}{Galerkin} & \multicolumn{1}{c|}{ROM} & \multicolumn{1}{c|}{{ROM}} \\
& & & & & & \multicolumn{1}{c|}{FEM} & \multicolumn{1}{c|}{offline phase} & \multicolumn{1}{c|}{online phase} \\\hline 
$\Omega_a$ & $42354$ &  $433288$ & $79$ & $\left[ 70, 80\right]$ & $100$ & $1214.3$ & $16825.9$ & $109.3$\\
\hline
$\Omega_b$ & $27398$ &  $280274$ & $79$ & $\left[ 45, 50\right]$ & $100$ & $634$ & $12106.8$ & $118$\\
\hline \hline

\end{tabular}
\caption{Computational performance of the proposed method for case I and case II.}
\label{tableForComputationalPerfomanceInCaseIAndCaseII}
\end{table}

\subsubsection{Two grafts connections}

\begin{figure}[h]
\centering
\begin{tabular}{llllll}
\hline
\multicolumn{2}{|c|}{} & \multicolumn{2}{c|}{} & \multicolumn{2}{c|}{}\\
\multicolumn{1}{|c}{\multirow{-12}{*}{\rotatebox{90}{(a).}}} & \multicolumn{1}{c|}{\includegraphics[width = 0.18\textwidth, keepaspectratio]{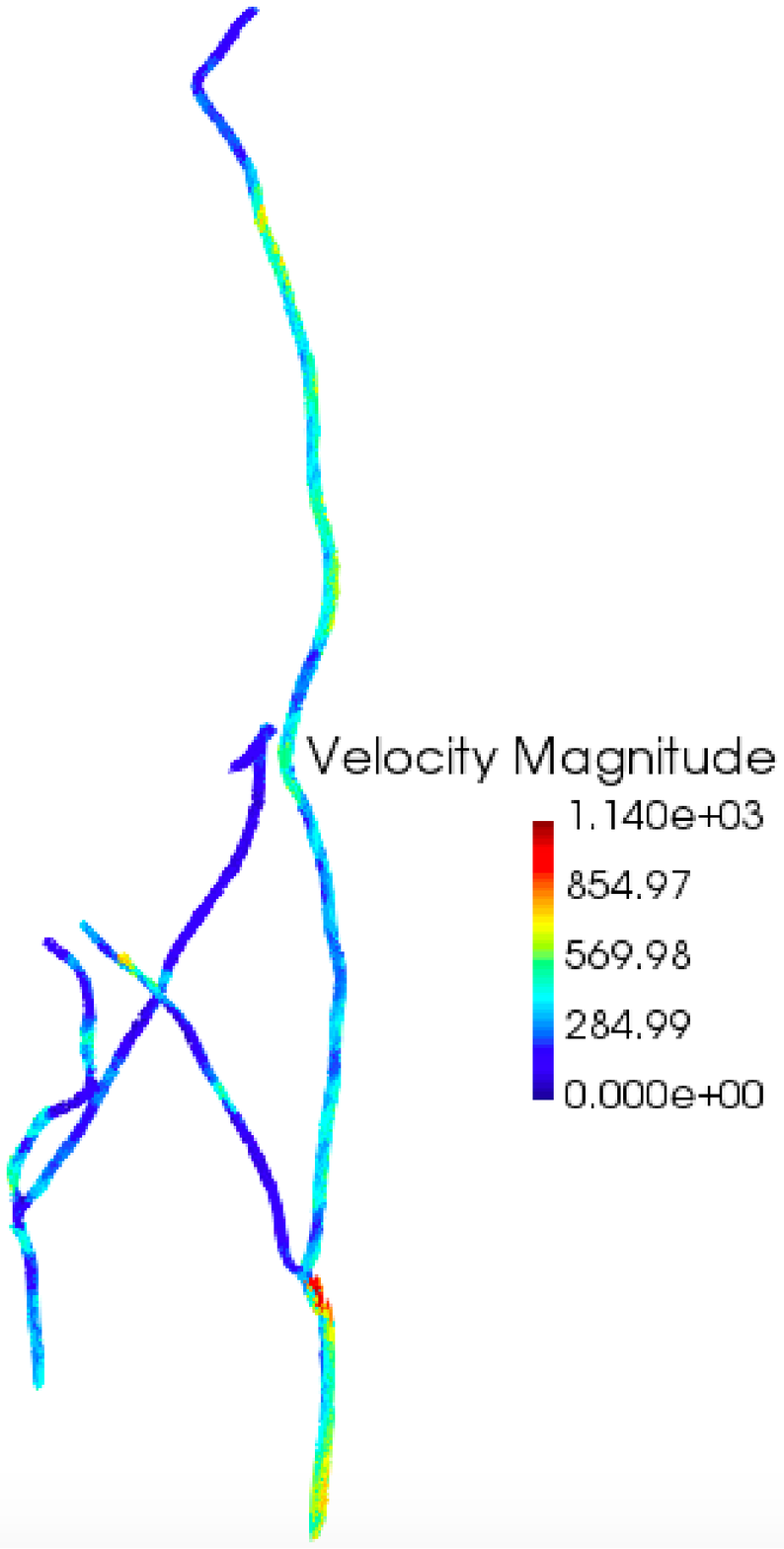}} &  \multicolumn{1}{c}{\multirow{-12}{*}{\rotatebox{90}{(b).}}} & \multicolumn{1}{c|}{\includegraphics[width = 0.2\textwidth, keepaspectratio]{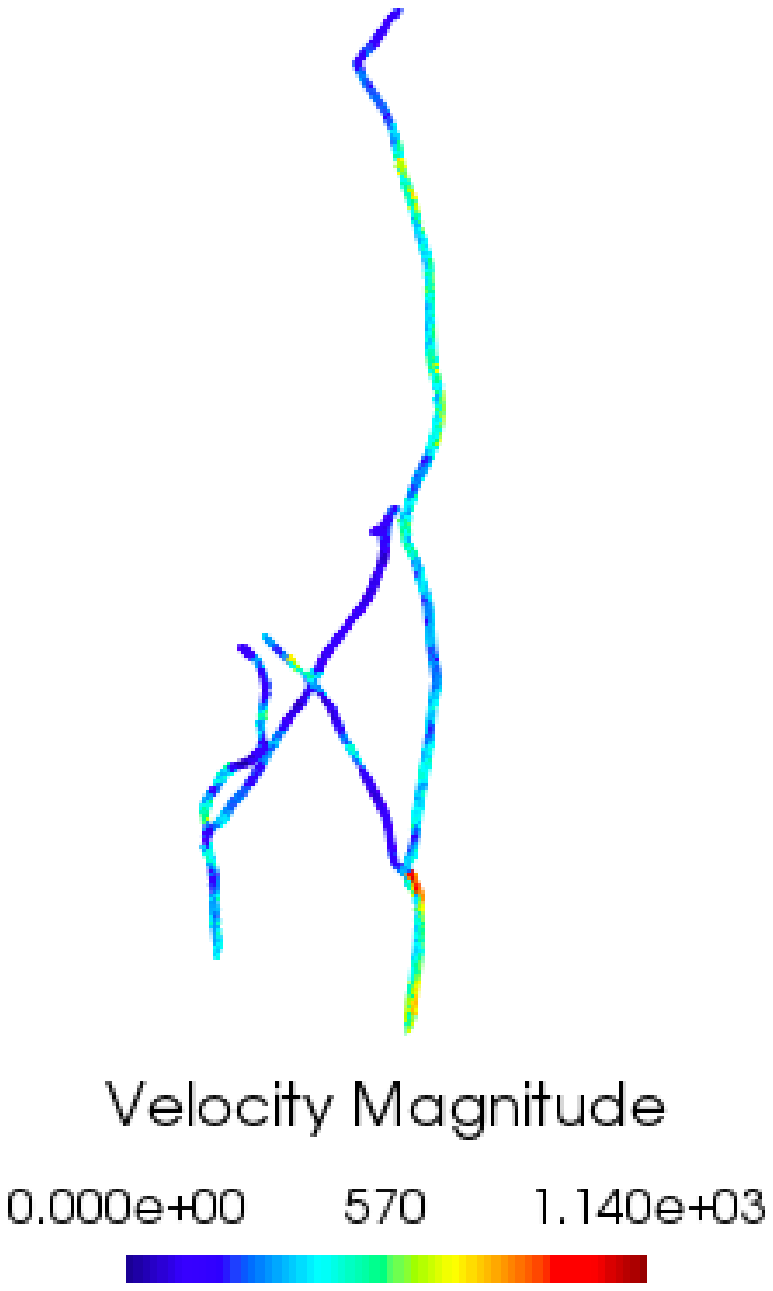}} & \multicolumn{1}{c}{\multirow{-12}{*}{\rotatebox{90}{(c).}}} & \multicolumn{1}{c|}{\rotatebox{90}{\includegraphics[width = 0.35\textwidth, keepaspectratio]{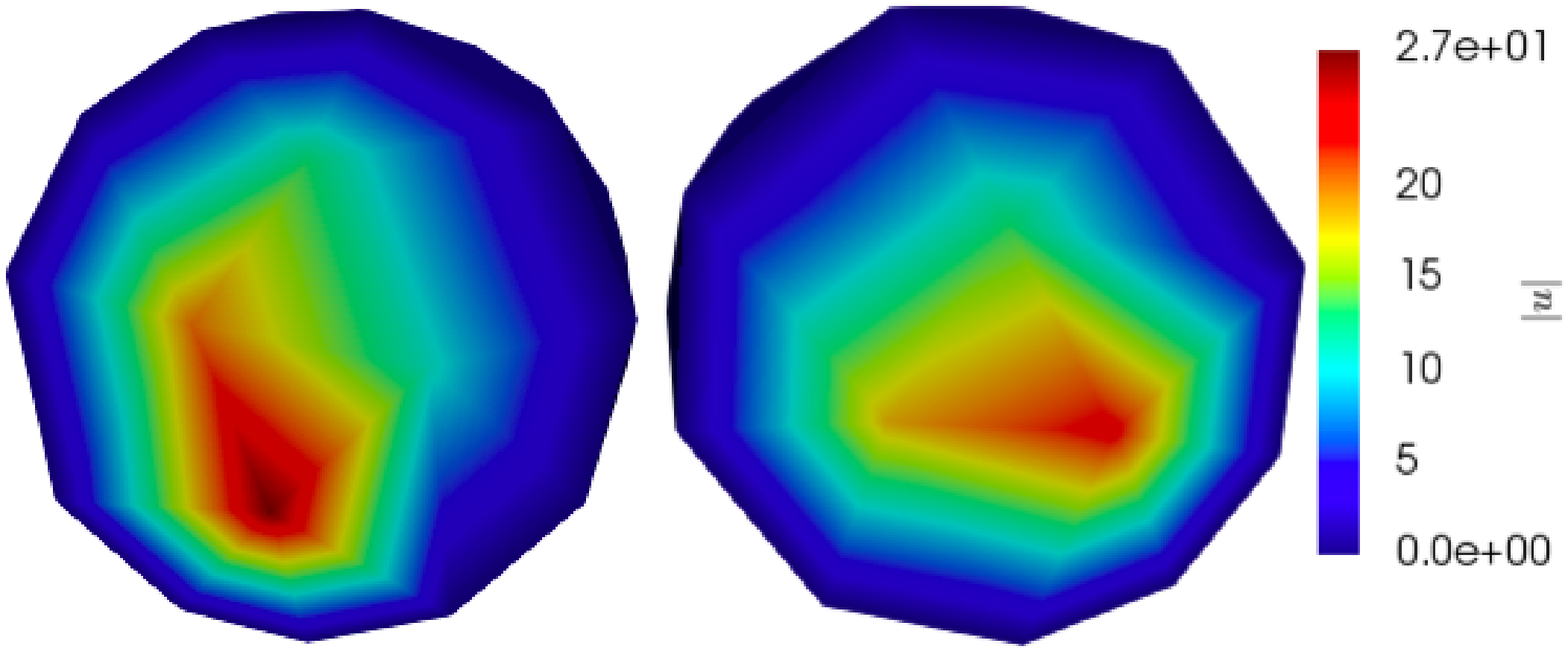}}}\\
\multicolumn{2}{|c|}{} & \multicolumn{2}{c|}{} & \multicolumn{2}{c|}{}\\
\hline
\end{tabular}
\caption{(a). Galerkin finite element approx. of velocity (mm/s). (b). {POD}--Galerkin approx. of velocity $\left( mm/s \right)$. (c). Boundary control magnitude $\left( mm^2/s^2\right)$, top: $\Gamma_{o_2}$, bottom: $\Gamma_{o_1}$.}
\label{RIMAOM1velAndControlApprox.}
\end{figure}

\begin{figure}[h]
\centering
\begin{tabular}{llllll}
\multicolumn{1}{c}{\multirow{-12}{*}{\rotatebox{90}{(a).}}} & \multicolumn{1}{c}{\includegraphics[width = 0.4\textwidth , keepaspectratio]{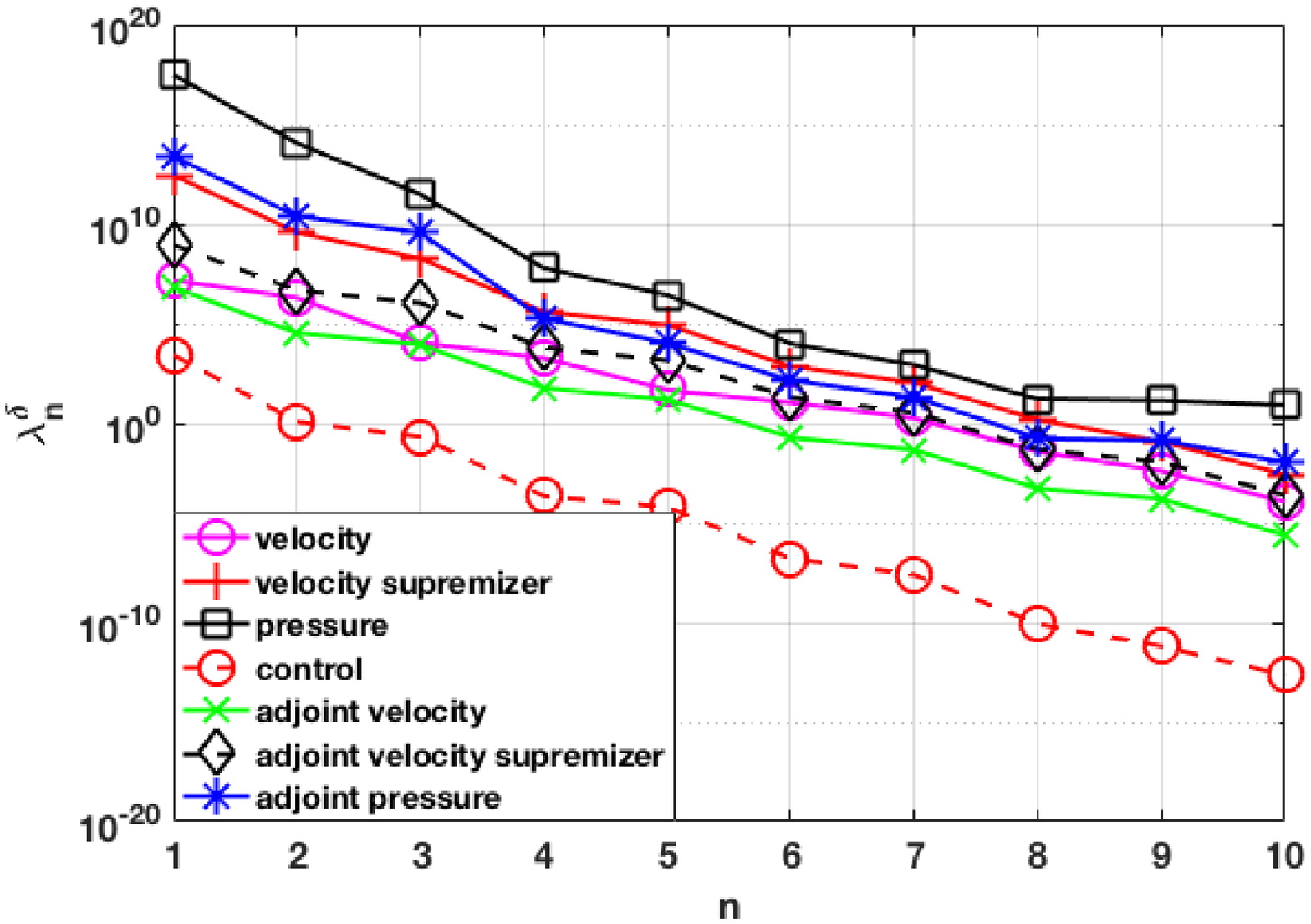}} & \multicolumn{1}{c}{\hspace*{0.5cm}} & \multicolumn{1}{c}{\multirow{-12}{*}{\rotatebox{90}{(b).}}} & \multicolumn{1}{c}{\includegraphics[width = 0.4\textwidth , keepaspectratio]{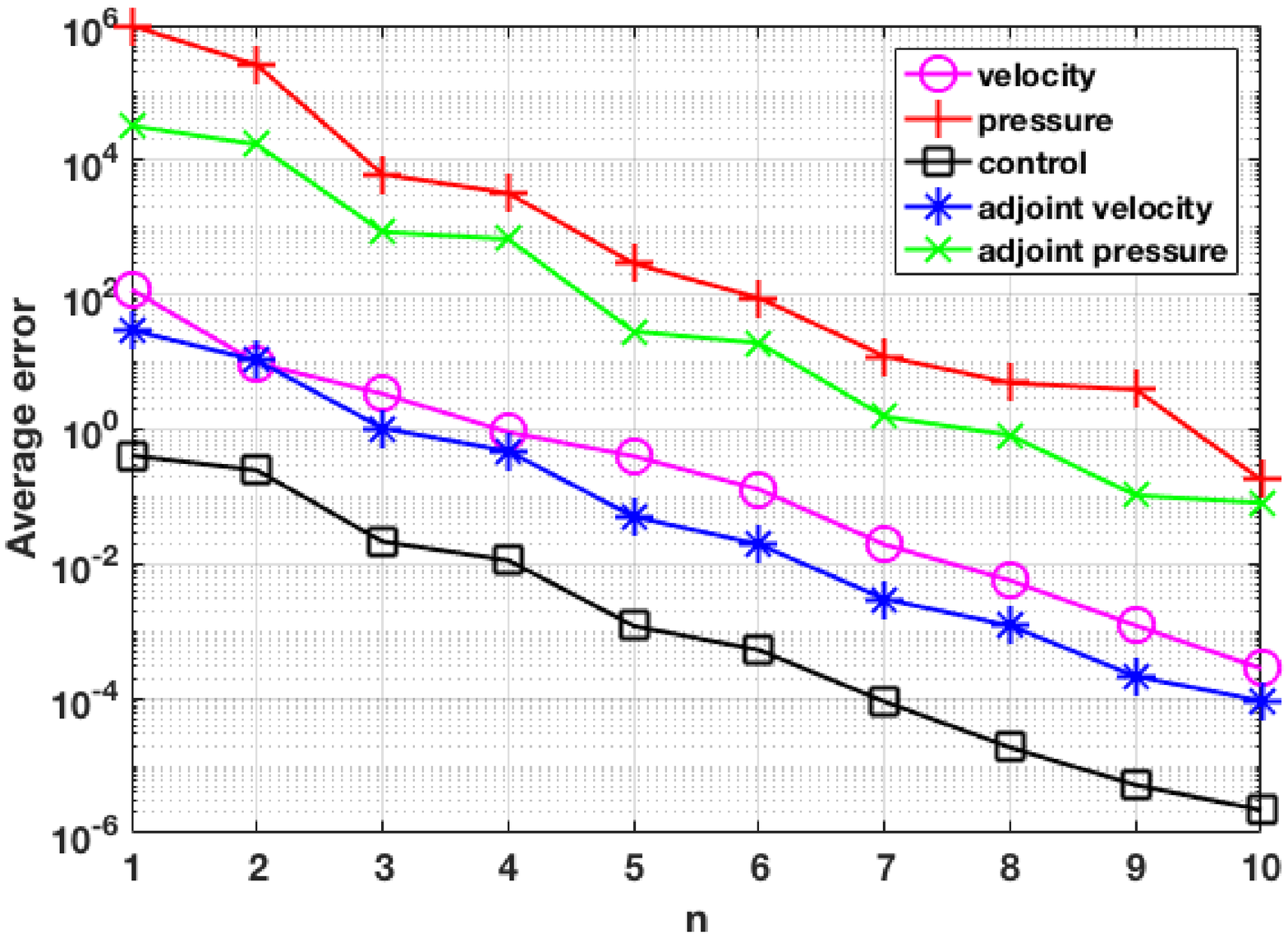}}\\
\multicolumn{1}{c}{\multirow{-12}{*}{\rotatebox{90}{(c).}}} & \multicolumn{1}{c}{\includegraphics[width = 0.4\textwidth , keepaspectratio]{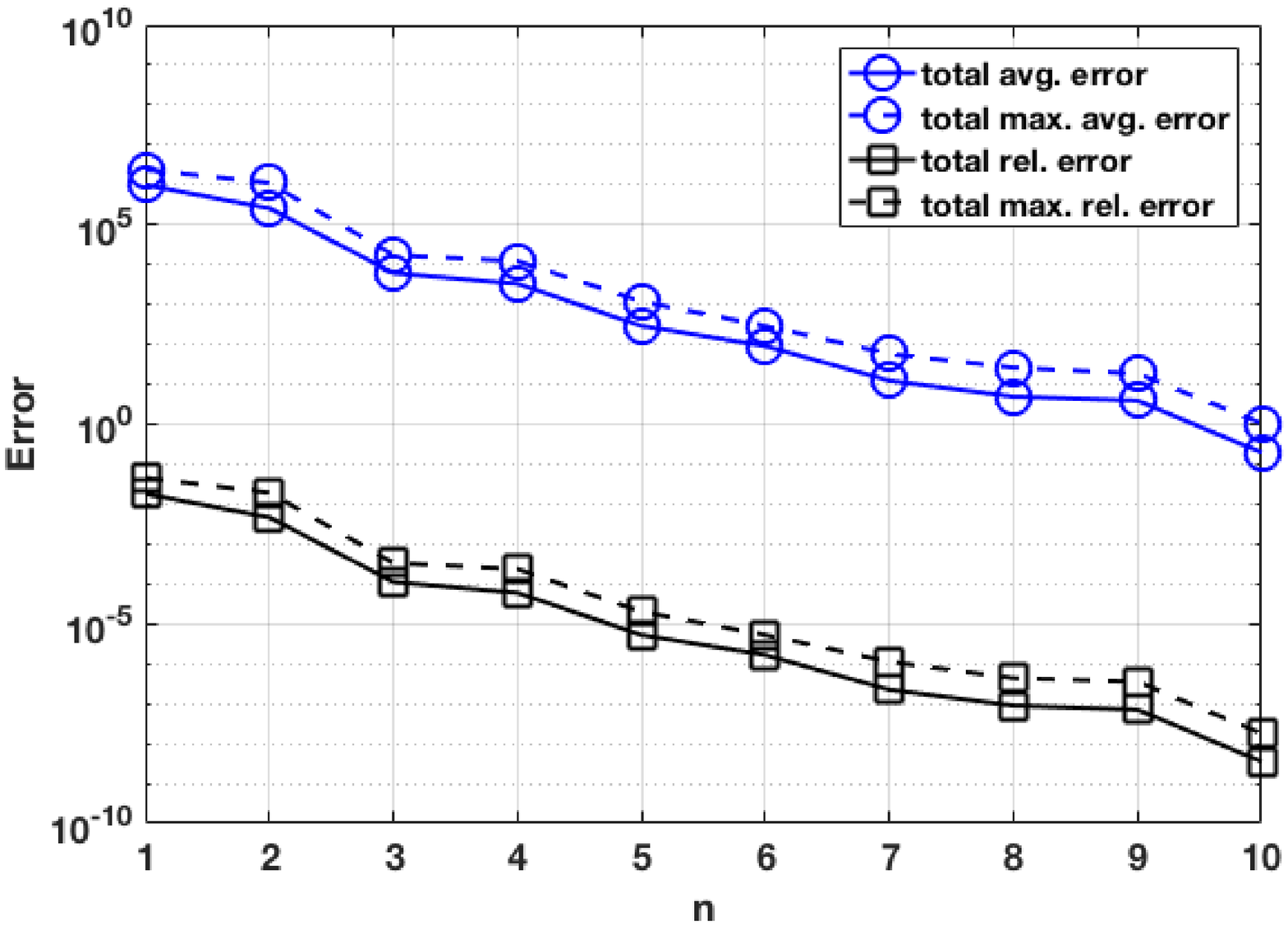}} &  \multicolumn{1}{c}{\hspace*{0.5cm}} & \multicolumn{1}{c}{\multirow{-12}{*}{\rotatebox{90}{(d).}}} &\multicolumn{1}{c}{\includegraphics[width = 0.4\textwidth , keepaspectratio]{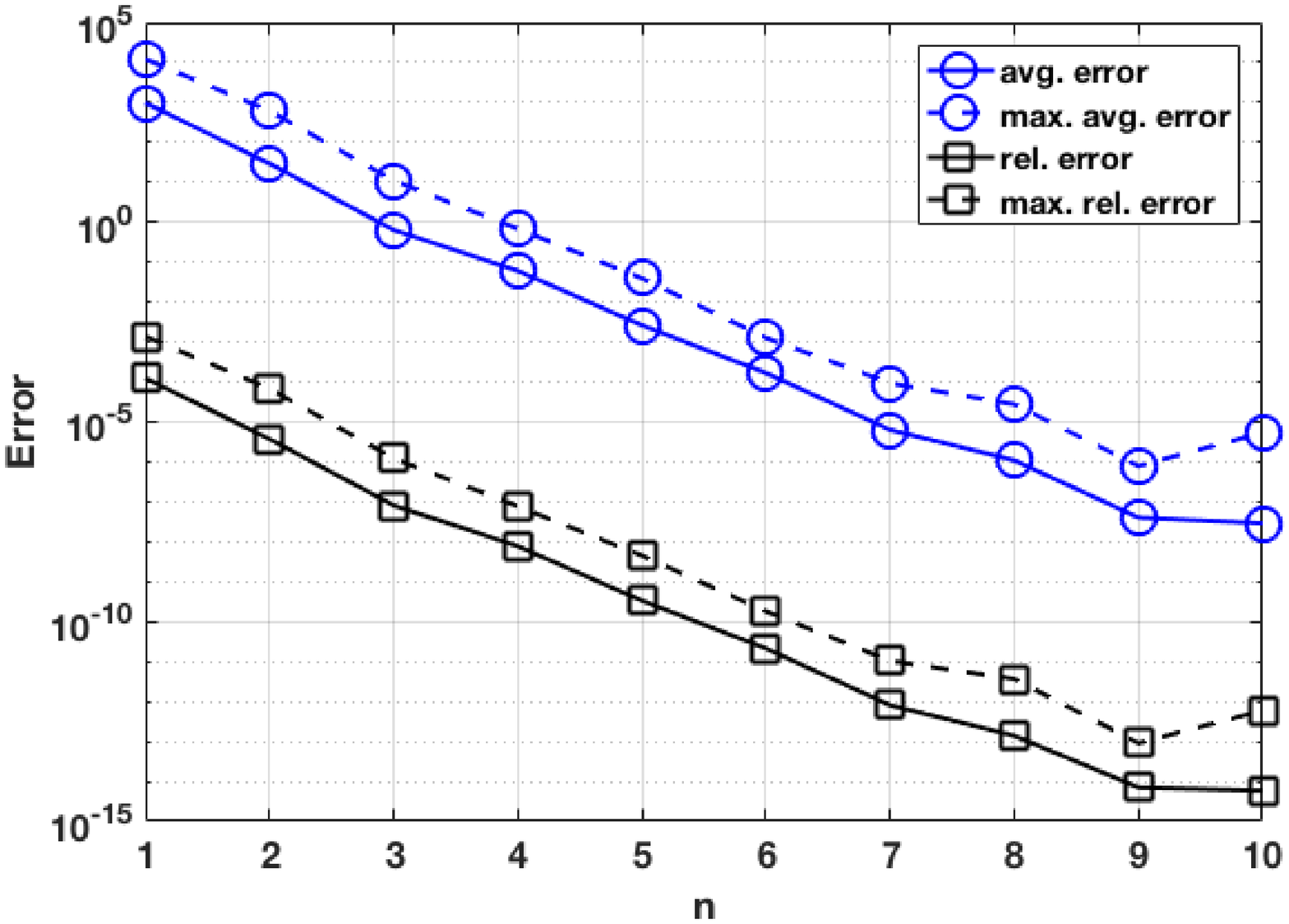}}\\
\end{tabular}
\caption{(a). Eigenvalues' reduction. (b). Average error between FE and {POD}--Galerkin approx. of $\bm{\delta} = \bm{v}, p, \bm{u}, \bm{w}, q$. (c). Total error between Galerkin FE and {POD}--Galerkin approximations. (d). Error between FE and {POD}--Galerkin reduction of $\mathcal{J}$.}
\label{ErrorsRIMAOM1}
\end{figure}

We further apply the reduced order optimal flow control framework (figure \ref{tikzfigure2}) to $\left(\Omega_c\right)$ (figure \ref{figure2}(c)) comprising of two grafts connected to two separate stenosed arteries. Mathematical problem is the same as considered for single graft connections, with the number of inlets and outlets doubled. We label the inlets of right internal mammary artery (RIMA) and left anterior descending artery (LAD) as $\Gamma_{{in}_1}$ and the inlets of saphenuous vein (SV) and first obtuse marginal artery (OM1) as $\Gamma_{{in}_2}$. Similarly respective outlets are marked as $\Gamma_{o_1}$ and $\Gamma_{o_2}$. We consider two physical parameters $\left( \bm{\mu}_1, \bm{\mu}_2\right)$, both being Reynolds numbers, that is, $\left( \bm{\mu}_1, \bm{\mu}_2\right) = \left( Re|_{\Gamma_{{in}_1}}, Re|_{\Gamma_{{in}_2}}\right)\in \left[70, 80\right]\times\left[45, 50\right]$ appearing in velocity profiles being defined at inlets. Moreover, we consider following solution spaces for velocity, pressure and control respectively:
\begin{equation*}
V\left(\Omega_c\right) =  H^1_{\Gamma_{{in}_1} \cup \Gamma_{{{in}_2}}\cup \Gamma_w} = \left[ \bm{v}\in \left[H^1\left(\Omega_c\right)\right]^3 : \bm{v}|_{\Gamma_{{in}_1}} = {\bm{v_{in}}}_1,\ \bm{v}|_{\Gamma_{{in}_2}} = {\bm{v_{in}}}_2\ \text{and}\ \bm{v}|_{\Gamma_w} = \bm{0} \right],\ P\left(\Omega_c\right)= L^2\left(\Omega_c\right),
\end{equation*}
and
\begin{equation*}
U\left(\Gamma_{o_1} \cup \Gamma_{o_2}\right)= \left[ L^2\left(\Gamma_{o_1} \cup \Gamma_{o_2}\right)\right]^3,
\end{equation*}
where, 
\begin{equation*}
{\bm{v_{in}}}_1\left(\bm{\mu}_1\right) = -\frac{\eta\bm{\mu}_1}{R_{{in}_1}}\left(1 - \frac{r^2}{R_{{in}_1}^2}\right){\bm{n}}_{{in}_1},\quad\quad {\bm{v_{in}}}_2\left(\bm{\mu}_2\right) = -\frac{\eta\bm{\mu}_2}{R_{{in}_2}}\left(1 - \frac{r^2}{R_{{in}_2}^2}\right){\bm{n}}_{{in}_2},
\end{equation*}
and $\bm{v_o}\in L^2\left(\Omega_c\right)$. We show Galerkin finite element and {POD}--Galerkin approximations of velocity in figures \ref{RIMAOM1velAndControlApprox.}(a) and \ref{RIMAOM1velAndControlApprox.}(b) respectively. Figure \ref{RIMAOM1velAndControlApprox.}(c) shows control magnitude calculated altogether over the two outlets. The computational performance of the proposed method for this test is reported in table \ref{tableForComputationalPerfomanceTwoGrafts}.\\
For $\left({\mu}_1, {\mu}_2\right) = \left( 80, 50\right)$, the full order problem has dimensions $\mathcal{N} = 715462$ and takes $1848.13$ seconds for one simulation. We use $\mathrm{N}_{max} = 10$ to construct {POD} bases and retain $99.99 \%$ energy of full order solution manifold (see figure \ref{ErrorsRIMAOM1}(a)). We attain $40$ reduced bases for velocity, $2$ for non-homogeneous Dirichlet conditions at the two inlets, $20$ reduced bases for pressure and $10$ for control. $40$ reduced bases are constructed for adjoint velocity and $20$ for adjoint pressure, thus reduced order spaces are constructed using $132$ bases. Moreover, $202.27$ seconds are spent in performing online phase each time with a computational cost of $26881.7$ seconds for offline phase (see table \ref{tableForComputationalPerfomanceTwoGrafts}). 
\begin{table}
\centering
\begin{tabular}{|cc|c|}
\hline\hline
\multicolumn{2}{|c|}{Mesh size} & $605451$\\
\hline
\multicolumn{2}{|c|}{Galerkin finite element dofs $\left(\mathcal{N}\right)$} & $715462$\\
\hline
\multicolumn{2}{|c|}{No. of reduced order bases $\mathrm{N}$} & $132$\\
\hline
\multicolumn{2}{|c|}{$\mathcal{D}$} & $\left[70, 80\right]\times\left[45, 50\right]$\\
\hline
\multicolumn{2}{|c|}{$| \Lambda |$} & $100$\\
\hline
\multicolumn{1}{|c|}{\multirow{3}{*}{Comp. time}} & Galerkin FE & $1848.13$ seconds \\\cline{2-3}
\multicolumn{1}{|c|}{} & offline phase & $26881.7$ seconds \\\cline{2-3}
\multicolumn{1}{|c|}{} &  online phase & $202.27$ seconds\\\hline
\hline

\end{tabular}
\caption{Computational performances for Navier-Stokes constrained optimal flow contol: two graft connections}
\label{tableForComputationalPerfomanceTwoGrafts}
\end{table}

\noindent We report a reduction to approximately $10^{-4}$ from $10^2$ and to $10^{-1}$ from $10^6$ in error for velocity and pressure approximations  (figure \ref{ErrorsRIMAOM1}(b)). $\mathcal{E}_{\bm{u}}$ decreases approximately to $10^{-6}$ and a similar behavior is observed for adjoint velocity and adjoint pressure (figure \ref{ErrorsRIMAOM1}(c)). Total error $\mathcal{E}_T$ reduces from $10^6$ to approximately $10^{-1}$ with a decrease from $10^{-1}$ to approximately $10^{-8}$ in corresponding total relative error (figure \ref{ErrorsRIMAOM1}(c)). Furthermore, the difference between Galerkin FE and {POD}--Galerkin approximations of $\mathcal{J}$ decreases to about $10^{-7}$ for $n = 10$ (figure \ref{ErrorsRIMAOM1}(d)).

\section{Summary and conclusion}
In this article, we focused on applications of reduced order methods in numerical modeling of patient-specific hemodynamics through parameterized optimal flow control problem with three-folds aim, $\left( i\right)$ automated quantification of unknown boundary conditions, $\left( ii\right)$ computational efficiency in terms of CPU time, and $\left( iii\right)$ matching known physiological data. In section \ref{sec: methodology}, we first summarized the procedure followed for extraction of patient-specific cardiovascular geometries from clinically acquired images and then we discussed mathematical details of the optimal control framework incorporating model order reduction techniques. In this work, we particularly made use of Galerkin finite element method at snapshots level and {POD}--Galerkin to reduce the solution manifold based upon the snapshots. In section \ref{sec: numerical results}, we applied the framework to three patient-specific geometries of coronary artery bypass grafts for different parameter-dependent inflow velocity scenarios, separately considering the cases comprising of one graft connection and two grafts connections. Among these applications, we demonstrated order of reduction achieved in objective functional for controlled physical flow and showed speedup attained by {POD}--Galerkin and its computational performance in comparison to the high-fidelity methods. We conclude this article with the following remarks:

\begin{itemize}
\item[ $\left( i\right)$.] In the applications of the reduced order parametrized optimal flow control model discussed in this work, we have considered different hemodynamics scenarios by tuning Reynolds number that generate laminar flow at the inlets. The optimal flow control part of the discussed framework automatically tunes the outflux required by the mathematical model to match the desired data. This is illustrated in the results (see figures \ref{Boundary control magnitudes and eigenvalues' reductions for case I and case II.}(a) case I, \ref{Boundary control magnitudes and eigenvalues' reductions for case I and case II.}(a) case II and \ref{RIMAOM1velAndControlApprox.}(c)).
\item[ $\left( ii\right)$.] We have shown that the full-order methods usually take $\mathcal{O}\left( 10^3\right)$ seconds for a single simulation. In parametrized settings, this time is increased by the number of parameter values taken into consideration, or in other words, considering only the full-order methods in cardiovascular problems, one has to repeatedly spend about $\mathcal{O}\left( 10^3\right)$ seconds. On the other hand, through the implementation of reduced order methods, one only has to repeat the online phase for different tunings of Reynolds number. We have shown that the time taken by this phase is in the range of $\mathcal{O}\left( 10^1\right)$ to $\mathcal{O}\left( 10^2\right)$ seconds, for instance see tables \ref{ComputationalPerformanceAndSpeedupForStokes}(a), \ref{tableForComputationalPerfomanceInCaseIAndCaseII} and  \ref{tableForComputationalPerfomanceTwoGrafts}. Thus, reduced order methods provide an inexpensive computational environment for such coupled problems dependent on different parametrized scenarios, as the ones arising in computational hemodynamics modelling.
\item[ $\left( iii\right)$.] As discussed in single graft connection cases, we reiterate that the computational efficiency, discussed in the previous remark, comes at the cost of performing the expensive offline phase. This is justifiable since the computational cost of this phase depends upon cardinality of sampled training set, the cost of high order solutions and in case of Navier-Stokes state constraints, it also depends upon the cost of reassembly of the non-linear and trilinear operators from scratch and also upon the cost of iterative numerical methods for non-linear PDEs-dependent problems. We remark that the expense made in offline phase is balanced out in biomedical problems since as the problem grows in complexity, for instance, in case of increased number of parameters, geometrical parameterization, denser mesh and larger cardiovascular structures, etc., the computational cost of the full-order methods grows too. Consequently, utilizing only the full-order methods for repetitive simulations leads to the impractical computational cost of order of days. In such cases, use of the inexpensive online phase becomes vital.
\item[ $\left( iv\right)$.] We further remark that we are using third order tensors to save non-linear convection operators for Navier-Stokes state constraints, where affine decomposition assumption does not hold true. In this work, the operator is reassembled at each iteration, impacting the computational cost of offline as well as online phase for non-linear state constraints. However, in the linear case, that is Stokes state constraints, the affine decomposition assumption is preserved. Figure \ref{ComputationalPerformanceAndSpeedupForStokes}(b) shows that maximum of $\mathcal{O}\left( 10^4\right)$ speedup is achieved by {POD}--Galerkin, both for objective functional reduction and solutions to optimal control problem for linear state constraints.
\item[ $\left( v\right)$.] Lastly, we remark that through the application of the proposed reduced order  parametrized optimal flow control model, we have approximated full order solution spaces comprising of $\mathcal{O}\left( 10^5\right)$ degrees of freedom through $\mathcal{O}\left( 10^2\right)$ reduced bases together with sufficient reduction in error between approximations, illustrated in figures \ref{Errors for case I and case II} and \ref{ErrorsRIMAOM1}.

\end{itemize}

\subsection{Future perspectives}
Thanks to {POD}--Galerkin, the computational framework proposed and applied in this work noticeably lowers the high computational cost of numerical approximations of optimal flow control problems, commonly reported by Romorowski et al. \citealp{Veneziani2018}, Tiago et al. \citealp{TiagoEtAl2017} and Koltukluo{\u g}lu et al. \citealp{KoltukluogluEtAl2018}. Furthermore, the application of this framework achieves the similar objectives as in the previously mentioned works, that are, quantification of boundary conditions and data assimilation. Nevertheless, the discussed framework is an initial step in creating a complete pipeline to be exploitable by medical community in clinical studies for coronary artery bypass grafts surgery. To serve the purpose, we propose some future extensions of the work discussed in this article.

A near-future extension of the proposed pipeline is to match patient-specific physiological data, extracted from 4D-flow MRIs. The patient-specific physiological quantities can be blood flow rate or flow velocity and can be imposed through the expression used to impose desired data in $\mathcal{J}$ in this work. Efforts in this regards are work-in-progress \citealp{Francesca}. Furthermore, we propose use of surrogate lumped parameter network models at the outlets to quantify the boundary conditions more accurately. However, as already mentioned, these models require manual tuning of a few parameters until user-desired accuracy is attained. We propose automated quantification of these parameters through the unknown control variables.

Moreover, we suggest to take fluid-structure interaction into account for the discussed framework for accurate patient-specific hemodynamics modeling. In such cases, we also suggest to consider time-dependent optimal flow control problems \citealp{StollEtAl2013,ZainibEtAl2} to model patient-specific hemodynamics with the discussed objectives in blood vessels with sufficiently large diameter such as in the case of aorto-coronary bypass grafts surgery.\\

\noindent {\bf Acknowledgements.}
We acknowledge the support by European Research Council Executive Agency, ERC CoG 2015 AROMA-CFD, GA 681447, 2016-2021. We also appreciate the support of INDAM-GNCS, project on ``Advanced intrusive and non-intrusive model order reduction techniques and applications''.

\bibliographystyle{abbrvnat}
\bibliography{ms}
\end{document}